\documentclass[a4paper,reqno,twoside,11pt]{amsart}

\usepackage{a4wide}
\usepackage{amssymb}
\usepackage{latexsym}
\usepackage{amsmath}
\usepackage{color,soul}

\theoremstyle{plain}
\newtheorem{propn}{\sc Proposition}[section]
\newtheorem{thm}[propn]{\sc Theorem}
\newtheorem{lemma}[propn]{\sc Lemma}
\newtheorem{cor}[propn]{\sc Corollary}

\theoremstyle{definition}

\newtheorem{example}[propn]{\sc Example}

\theoremstyle{remark}
\newtheorem*{rem}{\sc Remark}
\newtheorem*{rems}{\sc Remarks}

\newtheorem*{note}{\sc Note}

\newcommand{\ep}{\varepsilon}
\newcommand{\al}{\alpha}
\newcommand{\be}{\beta}

\newcommand{\Ga}{\Gamma}
\newcommand{\de}{\delta}
\newcommand{\De}{\Delta}
\newcommand{\ka}{\kappa}
\newcommand{\la}{\lambda}
\newcommand{\La}{\Lambda}
\newcommand{\om}{\omega}
\newcommand{\Om}{\Omega}
\newcommand{\ups}{\upsilon}
\newcommand{\Ups}{\Upsilon}
\newcommand{\si}{\sigma}
\newcommand{\Si}{\Sigma}

\newcommand{\Hil}{\mathsf{H}}
\newcommand{\hil}{\mathsf{h}}
\newcommand{\Kil}{\mathsf{K}}
\newcommand{\kil}{\mathsf{k}}
\newcommand{\Til}{\mathsf{T}}
\newcommand{\Al}{\mathsf{A}}
\newcommand{\uAl}{\Al^{\unitise}}
\newcommand{\Cpt}{K}
\newcommand{\Ou}{\mathsf{U}}
\newcommand{\Ov}{\mathsf{V}}
\newcommand{\Ow}{\mathsf{W}}
\newcommand{\Oe}{\mathsf{E}}
\newcommand{\vN}{\mathsf{M}}
\newcommand{\X}{\mathfrak{X}}

\newcommand{\Real}{\mathbb{R}}
\newcommand{\Rplus}{{\Real_+}}
\newcommand{\Comp}{\mathbb{C}}
\newcommand{\Nat}{\mathbb{N}}
\newcommand{\Int}{\mathbb{Z}}

\newcommand{\ip}[3][]{#1\langle #2, #3 #1\rangle}
\newcommand{\norm}[1]{\lVert #1 \rVert}
\newcommand{\cbnorm}[1]{\lVert #1 \rVert_\cb}
\newcommand{\abs}[1]{\lvert #1 \rvert}

\newcommand{\w}[1]{\varpi (#1)}
\newcommand{\bra}[1]{\langle #1 \vert}
\newcommand{\ket}[1]{\vert #1 \rangle}
\newcommand{\dyad}[2]{\ket{#1}\bra{#2}}

\newcommand{\init}{\mathfrak{h}}
\newcommand{\noise}{\mathsf{k}}
\newcommand{\khat}{{\wh{\noise}}}
\newcommand{\etaaug}{{\ol{\eta}}}
\newcommand{\etaaugdec}{{\etaaug \text{\tu{-dec}}}}
\newcommand{\Fock}{\mathcal{F}}
\newcommand{\Exps}{\mathcal{E}}
\newcommand{\Step}{\mathbb{S}}
\newcommand{\BoxJ}{\square_\Jx}
\newcommand{\BoxJt}{\square_{\Jx, t}}
\newcommand{\ind}{\mathbf{1}}
\newcommand{\Mat}{\mathrm{M}}

\newcommand{\Diag}{\mathrm{D}}
\newcommand{\Ix}{\mathcal{I}}
\newcommand{\Jx}{\mathcal{J}}
\newcommand{\Lx}{\mathcal{L}}
\newcommand{\Rx}{\mathcal{R}}
\newcommand{\Sx}{\mathcal{S}}

\newcommand{\rd}{\mathrm{d}}
\newcommand{\srd}{\,\rd}
\newcommand{\I}{\mathrm{i}}
\newcommand{\Lindbladian}{\mathcal{L}}
\newcommand{\Op}{\mathcal{O}}
\newcommand{\domain}{\mathcal{D}}

\newcommand{\bfx}{\mathbf{x}}

\newcommand{\unitise}{{\text{\tu{u}}}}
\newcommand{\cb}{{\text{\tu{cb}}}}
\newcommand{\bd}{{\text{\tu{b}}}}

\newcommand{\kadec}{{\ka \text{\tu{-dec}}}}

\newcommand{\qc}{{\text{\tu{qc}}}}
\newcommand{\cp}{{\text{\tu{cp}}}}
\newcommand{\cpqc}{{\text{\tu{cpqc}}}}
\newcommand{\cpc}{{\text{\tu{cpc}}}}
\newcommand{\cpu}{{\text{\tu{cpu}}}}

\newcommand{\h}{{\text{\tu{h}}}}
\newcommand{\unital}{{\text{\tu{u}}}}
\newcommand{\wmult}{{\text{\tu{mult}}}}
\newcommand{\starhom}{{{}^*\text{\tu{-hom}}}}

\newcommand{\wh}{\widehat}
\newcommand{\wt}{\widetilde}
\newcommand{\ol}{\overline}
\newcommand{\ot}{\otimes}
\newcommand{\aot}{\underline{\ot}\,}
\newcommand{\vot}{\overline{\ot}\,}
\newcommand{\otm}{\ot_\Mat}

\newcommand{\op}{\oplus}

\newcommand{\schur}{\boldsymbol{\cdot}}
\newcommand{\les}{\leqslant}
\newcommand{\ges}{\geqslant}
\newcommand{\To}{\rightarrow}
\newcommand{\Tends}{\rightarrow}
\newcommand{\Implies}{\Rightarrow}
\newcommand{\Iff}{\Leftrightarrow}
\newcommand{\fsubset}{\subset\subset}

\newcommand{\ti}{\textit}
\newcommand{\tu}{\textup}

\newcommand{\El}{El\text{-}}
\DeclareMathOperator{\Alg}{Alg}
\DeclareMathOperator{\Dom}{Dom}
\DeclareMathOperator{\Ran}{Ran}
\DeclareMathOperator{\Lin}{Lin}

\DeclareMathOperator{\id}{id}
\DeclareMathOperator{\re}{Re}
\DeclareMathOperator{\im}{Im}

\DeclareMathOperator{\inter}{inter}

\newcommand{\msgp}{\mathcal{P}}
\newcommand{\msgq}{\mathcal{Q}}
\newcommand{\zmsgp}{{}^0\!\msgp}
\newcommand{\zk}{{}^0\!k}
\newcommand{\bek}{{}^\be\! k}
\newcommand{\osgp}{P}

\newcommand{\gencp}{\mathfrak{cp}}
\newcommand{\gencpqc}{\mathfrak{cpqc}}
\newcommand{\gencpbeta}{\mathfrak{cpqc}_\be}
\newcommand{\gencpzero}{\mathfrak{cpqc}_0}
\newcommand{\gencpc}{\mathfrak{cpc}}
\newcommand{\gencpu}{\mathfrak{cpu}}
\newcommand{\genh}{\mathfrak{real}}
\newcommand{\genu}{\mathfrak{u}}
\newcommand{\genwmult}{\mathfrak{mult}}
\newcommand{\genstarhom}{{}^*\text{-}\mathfrak{hom}}
\newcommand{\genflow}{\mathfrak{flow}}
\newcommand{\genstruct}{\mathfrak{struct}}
\newcommand{\genstarhomone}{{}^*\text{-}\mathfrak{hom}_1}

\newcommand{\QSC}{\text{\tu{QSC}}}
\newcommand{\QSCbeta}{\QSC_\qc^\be}
\newcommand{\QSCqc}{\QSC_\qc}
\newcommand{\QSCcb}{\QSC_\cb}
\newcommand{\QSCcp}{\QSC_\cp}
\newcommand{\QSCcpc}{\QSC_\cpc}
\newcommand{\QSCcpu}{\QSC_\cpu}
\newcommand{\QSCcpqc}{\QSC_\cpqc}
\newcommand{\QSCcpbeta}{\QSC_\cpqc^\be}
\newcommand{\QSCb}{\QSC_\bd}
\newcommand{\QSCh}{\QSC_\h}
\newcommand{\QSCu}{\QSC_\unital}

\newcommand{\QSF}{\text{\tu{QSF}}}

\newenvironment{alist}
{

\begin{enumerate}}
{\end{enumerate}}

\newenvironment{rlist}
{

\begin{enumerate}}
{\end{enumerate}}

\numberwithin{equation}{section}
\begin{document}

\title[Quantum stochastic semigroups on operator spaces II]
{Quantum stochastic cocycles \\
and completely bounded semigroups \\
on operator spaces II}
\dedicatory{Dedicated to the memory of Ola Bratteli}
\author[Martin Lindsay]{J.\ Martin Lindsay}
\address{Department of Mathematics and Statistics \\ Lancaster
University \\ Lancaster LA1 4YF \\ UK}
\email{j.m.lindsay@lancs.ac.uk}
\author[Stephen Wills]{Stephen J.\ Wills}
\address{School of Mathematical Sciences \\ University College
Cork \\ Cork \\ Ireland}
\email{s.wills@ucc.ie}

\subjclass[2000]{
Primary
81S25,          
46L53;          
Secondary
46N50,          
47D06}          

\keywords{
Quantum dynamical semigroup,
quantum Markov semigroup,
  completely positive,
  quasicontractive,
  generator,
  operator space,
  operator system,
  matrix space,
  Schur-action,
  Markovian cocycle,
  stochastic semigroup,
  quantum exclusion process}

\begin{abstract}
Quantum stochastic cocycles provide a basic model for time-homogeneous
Markovian evolutions in a quantum setting, and a direct counterpart in
continuous time to quantum random walks, in both the Schr\"{o}dinger
and Heisenberg pictures. This paper is a sequel to one in which
correspondences were established between classes of quantum stochastic
cocycle on an operator space or $C^*$-algebra, and classes of
Schur-action `global' semigroup on associated matrix spaces over the
operator space. In this paper we investigate the stochastic generation
of cocycles via the generation of their corresponding global
semigroups, with the primary purpose of strengthening the scope of
applicability of semigroup theory to the analysis and construction of
quantum stochastic cocycles. An explicit description is given of the
affine relationship between the stochastic generator of a completely
bounded cocycle and the generator of any one of its associated global
semigroups. Using this, the structure of the stochastic generator of a
completely positive quasicontractive quantum stochastic cocycle on a
$C^*$-algebra whose expectation semigroup is norm continuous is
derived, giving a comprehensive stochastic generalisation of the
Christensen--Evans extension of the GKS\&L theorem of Gorini,
Kossakowski and Sudarshan, and Lindblad. The transformation also
provides a new existence theorem for cocycles with unbounded structure
map as stochastic generator. The latter is applied to a model of
interacting particles known as the quantum exclusion Markov process,
in particular on integer lattices in dimensions one and two.
\end{abstract}

	\maketitle

\vspace*{-15cm}
\begin{flushright}
\ti{\small To appear in:
\\
Communications in Mathematical Physics}
\end{flushright}
\vspace{15cm}

 \tableofcontents


\section*{Introduction}

This paper is a continuation of our operator space analysis of quantum
stochastic cocycles begun in~\cite{spawn}. In that paper we
established correspondences between classes of completely bounded
quantum stochastic cocycle on an operator space and associated
`global' semigroups of Schur-action maps on corresponding matrix
spaces. This was done for a variety of types of operator space,
including operator systems, $C^*$-algebras and column operator spaces,
and corresponding cocycle types, such as completely contractive,
completely positive and contractive operator cocycles. In the current
paper we analyse the relationship between the stochastic generation of
such cocycles via quantum stochastic differential equations, and the
generation of corresponding global semigroups, and identify the affine
transformation from one to the other. This is used to obtain a new
existence theorem for quantum stochastic differential equations with
unbounded coefficients. We also use it to characterise, directly and
in full generality, the stochastic generators of completely positive
quasicontractive quantum stochastic cocycles on a $C^*$-algebra whose
expectation semigroups are norm continuous. In view of the fact that
norm-continuous semigroups are automatically quasicontractive, this
result amounts to a full quantum stochastic extension of the
well-known characterisation of the generators of norm-continuous,
completely positive semigroups on a $C^*$-algebra (\cite{ChrisEvans}).
The new existence theorem for quantum stochastic cocycles is shown to
be applicable to `structure maps' on a $C^*$-algebra whose domain is
matricially square-root closed, by means of semigroup theory. The
structure relations are necessary conditions for any generated quantum
stochastic cocycle to be a quantum stochastic flow, that is to be
unital and ${}^*$-homomorphic. In turn, we construct such structure
maps to obtain quantum Markov exclusion processes on integer lattices
governed by quantum stochastic differential equations by applying a
theorem of Bratteli and Kishimoto.

In brief, the main results of the paper are as follows.
Theorem~\ref{thm: KGagen} gives the affine transformation from
stochastic generator to generator of associated global semigroup, for
elementary quantum stochastic cocycles with completely bounded
stochastic generator (\emph{elementary}, formerly \emph{cb-Markov
  regular}, means that each of its associated semigroups is cb-norm
continuous).  Theorem~\ref{thm: ABC} is a new existence theorem for
the generation of quantum stochastic cocycles. Theorem~\ref{thm: nonu}
characterises the stochastic generator of completely positive
quasicontractive elementary quantum stochastic cocycles.
Theorem~\ref{thm: stoch gen} applies the above-mentioned existence
theorem to structure maps, and Theorem~\ref{thm: 7.4} demonstrates its
applicability to quantum exclusion processes.

The prequel of this paper contains a long motivational
introduction, and extensive bibliography. A detailed outline of the
contents of the current paper follows.

In Section~\ref{MS+decomp} we develop some theory of
$\ka$-decomposable completely bounded maps between $\hil$-matrix
spaces over operator spaces (whose definition is recalled there), for
a Hilbert space $\hil$ with orthonormal basis $\ka$. Under the
identification of the matrix space with an actual space of matrices
with operator entries, determined by the basis, these
$\ka$-decomposable maps correspond to Schur-action maps. Also, each
$\ka$-decomposable completely bounded map restricts to a map between
spatial tensor products of the space of compact operators on $\hil$
with the corresponding operator space. We show that this restriction
is completely isometrically reversible (Proposition~\ref{new
  propn}). The theory is then applied to obtain a characterisation of
the generators of norm-continuous completely positive semigroups of
$\ka$-decomposable operators on the $\hil$-matrix space over a unital
$C^*$-algebra (Theorem~\ref{decomp CE}). It is worth recalling that
such matrix spaces are operator systems, but typically \emph{not}
$C^*$-algebras.

In Section~\ref{cocycle cty} we investigate in detail the
correspondence between continuity properties, in the time parameter
$t$, of completely bounded quantum stochastic cocycles and that of
their global semigroups. In Section~\ref{Section: QSC QSDE} we relate
quantum stochastic cocycles to a basic source, namely quantum
stochastic differential equations. In Section~\ref{Section: gen}, by
applying quantum stochastic calculus to the appropriate `diagonal Weyl
process' introduced in~\cite{spawn}, we compute the affine
transformation from the stochastic generator of a completely bounded
quantum stochastic cocycle to the generator of the corresponding
global semigroup, and give sufficient conditions for its injectivity
(Theorem~\ref{thm: KGagen}). Restricting to the global semigroups with
respect to an orthonormal basis of the noise dimension space, the
transformation is necessarily injective; we give the form of the
(partially defined) inverse (Theorem~\ref{Z}).  In Section~\ref{sec:
Mr CPC}, bijectivity is demonstrated for specific classes of cocycle.

Section~\ref{CPC cocs} concentrates on completely positive quantum
stochastic cocycles on operator systems and $C^*$-algebras. We first
establish a characterisation of contractivity for completely positive
quantum stochastic cocycles in terms of the generators of its
associated semigroups (Proposition~\ref{CPC+char}). This constitutes
an infinitesimal counterpart to one of the central results
of~\cite{spawn}. Using it we obtain an existence theorem for
completely positive quantum stochastic cocycles governed by a quantum
stochastic differential equation with unbounded coefficient
(Theorem~\ref{thm: ABC}). In Section~\ref{sec: Mr CPC},
Proposition~\ref{CPC+char} and the affine correspondence established
in Section~\ref{Section: gen} are employed to obtain a direct proof of
the full quantum stochastic extension of the characterisation theorem
for generators of norm-continuous quantum dynamical semigroups on a
$C^*$-algebra $\Al$ (\cite{ChrisEvans}, \cite{Lindblad}, \cite{GKS}),
namely a characterisation of the stochastic generators of completely
positive quasicontractive elementary quantum stochastic cocycles on
$\Al$ with arbitrary noise dimension space, in both the unital and
nonunital cases (Theorem~\ref{thm: nonu}), completing our earlier
results~\cite{father} and~\cite{gran}. As an application, all
${}^*$-homomorphic elementary quantum stochastic cocycles are shown to
be stochastically generated (Theorem~\ref{thm: 6.6}).

In Section~\ref{sec: structure} we prove an existence theorem for
completely positive unital quantum stochastic cocycles on a unital
$C^*$-algebra (Theorem~\ref{thm: stoch gen}). Our hypotheses include
necessary conditions on the coefficient map of the governing quantum
stochastic differential equation for the solution to be a quantum
stochastic flow, by which we mean a quantum stochastic cocycle which
is unital and ${}^*$-homomorphic. As a consequence, the existence
theorem dovetails very nicely with recent work on such flows
(\cite{DGS}) in which our conclusion coincides with their standing
hypothesis. A corollary of the theorem is that completely bounded
structure maps generate quantum stochastic cocycles that are
necessarily completely positive and contractive.

In the final section we demonstrate the applicability of the existence
theorem by constructing quantum exclusion processes as dilations of
the quantum Markov semigroups introduced in~\cite{Reb}; the methods
used in~\cite{BeW} and~\cite{BWFK} differ from, but are nicely
complementary to, those employed here.

\subsection{Notation}

The symbols $\aot$, $\ot$ and $\vot$ are used to denote respectively
the algebraic, spatial and ultraweak tensor products. The
$C^*$-algebra of compact operators on a Hilbert space $\hil$ is
denoted $\Cpt(\hil)$. For vectors $\xi \in \hil$ and $\xi' \in \hil'$,
using Dirac-inspired notation of bras and kets, we set
\[
E_\xi := I_\Hil \ot \ket{\xi} \in B(\Hil; \Hil \ot \hil), \
p(\xi) := I_\Hil \ot \dyad{\xi}{\xi } \in B(\Hil \ot \hil)
\text{ and }
E^{\xi'} := I_{\Hil'} \ot \bra{\xi} \in B(\Hil' \ot \hil'; \Hil'),
\]
where the Hilbert spaces $\Hil$ and $\Hil'$ should always be clear
from context. Thus $E^\xi = (E_\xi)^*$, $p(\xi) = E_\xi E^\xi$, and,
for $T \in B(\Hil \ot \hil; \Hil' \ot \hil')$, $E^{\xi'} T E_\xi =
\big( \id_{B( \Hil; \Hil' )} \vot \omega_{\xi', \xi} \big) (T)$ for
the vector functional $\omega_{\xi', \xi}: B( \hil; \hil' ) \to \Comp$,
$A \mapsto \ip{ \xi' }{ A \xi }$. The indicator function of a set $A$
is denoted $\ind_A$. This notation is extended to vector-valued
functions $f:\Rplus \To V$ and subintervals $I$ of $\Rplus$ as
follows: $f_I$ denotes the function from $\Rplus$ to $V$ which agrees
with $f$ on $I$ and is zero outside $I$. In particular we apply this
to vectors in $V$, by viewing them as constant functions from $\Rplus$
to $V$. Thus, for $v \in V$ and $t \in \Rplus$, $v_{[0,t[}$ is the
function equal to $v$ on $[0,t[$ and zero on $[t, \infty[$. We use
the symbol $\fsubset$ to denote subset of finite cardinality.

\subsection{Terminology}

Quantum stochastic analysis is the analysis of noncommutative
processes adapted to the intrinsic operator filtration of a symmetric
Fock space over $L^2( \Rplus; \noise )$ for a (noise, or multiplicity)
Hilbert space $\noise$ (\cite{Lgreifswald}). In particular, quantum
stochastic differential equations are with respect to Fock space
creation, preservation and annihilation processes in the sense of
Hudson and Parthasarathy (\cite{HuP},~\cite{Partha}). This class of
filtrations is universal for tensor independence in various senses
(see \emph{e.g.}~\cite{SSS}). Thus \emph{quantum stochastic cocycle}
is a shorthand for Markovian cocycle, in the sense of Accardi
(\cite{Accardi},~\cite{AFL}), specific to the context of Fock
filtrations. We stress that `classical' processes are not excluded.
Indeed these arise naturally when the cocycle is governed by a quantum
stochastic differential equation whose coefficient map takes a
particular form; this is nicely illustrated in~\cite{KuM}.

\section{Matrix spaces and $\ka$-decomposability}
\label{MS+decomp}

In this section we recall some of the results concerning matrix spaces
developed in~\cite{spawn} and further extend these. Throughout we work
with \emph{concrete} operator spaces, that is, norm-closed subspaces
of $B(\Hil;\Hil')$ for some Hilbert spaces $\Hil$ and $\Hil'$. The
exceptions are Propositions~\ref{new propn} and~\ref{aaa} which
concern operator spaces of completely bounded maps. Excellent
references for operator spaces, operator systems and completely
bounded maps are~\cite{EfR},~\cite{Paulsen} and~\cite{Pisier}; see
also~\cite{BlL} for operator algebras from this perspective. For more
on matrix spaces, see~\cite{Lgreifswald}.

\emph{For the rest of the paper} $\Ov$ and $\Ow$ denote generic
(concrete) operator spaces.

\subsection{Spaces of matrices}

For an index set $\Ix$,
\[
\Mat_\Ix (\Ov)_\bd := \Big\{ A \in \Mat_\Ix (\Ov): \sup_{\Jx \fsubset
\Ix} \norm{A^{[\Jx]}} < \infty \Big\},
\]
where $\Mat_\Ix (\Ov)$ denotes the vector space of matrices
$A = [a^i_j]_{i,j \in \Ix}$ with entries in $\Ov$, and, for
$\Jx \fsubset \Ix$, $A^{[\Jx]} \in \Mat_\Jx (\Ov)$ denotes the
$\Jx \times \Jx$-truncation of $A$ to a finite matrix. This has a
natural operator space structure (\cite{EfR}, Chapter 10). When $\Ov$
is an operator algebra the \emph{Schur product} on $\Mat_\Ix (\Ov)$
\[
[a^i_j] \schur [b^i_j] := [a^i_j b^i_j]
\]
restricts to $\Mat_\Ix (\Ov)_\bd$. Let $\Diag_\Ix (\Ov)$ and
$\Diag_\Ix (\Ov)_\bd$ denote the corresponding subspaces of diagonal
matrices, so that
\[
\Diag_\Ix (\Ov)_\bd := \Big\{ [a^i_j] \in \Mat_\Ix (\Ov): a^i_j = 0
\text{ if } i \neq j \text{ and } \sup_{i \in \Ix} \norm{a^i_i} <
\infty \Big\}.
\]

\subsection{Matrix spaces}

If $B(\Hil; \Hil')$ is the ambient full operator space of $\Ov$ then,
for any Hilbert spaces $\hil$ and $\hil'$, we call the operator space
\[
\Ov \otm B(\hil; \hil') := \big\{ T \in B(\Hil \ot \hil; \Hil' \ot \hil'):
E^x T E_y \in \Ov \text{ for all } x \in \hil', y \in \hil \big\}
\]
the $\hil$-$\hil'$-\emph{matrix space} over $\Ov$ (or
$\hil$-\emph{matrix space} when $\hil' = \hil$). Note that it suffices
to verify the membership condition for $x$ and $y$ running through
total subsets of $\hil'$ and $\hil$ such as orthonormal bases. Thus
$\Ov \otm B(\hil)$ contains the operator space $\Ov \ot \Cpt (\hil)$,
and is an operator system when $\Ov$ is. Note, however, that it is
typically \emph{not} a $C^*$-algebra when $\Ov$ is, unless $\hil$ is
finite dimensional. To any completely bounded operator
$\varphi: \Ov \To \Ow$, and Hilbert space $\hil$, the map
$\varphi \aot \id_{B (\hil)}$ extends uniquely to a bounded operator
$\varphi \otm \id_{B(\hil)}: \Ov \otm B(\hil) \To \Ow \otm B(\hil)$,
satisfying
\[
E^x (\varphi \otm \id_{B(\hil)}(T)) E_y = \varphi (E^x T E_y), \quad T
\in \Ov \otm B(\hil), \ x,y \in \hil.
\]
This is consistent with two well-known extensions, namely
$\varphi \ot \id_{\Cpt (\hil)}$ and, when $\varphi$ is an ultraweakly
continuous completely bounded map between dual operator spaces, also
$\varphi \vot \id_{B(\hil)}$ (\cite{EfR}). These $\hil$-\emph{matrix
  liftings} are completely bounded and satisfy
$\cbnorm{\varphi \otm \id_{B(\hil)}} = \cbnorm{\varphi}$. If $\Ow$ is
explicitly of the form $\Ou \otm B(\Kil; \Kil')$, for some operator
space $\Ou$ and Hilbert spaces $\Kil$ and $\Kil'$, we write
\begin{equation}
\label{flipped lift}
\varphi^\hil \ \text{ for } \ \Pi \circ (\varphi \otm
\id_{B(\hil)}): \Ov \otm B(\hil) \To \Ou \otm B(\hil \ot \Kil; \hil
\ot \Kil'),
\end{equation}
$\Pi$ being the tensor flip $\Ou \otm B(\Kil \ot \hil; \Kil' \ot
\hil) \To \Ou \otm B(\hil \ot \Kil; \hil \ot \Kil')$. If $B(\Hil;
\Hil')$ is the ambient full operator space of $\Ov$ then, for any
ultraweakly continuous completely bounded map $\chi: B(\Kil; \Kil')
\To B(\hil; \hil')$, the map $\id_{B(\Hil; \Hil')} \vot \chi$
restricts to a map
\[
\id_\Ov \otm \, \chi: \Ov \otm B(\Kil; \Kil') \To \Ov \otm B(\hil;
\hil').
\]

Matrix spaces are thus an abstraction of the above operator spaces of
matrices. Specifically, any choice of orthonormal basis
$\ka = (e_i)_{i \in \Ix}$ for $\hil$ determines the completely
isometric isomorphism (\emph{i.e.}~linear isomorphism each of whose
matrix liftings is an isometry)
\begin{equation}
\label{op matrix}
\Ov \otm B(\hil) \cong \Mat_\Ix (\Ov)_\bd, \quad T \longleftrightarrow
{}^\ka T = [T^i_j]_{i,j \in \Ix}, \text{ where } T^i_j := E^{e_i} T
E_{e_j} \in \Ov.
\end{equation}
Let $\Ov \otm \Diag_\ka (\hil)$ denote the subspace of the matrix
space tensor product that corresponds to $\Diag_\Ix (\Ov)_\bd$ under
the above identification. In contrast to $\Ov \otm B(\hil)$, this
operator space is manifestly $\ka$-dependent. Note that
\[
\Ov \otm \Diag_\ka (\hil)
=
\bigl\{
T \in \Ov \otm B(\hil):
T p(e_j) = p(e_j) T \text{ for all } j \in \Ix
\bigr\}.
\]
We have frequent recourse to ampliations of the following kind, for an
operator space $\Ov$ and Hilbert space $\hil$:
\begin{equation}
\label{eqn: iota}
\iota^\Ov_\hil : \Ov \To \Ov \otm B( \hil ), \quad a \mapsto a \ot
I_\hil.
\end{equation}

\subsection{Schur isometries and products}

For a Hilbert space $\hil$ with orthonormal basis $\ka= (e_i)_{i \in
\Ix}$, let $S_\ka: \hil \To \hil \ot \hil$ denote the \emph{Schur
isometry} defined by continuous linear extension of the map $e_i
\mapsto e_i \ot e_i$, and let $\Si_\ka: B(\hil) \To B(\hil \ot \hil)$
be the corresponding \emph{Schur homomorphism}, $T \mapsto S_\ka T
S^*_\ka$. Thus $\Si_\ka$ is an injective normal ${}^*$-homomorphism;
its natural left inverse is the normal completely positive unital map
$\Ups_\ka: R \mapsto S^*_\ka R S_\ka$. Note that
\[
\Si_\ka ( \Cpt (\hil) ) \subset \Cpt ( \hil \ot \hil ) \ \text{ and }
\ \Ups_\ka ( \Cpt ( \hil \ot \hil ) ) = \Cpt ( \hil ).
\]
\begin{lemma}
\label{Ranges}
The following hold\tu{:}
\begin{align}
\big(\id_\Ov \otm\, \Si_\ka\big) \big( \Ov \otm B(\hil) \big) &\subset
\big( \Ov \ot \Cpt (\hil) \big) \otm B(\hil), \text{ and
} \label{Aa} \\
\Ov \otm B(\hil) &= \big(\id_\Ov \otm \Ups_\ka\big) \big( (\Ov \ot
\Cpt (\hil)) \otm B(\hil) \big). \label{Ac}
\end{align}
\end{lemma}

\begin{proof}
Let $T \in \Ov \otm B(\hil)$ and $i, j, l \in \Ix$, and let
$B(\Kil; \Kil')$ be the ambient full operator space of $\Ov$. The
identity
\[
(I_\hil \ot \bra{e_l}) S_\ka = \dyad{e_l}{e_l}
\]
gives
\begin{equation}
\label{for id}
E^{e_i} \big( \id_\Ov \otm \Si_\ka \big) (T) E_{e_j}  = E^{e_i}
(I_{\Kil'} \ot S_\ka) T (I_\Kil \ot S^*_\ka) E_{e_j} = p(e_i) T p(e_j).
\end{equation}
However
\[
p(e_i) T p(e_j) = E^{e_i} T E_{e_j} \ot \dyad{e_i}{e_j} \in \Ov \ot
\Cpt(\hil),
\]
and so inclusion~\eqref{Aa} holds. Since $\id_\Ov \otm \Ups_\ka$ is a
left inverse for $\id_\Ov \otm \Si_\ka$,~\eqref{Ac} follows
from~\eqref{Aa}.
\end{proof}

\begin{rems}
Note that under the identification~\eqref{op matrix},
equation~\eqref{for id} says that $\big( \id_\Ov \otm \Si_\ka \big)
(T)$ is a matrix of matrices in which the $(i,j)$-block has only one
nonzero component, namely $T^i_j$ in the $(i,j)$-place of that block.

The map $\id_{B(\hil)} \vot \Ups_\ka$ takes a matrix of matrices and
transforms it into a matrix by selecting appropriately from each
block. This is an instructive viewpoint for appreciating the next few
results.
\end{rems}

The $\ka$-\emph{Schur product} on $B(\Kil \ot \hil)$ is defined by
\[
R \schur_\ka T := \big( \id_{B(\Kil)} \vot \Ups_\ka \big) \big( (R \ot
I_{\hil}) (I_{\Kil} \ot \Pi) (T \ot I_{\hil}) (I_\Kil \ot \Pi)
\big)
\]
where $\Pi$ denotes the unitary tensor flip on $\hil \ot \hil$. The
terminology is justified since, under the identification~\eqref{op
matrix},
\[
(R \schur_\ka T)^i_j = R^i_j T^i_j.
\]

The following generalisation of a well-known property of Schur
products of complex matrices is easily verified.

\begin{lemma}
\label{positivity for Schur}
Let $R \in B(\Kil \ot \hil)$ and $T \in I_\Kil \ot B(\hil)$. If $R$
and $T$ are both nonnegative then so is $R \schur_\ka T$.
\end{lemma}

\subsection{$\ka$-decomposability}

Let $\hil$ be a Hilbert space with orthonormal basis
$\ka = (e_i)_{i \in \Ix}$. A map $\phi \in B(\Ov \otm B(\hil); \Ow
\otm B(\hil))$ is called $\ka$-\emph{decomposable} if it satisfies any
of the following equivalent conditions: for each $T \in \Ov \otm
B(\hil)$ and $i,j \in \Ix$,
\begin{rlist}
\item
$E^{e_i} \phi (T) E_{e_j} = E^{e_i} \phi (p(e_i) T p(e_j)) E_{e_j}$;
or
\smallskip
\item
$p(e_i) \phi (T) p(e_j) = p(e_i) \phi (p(e_i) T p(e_j)) p(e_j)$;
or
\smallskip
\item
$p(e_i) \phi (T) p(e_j) = \phi (p(e_i) T p(e_j))$.
\end{rlist}
The same definition applies to maps $\psi \in B( \Ov \ot \Cpt (\hil);
\Ow \ot \Cpt (\hil))$. We use the respective notations
\[
B_\kadec \big(\Ov \otm B(\hil); \Ow \otm B(\hil) \big) \quad
\text{and} \quad B_\kadec \big(\Ov \ot \Cpt (\hil); \Ow \ot \Cpt
(\hil) \big)
\]
and similarly for spaces of completely bounded $\ka$-decomposable
maps. Note that these subspaces are norm closed.

Thus a $\ka$-decomposable map $\phi: \Ov \otm B(\hil) \To \Ow \otm
B(\hil)$ has \emph{Schur action}. That is, under the
identifications~\eqref{op matrix} determined by $\ka$, for each $i,j
\in \Ix$,
\begin{equation}
\label{map matrix}
\begin{gathered}
\phi(T)^i_j = \phi^i_j (T^i_j), \ \text{ for the map } \\
\phi^i_j: a \mapsto E^{e_i} \phi \big( E_{e_i} a E^{e_j} \big) E_{e_j}.
\end{gathered}
\end{equation}
The resulting Schur-action map $[\phi^i_j] \schur : \Mat_\Ix (\Ov)_\bd
\To \Mat_\Ix (\Ow)_\bd$ is denoted ${}^\ka\! \phi$. Note that if
$\varphi \in CB (\Ov; \Ow)$ then $\varphi \otm \id_{B(\hil)}$ is
$\ka$-decomposable, with $(\varphi \otm \id_{B(\hil)})^i_j = \varphi$
for every~$i$ and~$j$.

Note also that if $\phi \in B_\kadec \big( \Ov \otm B(\hil); \Ow \otm
B(\hil) \big)$ then $\phi \big( a \ot \dyad{e_i}{e_j} \big) = \phi^i_j
(a) \ot \dyad{e_i}{e_j}$, for $a \in \Ov$ and $i,j \in \Ix$.
It follows that
\[
\phi \big( \Ov \ot \Cpt (\hil) \big) \subset \Ow \ot \Cpt (\hil).
\]
Restriction therefore induces a map
\begin{equation}
\label{PSI}
\Phi: CB_\kadec \big( \Ov \otm B(\hil); \Ow \otm B(\hil) \big) \To
CB_\kadec \big( \Ov \ot \Cpt (\hil); \Ow \ot \Cpt (\hil) \big).
\end{equation}
In Proposition~\ref{new propn} below, we see that there is a
satisfactory extension map inverting $\Phi$.

\begin{lemma}
\label{A}
Let $\phi \in CB_\kadec (\Ov \otm B(\hil); \Ow \otm B(\hil))$. Then
\begin{equation}
\label{phi circ}
\big( \phi \otm \id_{B(\hil)} \big) \circ \big( \id_\Ov \otm \Si_\ka \big)
= \big( \id_\Ow \otm \Si_\ka \big) \circ \phi.
\end{equation}
\end{lemma}

\begin{proof}
For $T \in \Ov \otm B(\hil)$ and $i, j \in \Ix$, applying~\eqref{for
id} gives both
\[
E^{e_i} \big( \phi \otm \id_{B(\hil)} \big) \circ \big( \id_\Ov \otm
\Si_\ka \big) (T) E_{e_j} = \phi \big( E^{e_i} (\id_\Ov \otm \Si_\ka)
(T) E_{e_j} \big) = \phi (p(e_i) T p(e_j))
\]
and
\[
E^{e_i} \big( (\id_\Ow \otm \Si_\ka) \circ \phi \big) (T) E_{e_j} =
p(e_i) \phi(T) p(e_j).
\]
Thus~\eqref{phi circ} follows by $\ka$-decomposability of $\phi$.
\end{proof}

By~\eqref{Aa}, the map $ \id_\Ov \otm \Si_\ka$ co-restricts to a map
$\big( \id_\Ov \otm \Si_\ka \big)' \in CB \big( \Ov \otm B(\hil); (\Ov
\ot \Cpt(\hil) ) \otm B(\hil) \big)$.

\begin{propn}
\label{new propn}
The map $\Phi$ defined in~\eqref{PSI} is a completely isometric
isomorphism whose inverse $\Psi$ is the extension map given by
\begin{equation}
\label{eqn: Psi}
\Psi: \psi \mapsto \big( \id_\Ow \otm \Ups_\ka \big) \circ \big( \psi
\otm \id_{B(\hil)} \big) \circ \big( \id_\Ov \otm \Si_\ka \big)'.
\end{equation}
\end{propn}

\begin{proof}
Let $\Psi$ be the map~\eqref{eqn: Psi}. Then, for $\psi \in CB_\kadec
\big( \Ov \ot \Cpt (\hil); \Ow \ot \Cpt (\hil) \big)$, $T \in \Ov \otm
B(\hil)$ and $i, j \in \Ix$,
\begin{align*}
E^{e_i} \Psi (\psi) (T) E_{e_j} &= E^{e_i \ot e_i} (\psi \otm
\id_{B(\hil)}) \big( (\id_\Ov \otm \Si_\ka) (T) \big) E_{e_j \ot e_j} \\
&= E^{e_i} \psi \big( E^{e_i} (\id_\Ov \otm \Si_\ka ) (T) E_{e_j}
\big) E_{e_j} \\
&= E^{e_i} \psi \big( p(e_i) T p(e_j) \big) E_{e_j}
\end{align*}
by~\eqref{for id}. Thus $\Psi (\psi)$ is $\ka$-decomposable. To see
that $\Psi$ inverts $\Phi$, pick $\phi \in CB_\kadec \big( \Ov \otm
B(\hil); \Ow \otm B(\hil) \big)$ and set $\psi = \Phi (\phi)$ then
$\psi (p(e_i) T p(e_j)) = \phi (p(e_i) T p(e_j))$ and so
\[
E^{e_i} \Psi (\psi) (T) E_{e_j} = E^{e_i} \phi(T) E_{e_j}
\]
by $\ka$-decomposability of $\phi$ and $\psi$. Thus $\Psi (\psi) =
\phi$, hence $\Psi$ is a left inverse for $\Phi$. On the other hand, let
$\psi \in CB_\kadec \big( \Ov \ot \Cpt (\hil); \Ow \ot \Cpt (\hil) \big)$
and $T \in \Ov \ot \Cpt(\hil)$ then $\psi(T) \in \Ow \ot \Cpt(\hil)$ and,
by $\ka$-decomposability of
$\psi$,
\[
E^{e_i} \psi(T) E_{e_j} = E^{e_i} \psi (p(e_i) T p(e_j)) E_{e_j} =
E^{e_i} \Psi (\psi) (T) E_{e_j}.
\]
It follows that $\Phi( \Psi(\psi) )= \psi$. Thus $\Psi$ is a right
inverse too.

Now $\Psi$ is a composition of complete contractions and $\Phi$ is
clearly completely contractive. Therefore $\Phi$ is a complete
isometry and the result follows.
\end{proof}

\begin{rem}
Clearly the extension map $\Psi$ preserves complete positivity when
$\Ov$ and $\Ow$ are $C^*$-algebras or operator systems.
\end{rem}

We use Proposition~\ref{new propn} at the end of this section. A
generalisation of the extension procedure in the following result
underlies the transformations discussed in Theorems~\ref{thm: KGagen}
and~\ref{Z}.

\begin{propn}
\label{aaa}
The prescription $\varphi \mapsto (\id_{\Ow} \otm \Ups_\ka) \circ
(\varphi \otm \id_{B(\hil)})$ defines an injective complete
contraction
\begin{equation}
\label{eqn: inj CC}
\Xi: CB \big( \Ov; \Ow \otm B(\hil) \big) \To
CB_\kadec \big( \Ov \otm B(\hil); \Ow \otm B(\hil) \big).
\end{equation}
\end{propn}

\begin{proof}
Set $\phi = (\id_{\Ow} \otm \Ups_\ka) \circ (\varphi \otm
\id_{B(\hil)})$ where $\varphi \in CB \big( \Ov; \Ow \otm B(\hil)
\big)$. Then
\[
E^{e_i} \phi(T) E_{e_j} = E^{e_i} \varphi( E^{e_i} T E_{e_j} )
E_{e_j}, \quad T \in \Ov \otm B(\hil), i, j \in \Ix.
\]
It follows that $\phi$ is $\ka$-decomposable, and that
\[
E^{e_i} \varphi(a) E_{e_j} = E^{e_i} \phi (a \ot \dyad{e_i}{e_j})
E_{e_j}, \quad a \in \Ov, i, j \in \Ix,
\]
which implies that $\varphi = 0$ if $\phi = 0$. Thus the
prescription~\eqref{eqn: inj CC} defines a linear injection $\Xi$
between the given spaces. Being a composition of complete
contractions, $\Xi$ is a complete contraction.
\end{proof}

\begin{rems}
(i) Reverting to the matrix viewpoint, for $\varphi \in CB(\Ov; \Ow
\otm B(\hil))$
\begin{equation}
\label{genkamat}
\Xi(\varphi) \longleftrightarrow \big( [a_j^i] \mapsto [\varphi_j^i
(a^i_j)] \big) \quad \text{where} \quad \varphi \longleftrightarrow
\big[ \varphi_j^i := E^{e_i} \varphi ( \cdot ) E_{e_j} \big].
\end{equation}

(ii) If $\hil$ is finite dimensional then $\Xi$ is bijective since it
has right inverse $\phi \mapsto \phi \circ \be$ for the map $\be: \Ov
\To \Ov \otm B(\hil)$, $a \mapsto a \ot \dyad{e_\Ix}{e_\Ix}$ where
$e_\Ix := \sum_{i \in \Ix} e_i$.  However, if $\hil$ is infinite
dimensional then $\Xi$ need not be bijective; for example if $\Ow =
\Ov \neq \{ 0 \}$ then $\id_{\Ov \otm B(\hil)} \notin \Ran \Xi$.
\end{rems}

\subsection{Representations and semigroup generators}

We conclude this section by applying the above techniques to obtain an
extension of the Christensen--Evans characterisation of the generators
of norm-continuous, completely positive semigroups on a $C^*$-algebra
to norm-continuous, $\ka$-decomposable semigroups on a matrix space
over a $C^*$-algebra.  The argument is based on the following slight
generalisation of the well-known representation theory of $\Cpt
(\kil)$.

\begin{lemma}
\label{decomposing repns}
Let $(\Pi, \Hil)$ be a representation of $\Al \ot \Cpt (\hil)$ for a
unital $C^*$-algebra $\Al$ and Hilbert space $\hil$. Then there is a
unital representation $(\pi, \Kil)$ of $\Al$ and isometry $V \in
B(\Kil \ot \hil; \Hil)$ satisfying
\begin{equation}
\label{V pi compact}
\Pi (A) = V \big( \pi \ot \id_{\Cpt (\hil)} \big) (A) V^*, \quad A \in
\Al \ot \Cpt (\hil).
\end{equation}
\end{lemma}

\begin{proof}
We may assume that $\hil \neq \{0\}$ and so choose a unit vector $e$
in $\hil$. Set $\Kil := P\Hil$, where $P \in B(\Hil)$ is the
orthogonal projection $\Pi (1_\Al \ot \dyad{e}{e})$, and let $J$ be
the inclusion $\Kil \To \Hil$. Then $JJ^* = P$ so the prescription
$a \mapsto J^* \Pi \big( a \ot \dyad{e}{e} \big) J$ defines a unital
representation $\pi: \Al \To B(\Kil)$. Since
\[
\ip[\Big]{\Pi \big(1_\Al \ot \dyad{x}{e} \big) J\xi}{\Pi \big(1_\Al
\ot \dyad{x'}{e} \big) J\xi'} = \ip{x}{x'} \ip{\xi}{\xi'}
\]
for $\xi, \xi' \in \Kil$, $x, x' \in \hil$, there is a unique isometry
$V \in B( \Kil \ot \hil; \Hil)$ satisfying
\[
V( \xi \ot x ) = \Pi \big( 1_\Al \ot \dyad{x}{e} \big) J \xi,
\]
for all $\xi \in \kil$ and $x \in \hil$. To establish~\eqref{V pi
  compact} it suffices, by linearity and continuity, to show that
\begin{equation}
\label{matelforrepn}
\ip[\big]{\zeta}{\Pi (a \ot \dyad{x}{x'}) \zeta'} = \ip{\zeta}{V
(\pi(a) \ot \dyad{x}{x'}) V^* \zeta'}
\end{equation}
for all $a \in \Al$, $x,x' \in \hil$ and $\zeta, \zeta' \in \Hil$.
But, for such elements, $\Pi (a \ot \dyad{x}{x'}) \zeta' = V(\xi' \ot
x)$ where $\xi' = J^* \Pi (a \ot \dyad{e}{x'}) \zeta'$, so it is in
fact sufficient to take $\zeta$ and $\zeta'$ of the form $V(\xi \ot
y)$ and $V(\xi' \ot y')$ where $\xi, \xi' \in \Kil$ and $y, y' \in
\hil$. For such choices it is easily verified that equality holds
in~\eqref{matelforrepn} with common value
$\ip{\xi}{\pi(a)\xi'}\ip{y}{x}\ip{x'}{y'}$.
\end{proof}

\begin{thm}
\label{decomp CE}
Let $\Al$ be a unital $C^*$-algebra acting nondegenerately on
$\init$, let $\hil$ be a Hilbert space with orthonormal basis $\ka$
and let $P = (P_t)_{t \ges 0}$ be a norm-continuous semigroup of
$\ka$-decomposable bounded operators on the operator system $\Al \otm
B(\hil)$ with generator $\psi$. Then $P$ is completely positive if and
only if $\psi$ has the form
\begin{equation}
\label{CP gen}
\psi: A \mapsto T ( \pi \otm \id_{B(\hil)} ) (A) T^* + NA + AN^*
\end{equation}
for a unital representation $(\pi, \Kil)$ of $\Al$ and operators $T
\in B(\Kil; \init) \vot B(\hil)$ and $N \in B(\init \ot \hil)$.
Moreover, $T$ and $N$ are necessarily $\ka$-diagonal, and may be
chosen so that $N \in \Al'' \, \vot B(\hil)$.
\end{thm}

\begin{proof}
Let $\psi^0$ and $P^0$ denote the $\ka$-decomposable operator and
semigroup on $\Al \ot \Cpt (\hil)$ obtained by restriction of $\psi$
and $P$. Since completely positive maps as well as maps of the
form~\eqref{CP gen} are completely bounded, Proposition~\ref{new
propn} and its accompanying remark imply that $P$ is a completely
positive semigroup if and only if $P^0$ is, and $\psi$ is given
by~\eqref{CP gen}, for some unital representation $(\pi, \Kil)$ of
$\Al$ and ($\ka$-diagonal) operators $T \in B(\Kil; \init) \vot
B(\hil)$ and $N \in B(\init \ot \hil)$, if and only if $\psi^0$ has
the form
\begin{equation}
\label{psi zero}
\psi^0: A \mapsto T \big( \pi \ot \id_{\Cpt (\hil)} \big) (A) T^* +NA
+AN^*
\end{equation}
for the same $\pi$, $\Kil$, $T$ and $N$.

Now if $\psi^0$ has the form~\eqref{psi zero} then it is conditionally
completely positive, and hence generates a completely positive
semigroup (\cite{QuTome}, Theorem~4.27).

Conversely, if $\psi^0$ is the generator of a completely positive
semigroup then, by Theorem~3.1 of~\cite{ChrisEvans}, Stinespring's
Theorem and Lemma~\ref{decomposing repns}, it must have the
form~\eqref{psi zero}, for a unital representation $(\pi, \Kil)$ of
$\Al$ and operators $T \in B(\Kil; \init) \vot B(\hil)$ and $N \in
\big(\Al \ot \Cpt (\hil)\big)'' = \Al'' \, \vot B(\hil)$. It therefore
only remains to show that the $\ka$-decomposability of $\psi^0$
implies that $T$ and $N$ must be $\ka$-diagonal. We may assume that
$\dim \hil \ges 2$. Let $\ka = (e_i)_{i \in \Ix}$ and pick $j \in
\Ix$. Then
\[
\psi^0 \big( 1_\Al \ot \dyad{e_j}{e_j} \big) = T E_{e_j} E^{e_j} T^* +
N E_{e_j} E^{e_j} + E_{e_j} E^{e_j} N^*
\]
and so, by the $\ka$-decomposability of $\psi^0$,
\[
E^{e_i} T E_{e_j} E^{e_j} T^* E_{e_k}+ \de_{j, k} E^{e_i} N E_{e_j} +
\de_{i,j} E^{e_j} N^* E_{e_k} = 0
\]
unless $i = k = j$, where $\de$ is the Kronecker delta. Setting $i
= k \neq j$ shows that $T$ is $\ka$-diagonal. Setting $i \neq j = k$
shows that $N$ is $\ka$-diagonal too.
\end{proof}

\begin{rem}
Under the identification $\Al \otm B(\hil) \cong \Mat_\Ix (\Al)_\bd$,
induced by the choice of basis $\ka$, and corresponding
identifications for the $\ka$-decomposable maps (as in~\eqref{op
matrix} and~\eqref{map matrix}), $\psi$ is given by ${}^\ka \psi =
[\psi^i_j] \schur \,$ where
\[
\psi^i_j: a \mapsto t_i \pi (a) t^*_j + n_i a + a n_j^*
\]
for $t_i =E^{e_i} T E_{e_i}$ and $n_i = E^{e_i} N E_{e_i}$.
\end{rem}

\section{Quantum stochastic cocycles and associated semigroups}
\label{cocycle cty}

For the rest of the paper we let $B(\init; \init')$ be the ambient
full operator space of $\Ov$ (or $B(\init)$ in case $\Ov$ is an
operator system or $C^*$-algebra). We also let $\noise$ be a fixed
Hilbert space, referred to as the \emph{noise dimension space}. Set
\begin{equation}
\label{khat defn}
\khat = \Comp \op \noise \quad \text{and, for each } x \in \noise,
\quad \wh{x} := \begin{pmatrix} 1 \\ x \end{pmatrix}.
\end{equation}
The \emph{quantum It\^{o} projection} is the orthogonal projection
$\De \in B(\khat)$ with range $\{0\} \oplus \noise$; whenever there is
no danger of confusion, its ampliations are denoted by the same
symbol. When $\eta = (d_i)_{i \in \Ix_0}$ is an orthonormal basis for
$\noise$ (with $0 \notin \Ix_0$), we set $\Ix := \{0\} \cup \Ix_0$ and
define an orthonormal basis $\etaaug = (e_\al)_{\al \in \Ix}$ for
$\khat$ by
\begin{equation}
\label{eta bar defn}
e_0 := \wh{d_0} = \begin{pmatrix} 1 \\ 0 \end{pmatrix} \text{ where }
d_0 := 0 \in \noise \text{ and, for } \ i \in \Ix_0, \ e_i
:= \begin{pmatrix} 0 \\ d_i \end{pmatrix};
\end{equation}
thus
\begin{equation}
\label{eqn: 2.3}
\Til(\eta) :=
\{ d_\al: \al \in \Ix \} =  \{0\} \cup \{ d_i: i \in \Ix_0 \}
\end{equation}
is a total subset of $\noise$ containing $0$.

For each subinterval $J$ of $\Rplus$ let $\Fock_J$ denote the
symmetric Fock space over $L^2 (J; \noise)$, dropping the subscript
when $J = \Rplus$. We use \emph{normalised} exponential vectors
\[
\w{f} := e^{-\frac{1}{2} \norm{f}^2} \big( (n!)^{-1/2} f^{\ot n}
\big)_{n \ges 0}, \quad f \in L^2 (\Rplus; \noise);
\]
these are linearly independent in $\Fock$. For any subset $\Til$ of
$\noise$ containing $0$, we set
\[
\Exps_\Til := \Lin \{ \w{f}: f \in \Step_\Til \},
\]
where $\Step_\Til$ denotes the collection of (right continuous)
$\Til$-valued step functions in $L^2 (\Rplus; \noise)$. The subscript
is dropped when $\Til = \noise$. If $\Til$ is total in $\noise$ then
$\Exps_\Til$ is total in $\Fock$ (\cite{Skeide}). For a proof of this,
together with the basics of quantum stochastic analysis and QS
cocycles, and extensive references, we refer to~\cite{Lgreifswald};
see also~\cite{FagPryc}. The \emph{Fock--Weyl operators} are the
unitary operators on $\Fock$ defined by continuous linear extension of
the following prescription:
\[
W(f): \w{g} \mapsto e^{-\I \im \ip{f}{g}} \w{f+g}, \qquad f, g \in L^2
(\Rplus; \noise).
\]

\subsection{Quantum stochastic cocycles}

By a \ti{quantum stochastic process on $\Ov$ with noise dimension
space $\noise$} we mean a family of linear maps
\[
k_t: \Ov \To L(\Exps; \Ov \otm \ket{\Fock}) \subset L(\init \aot
\Exps; \init' \ot \Fock), \quad t \in \Rplus,
\]
which is adapted and pointwise weakly measurable:
\begin{align*}
& k_t (a) E_{\w{g_{[0,t[}}} \in \Ov \otm \ket{\Fock_{[0,t[}} \ot
\ket{\w{0|_{[t,\infty[}}}, \quad \text{and} \\
& s \mapsto \ip{\zeta}{k_s (a) u \w{g}} \ \text{is measurable}
\end{align*}
for all $a \in \Ov$, $ t \ges 0$, $g \in \Step$, $u \in \init$ and
$\zeta \in \init' \ot \Fock$.  We refer to $\init \aot \Exps$ as the
\emph{exponential domain} for the process $k$. A (\ti{completely})
\ti{bounded}, \ti{completely positive} or \ti{completely contractive}
QS process on $\Ov$ is a QS process $k$ on $\Ov$ for which each $k_t$
has that property, in which case $k_t(a)$ will also denote the
continuous extension of this bounded operator to all of $\init \ot
\Fock$. Here we invoke the natural identification (cf.\ \cite{spawn},
Proposition~2.4) and inclusion
\begin{align*}
CB (\Ov; \Ov \otm B(\Fock)) &= CB (\Ov; CB (\ket{\Fock}; \Ov \otm
\ket{\Fock})) \\
&\subset L (\Ov; L (\Exps; \Ov \otm \ket{\Fock}).
\end{align*}
Adaptedness implies that a process $k$ is determined by the family of
functions $\{k^{f,g}: f,g \in \Step\}$ in $L(\Ov)$ defined by
\[
k^{f,g}_t := E^{\w{f_{[0,t[}}} k_t (\cdot) E_{\w{g_{[0,t[}}}, \quad t
\ges 0.
\]
A QS process $k$ on $\Ov$ is a (\ti{weak}) \ti{quantum stochastic
cocycle on $\Ov$} if
\begin{equation}
\label{weak cocycle}
k^{f,g}_0 = \id_\Ov \quad \text{and} \quad k^{f,g}_{r+t} = k^{f,g}_r
\circ k^{s^*_r f, s^*_r g}_t
\end{equation}
for all $f, g \in \Step$ and $r,t \ges 0$
(\cite{father},~\cite{sesquiBanach}). Here $(s^*_r)_{r \ges 0}$
denotes the semigroup of left shifts on $L^2 (\Rplus; \noise)$. The
collection of QS cocycles on $\Ov$ with noise dimension space $\noise$
is denoted $\QSC(\Ov, \noise)$. We say that $k \in \QSC(\Ov, \noise)$
has a \emph{hermitian conjugate QS cocycle} if
\[
\Dom k_t(a)^* \supset \init' \aot \Exps \text{ for all } t \in \Rplus,
a \in \Ov, \text{ and } k^\dagger: \Bigl( a^* \mapsto k_t(a)^*|_{\init'
\aot \Exps} \Bigr)_{t \ges 0} \in \QSC( \Ov^\dagger, \noise)
\]
where $\Ov^\dagger$ is the adjoint operator space $\{a^*: a \in
\Ov\}$. When $\Ov$ is an operator system or $C^*$-algebra (or, more
generally, when $\Ov$ is \emph{hermitian}: $\Ov^\dagger = \Ov$), we set
\begin{equation}
\label{eqn: QSCh}
\QSCh(\Ov, \noise) :=
\{ k \in \QSC(\Ov, \noise): k \text{ is hermitian}\}
\end{equation}
where ``$k$ is hermitian'' amounts to $k_t(a)^* \supset k_t (a^*)$,
for all $t \in \Rplus$ and $a \in \Ov$. When $\Ov$ is an operator
system we also set
\begin{equation}
\label{eqn: QSCu}
\QSCu(\Ov, \noise) :=
\{ k \in \QSC(\Ov, \noise): k \text{ is unital}\}
\end{equation}
where ``$k$ is unital'' means that $k_t(1) \subset I_{\init \ot \Fock}$
for all $t \in \Rplus$.

When $k$ is a completely bounded QS process, condition~\eqref{weak
cocycle} is equivalent to the following, which is more recognisable as
a cocycle identity:
\begin{equation}
\label{coc defn}
k_0 = \iota_\Fock^\Ov \quad \text{and} \quad k_{r+t} = \wh{k}_r \circ
\si_r \circ k_t, \quad r,t \ges 0,
\end{equation}
where $\iota_\Fock^\Ov$ is the ampliation introduced in~\eqref{eqn:
iota} and $\wh{k}_r$ denotes the natural (matrix-space) tensor
extension of $k_r$ to the range of $\si_r$ (see \cite{spawn},
Section~5). We refer to such processes as \ti{completely bounded QS
  cocycles}, and denote the class of these by $\QSCcb(\Ov,
\noise)$. When $\Ov$ is a $C^*$-algebra or operator system, we
similarly write
\[
\QSCcp(\Ov, \noise), \quad \QSCcpc(\Ov, \noise) \ \text{ and } \
\QSCcpu(\Ov, \noise)
\]
for the respective subclasses of completely positive, completely
positive contractive and completely positive unital QS cocycles.

\begin{rem}
In \cite{spawn} all QS cocycles were assumed to be completely bounded
(and the measurability condition was not imposed). The reason for
dealing with weak QS cocycles here, and elsewhere, is that solutions
of the QS differential equation~\eqref{QSDE} with completely bounded
coefficients are of this type --- but $k$ need not be completely
bounded. Indeed, $k_t (a)$ need not even be a bounded operator.
\end{rem}

QS cocycles $k$ have \emph{associated semigroups}, defined by
\[
\msgp^{x,y}_t = E^{\w{x_{[0,t[}}} k_t (\cdot) E_{\w{y_{[0,t[}}}, \quad x,y
\in \noise,
\]
and, for each $f,g \in \Step$, the semigroup decomposition of QS
cocycles represents $k^{f,g}_t$ in terms of these semigroups
(\cite{spawn}, Proposition 5.1). When it is bounded a QS cocycle is
thereby determined by the family of semigroups $\{ \msgp^{x,y}: x,y
\in \Til \}$ for any total subset $\Til$ of $\noise$ containing $0$. A
QS cocycle is \emph{elementary} if each of its associated semigroups
is completely bounded and cb-norm continuous. \emph{Note.} We have
previously used the terminology `cb-Markov regular'. We write
$\El\QSC(\Ov, \noise)$ for this class of cocycle. For a completely
bounded QS cocycle $k$ which is locally bounded in cb-norm, cb-norm
continuity for any one of the associated semigroups, such as the
\emph{vacuum-expectation semigroup} $\msgp^{0,0}$, implies that $k$ is
elementary. This is a simple consequence of the continuity of the
normalised exponential map $\varpi: L^2(\Rplus; \noise) \To
\Fock$. Similarly for such cocycles either all or none of the
associated semigroups are $C_0$-semigroups. Part~(c) of
Theorem~\ref{M} below includes the analogous result for continuity in
cb-norm.

A completely bounded QS cocycle $k$ with locally bounded cb-norm is
necessarily exponentially bounded in cb-norm: there is $M \ges 1$
and $\be \in \Real$ such that
\begin{equation}
\label{eqn: ktcbM}
\cbnorm{k_t} \les M e^{\be t}, \qquad t \in \Rplus
\end{equation}
(\cite{spawn}, Proposition~5.4). Following the terminology of
semigroup theory, we refer to the quantity $\inf \big\{ \beta \in
\Real: \sup_{t \ges 0} e^{- \beta t} \norm{k_t}_{\cb} < \infty \big\}
\in [-\infty, \infty]$ as the (\emph{exponential}) \emph{growth bound}
of $k$, set
\begin{align*}
&\QSCbeta (\Ov, \noise) := \big\{ k \in \QSCcb (\Ov, \noise):
\eqref{eqn: ktcbM} \text{ holds with } M=1 \big\}, \text{ and } \\
&\QSCqc (\Ov, \noise) := \bigcup_{\beta \in \Real} \QSCbeta (\Ov, \noise),
\end{align*}
and call elements of the latter class of cocycles
\emph{cb-quasicontractive}. Thus, for a $C^*$-algebra $\Al$,
$\El\QSCcpqc(\Al,\noise)$ denotes the class of completely positive
quasicontractive elementary QS cocycles on $\Al$ with noise dimension
space $\noise$.

\begin{rem}
Let $k \in \QSCcpqc(\Al,\noise)$ for a $C^*$-algebra $\Al$. If its
vacuum expectation semigroup is norm continuous then, by the
Christensen--Evans theorem (\cite{ChrisEvans}), it has completely
bounded generator and so is cb-norm continuous, and therefore (being
locally bounded in cb-norm) the QS cocycle $k$ is elementary.
\end{rem}

Let $\init$ be a Hilbert space. A suitably measurable, Fock-adapted
family of operators $X = (X_t)_{t \ges 0}$ in $B(\init \ot \Fock)$ is
a \emph{bounded left quantum stochastic operator cocycle} on $\init$
(with noise dimension space $\noise$) if
\[
X_0 = I_{\init \ot \Fock} \text{ and } X_{r+t} = X_r \si_r (X_t),
\quad r,t \in \Rplus.
\]
For $x, y \in \noise$,
$\osgp^{x,y} := \big( E^{\w{x_{[0,t[}}} X_t E_{\w{y_{[0,t[}}} \big)_{t
\ges 0}$ defines its $(x,y)$-\emph{associated semigroup} of operators
on $\init$; $X$ is \emph{elementary} if all of these are norm
continuous. Again, if $X$ is locally bounded then norm continuity for
any one of its associated semigroups implies that $X$ is
elementary. Important examples are the \emph{Weyl cocycles},
$(W(z_{[0,t[}))_{t \ges 0}$ ($z\in\noise$), much used in quantum
stochastic calculus (see \cite{Lgreifswald}, pages 247--9). We use the
notations $\El\QSCb (\init, \noise)$, $\QSCqc (\init, \noise)$ etc.

Connections between operator cocycles and mapping cocycles are
described in Proposition~5.5 of~\cite{spawn}. In this paper we shall
use the following one: given bounded left QS operator cocycles $X$ and
$Y$ on $\init$, the prescription
\[
k_t(T) := X_t ( T \ot I_\Fock ) Y_t^* \quad (T \in B(\init), t \in
\Rplus)
\]
defines a completely bounded cocycle $k$ on $B(\init)$.

\subsection{Associated $\Ga$-cocycle and global $\Ga$-semigroup}

An important idea of Accardi and Koz\-y\-rev (\cite{AK}), expanded on
at length in~\cite{spawn}, is to gather associated semigroups of a QS
cocycle into a matrix. For any index set $\Ix$, let $(\de^\al)_{\al
  \in \Ix}$ denote the standard orthonormal basis for $l^2(\Ix)$.
Then, given a map $\Ga: \Ix \To \noise$, define unitaries $W^\Ga_t \in
B(l^2(\Ix) \ot \Fock)$, for $t \in \Rplus$, by continuous linear
extension of the rule
\begin{equation}
\label{diag Weyl}
\de^\al \ot \xi \mapsto \de^\al \ot W \big( \Ga (\al)_{[0,t[}
\big) \xi, \quad \al \in \Ix, \xi \in \Fock.
\end{equation}
Note that $W^\Ga_t \in \Diag_\Ix \big( B(\Fock) \big)_\bd$ under the
identification~\eqref{op matrix}. We adapt this notation in two
special cases: if $\Ix = \{ 0 \} \cup \Ix_0$ and $\eta = (d_i)_{i \in
\Ix_0}$ is a basis for $\noise$ then we write $W^\eta = W^\Ga$ for the
map $\Ga: \Ix \To \noise$, $\al \mapsto d_\al$ (recall our
convention~\eqref{eta bar defn}); if $n \in \Nat$ and $\bfx = (x_1,
\cdots, x_n) \in \noise^n$, set $W^\bfx = W^\Ga$ for the map $\Ga:
\{1, \cdots, n\} \To \noise$, $i \mapsto x_i$.

The following is Proposition~5.9 of~\cite{spawn}. Recall the
notation~\eqref{flipped lift}.

\begin{propn}
\label{kGaPGa}
Let $k \in \QSCcb(\Ov, \noise)$, with associated semigroups
$\{\msgp^{x,y}: x,y \in \noise\}$, and let $B(\init;\init')$ be the
ambient full operator space of $\Ov$. Then, for any set $\Ix$ and map
$\Ga: \Ix \To \noise$, setting $\hil = l^2(\Ix)$,
\[
\big( (I_{\init'} \ot W^\Ga_t)^* k^{\hil}_t (\cdot) (I_\init \ot
W^\Ga_t) \big)_{t \ges 0}
\]
defines a cocycle $k^\Ga \in \QSCcb \big( \Ov \otm B (\hil),
\noise \big)$ that has Schur-action under the
identifications~\eqref{op matrix}\tu{;} its vacuum-expectation
semigroup is given by the Schur-action prescription
\[
\msgp^\Ga_t : [ a^\al_\be ] \mapsto \Big[ \msgp^{\Ga(\al), \Ga(\be)}_t
(a^\al_\be) \Big]
\]
where, for each $x, y \in \noise$, $\msgp^{x,y}$ is the
$(x,y)$-associated semigroup of $k$.
\end{propn}

\begin{rems}
We use the notation $k^\eta$, $\msgp^\bfx$, etc.\ in the special cases
noted above for $\Ga$. Also, for each set $\Ix$ and map $\Ga: \Ix \To
\noise$, the semigroup $\msgp^\Ga$ clearly leaves $\Ov \ot \Cpt(\hil)$
invariant; we denote the resulting restriction by $\msgp^{\Ga, \Cpt}$.
\end{rems}

\subsection{Continuity}

We now collect together the basic relationships between continuity for
QS cocycles and continuity for associated semigroups. We also
highlight some implications of strengthened continuity hypotheses
which render them inappropriate. Some hybrid locally convex topologies
on spaces of the form $B(\Hil \ot \hil; \Hil' \ot \hil)$ are
needed. These are the $\hil$-\emph{ultraweak topology} defined by the
seminorms
\[
p^\om: T \mapsto \norm{ (\id_{B(\Hil;\Hil')} \vot \om) (T)}, \quad \om
\in B(\hil)_*,
\]
and the $\hil$-\emph{weak topology} defined by the seminorms
\[
p_{x,y}: T \mapsto \norm{E^x T E_y},
\quad x,y \in \hil.
\]
Thus $p_{x,y} = p^{\omega_{x,y}}$ for the vector functional
\begin{equation}
\label{vec fnl}
\om_{x,y}: R \mapsto \ip{x}{Ry},
\quad x,y \in \hil.
\end{equation}
Note that, for any operator space $\Ov$ in $B(\Hil; \Hil')$, the
matrix space $\Ov \otm B(\hil)$ is the closure of the algebraic tensor
product $\Ov \aot B(\hil)$ with respect to either of these topologies,
and the two topologies coincide on norm-bounded subsets of
$B(\Hil \ot \hil; \Hil' \ot \hil)$. See Sections~2 and~3
of~\cite{spawn} for more details. Moreover, $\Ov \aot \Cpt (\hil)$ is
dense in $\Ov \otm B(\hil)$ with respect to the weak and strong
operator topologies on $B(\Hil \ot \hil; \Hil' \ot \hil)$.

For a set $\Ix$ and \ti{bounded} map $\Ga: \Ix \To \noise$, set
\[
w^\Ga_t:= W^\Ga_t E_{\w{0}} \in B (l^2(\Ix); l^2(\Ix) \ot \Fock),
\quad t \in \Rplus.
\]

\begin{lemma}
\label{lemma: F Gamma}
For all $t \in \Rplus$, $\norm{w^\Ga_t - w^\Ga_0} = \sqrt{2 (1 -
e^{-tM})}$, where $\displaystyle{M = \sup_{\al \in \Ix} \norm{\Ga
(\al)}^2/2}$.
\end{lemma}

\begin{proof}
Let $t \in \Rplus$. Then $w^\Ga_t -w^\Ga_0 \in \Diag_\Ix
(\ket{\Fock})_\bd$ and, for all $\al \in \Ix$,
\[
\norm{(w^\Ga_t -w^\Ga_0) \de^\al}^2 = \norm{\de^\al \ot \big(
\w{\Ga(\al)_{[0,t[}} -\w{0} \big)}^2 = 2(1 -e^{-t \norm{\Ga(\al)}^2/2}).
\]
The result follows.
\end{proof}

\begin{thm}
\label{M}
Let $k \in \QSCcb(\Ov, \noise)$, with locally bounded cb-norm.
Then the following sets of equivalences hold\tu{:}
\begin{alist}
\item
\begin{alist}
\label{wc}
\item $k$ is pointwise ultraweakly continuous.
\item $\msgp^\Ga$ is pointwise ultraweakly continuous for every
$\noise$-valued map $\Ga$.
\item $\msgp^{x,y}$ is pointwise ultraweakly continuous for all $x, y
\in \noise$.
\item $\msgp^{x,y}$ is pointwise ultraweakly continuous for all
$x, y \in \Til$, for some total subset $\Til$ of $\noise$ containing
$0$.
\end{alist}
\smallskip
\noindent
Suppose that $\Ov$ is ultraweakly closed and each map $\msgp^{x,y}_t$
is ultraweakly continuous. Then there is a further equivalence\tu{:}
\begin{alist}\setcounter{enumii}{4}
\item $\msgp^{x,y}$ is pointwise ultraweakly continuous at $0$ for
some $x, y \in \noise$.
\end{alist}

\smallskip
\item
\begin{alist}
\label{sc}
\item $k$ is pointwise $\Fock$-ultraweakly continuous.
\item $\msgp^{\Ga,\Cpt}$ is a $C_0$-semigroup for every
$\noise$-valued map $\Ga$.
\item $\msgp^{x,y}$ is a $C_0$-semigroup for all $x, y \in \noise$.
\item $\msgp^{x,y}$ is a $C_0$-semigroup for some $x, y \in \noise$.
\end{alist}

\smallskip
\item
\begin{alist}
\label{cbc}
\item
The map $\Rplus \To CB(\Ov)$, $s \mapsto (\id_\Ov \otm\, \om) \circ
k_s$ is continuous, for all $\om \in B(\Fock)_*$.
\item
For every \emph{bounded} $\noise$-valued map $\Ga: \Ix \To \noise$,
$\msgp^\Ga$ is cb-norm continuous.
\item $k$ is elementary.
\item $\msgp^{x,y}$ is cb-norm continuous for some $x, y \in \noise$.
\end{alist}
\end{alist}
\end{thm}

\begin{proof}
The implication (ii)$\Rightarrow$(iii) follows in~\eqref{wc},
\eqref{sc} and~\eqref{cbc} since if $\mathbf{x} \in \noise \times
\noise$ then
\[
\big[ \msgp^{x_i,x_j} \big] = \msgp^{\mathbf{x}} = \msgp^\Ga,
\]
for the map $\Ga: \{1,2 \} \To \noise$, $i \mapsto x_i$; the
implication (iii)$\Rightarrow$(iv) is obvious in~\eqref{wc},
\eqref{sc} and~\eqref{cbc}

\smallskip
\eqref{wc}
(i)$\Rightarrow$(ii):
Assume that (i) holds, fix a map $\Ga: \Ix \to \noise$ defined on some
set $\Ix$, and set $\hil := l^2(\Ix)$, with standard orthonormal basis
$(\de^\al)_{\al \in \Ix}$. The key identity is
\[
\ip{u' \ot \de^\al}{\msgp^\Ga_t(A) u \ot \de^\be} = \ip{u' \ot \w{\Ga
(\al)_{[0,t[}}}{k_t (E^{e_\al} A E_{e_\be}) u \ot \w{\Ga
(\be)_{[0,t[}}},
\]
valid for all $A \in \Ov \otm B (\hil)$, $t \in \Rplus$, $u' \in
\init'$, $u \in \init$ and $\al, \be \in \Ix$. The result follows by
local norm-boundedness of $\msgp^\Ga$, norm-totality of the set $\{
\om_{u' \ot \de^\al, u \ot \de^\be}: u' \in \init', u \in \init, \al,
\be \in \Ix \}$ in $B(\init \ot \hil; \init' \ot \hil)_*$,
adaptedness of $k$ and continuity of the map $f \mapsto \w{f}$.

(iv)$\Rightarrow$(v):
Obvious.

(iv)$\Rightarrow$(i): This follows from the local norm-boundedness of
$k$, the semigroup decomposition of QS cocycles, and the norm-totality
in the predual space $B(\init \ot \Fock; \init' \ot \Fock)_*$ of the
set $\{\om_{u' \ot \ep', u \ot \ep}: u' \in \init', u \in \init, \ep',
\ep \in \Exps_\Til\}$. See Proposition~4.2 of~\cite{father} for full
details.

(v)$\Rightarrow$(iii):
Now suppose that $\Ov$ is ultraweakly closed, and let $\Ov_*$ denote
the space of ultraweakly continuous linear functionals on $\Ov$.
Assume that each map $\msgp^{x,y}_t$ is ultraweakly continuous and the
associated semigroup $\msgp^{x',y'}$ is pointwise ultraweakly
continuous at $0$. For each $x,y \in \noise$, let $Q^{x,y}$ denote the
semigroup on $\Ov_*$ induced by $\msgp^{x,y}$. Pointwise ultraweak
continuity for $\msgp^{x,y}$ amounts to pointwise weak continuity for
$Q^{x,y}$. Therefore, by semigroup theory (\cite{Davies}, in
particular Theorem~6.2.6), (iii) holds if each $\msgp^{x,y}$ is
pointwise ultraweakly continuous at $0$. Thus let $x,y\in \noise$,
$a \in \Ov$ and $\om \in \Ov_*$.  Then, for $t \in \Rplus$,
\begin{equation*}
\om (\msgp_t^{x,y} a) -\om (\msgp_t^{x',y'} a) =\ip{\w{x_{[0,t[}}}{g_t
(a) \w{y_{[0,t[}}} -\ip{\w{x'_{[0,t[}}}{g_t (a) \w{y'_{[0,t[}}}
\end{equation*}
where $g_t := (\om \vot \id_{B(\Fock)}) \circ k_t \in CB(\Ov;
B(\Fock))$. Therefore the required pointwise continuity at $0$ of
$\msgp^{x,y}$ follows from the local boundedness of the family
$\big\{g_t(a): t \in \Rplus \big\}$ and the continuity of the
normalised exponential map $\varpi: L^2(\Rplus; \noise) \To \Fock$.

\medskip
\eqref{sc}
(i)$\Rightarrow$(ii):
Let $\Ix$ be a set, let $\Ga: \Ix \To \noise$ be a map and let
$(\de^\al)_{\al \in \Ix}$ be the standard orthonormal basis of
$l^2(\Ix)$. For all $a \in \Ov$, $\al, \be \in \Ix$ and $t \in [0,1[$,
\[
\msgp^{\Ga, \Cpt}_t (a \ot \dyad{\de^\al}{\de^\be}) =
\ip{\w{\Ga(\al)_{[t,1[}}}{\w{\Ga(\be)_{[t,1[}}}^{-1} (\id_\Ov \otm
\om_{\xi, \zeta}) \big( k_t (a) \big) \ot \dyad{\de^\al}{\de^\be}
\]
for $\xi = \w{\Ga(\al)_{[0,1[}}$ and $\zeta =
\w{\Ga(\be)_{[0,1[}}$. Since $\msgp^{\Ga, \Cpt}$ is locally norm
bounded, it follows by totality of $\{a \ot \dyad{\de^\al}{\de^\be}: a
\in \Ov, \al, \be \in \Ix\}$ that $\msgp^{\Ga, \Cpt}$ is strongly
continuous at $0$, and this extends to all of $\Rplus$ by standard
semigroup theory.

(iii)$\Rightarrow$(i):
This follows from the semigroup decomposition of QS cocycles and the
norm-totality of the family $\{\om_{\w{f}, \w{g}}: f,g \in \Step\}$ in
$B(\Fock)_*$.

(iv)$\Rightarrow$(iii): This has been remarked upon earlier.

\medskip
\eqref{cbc}
(i)$\Rightarrow$(ii):
Fix a set $\Ix$ and map $\Ga: \Ix \to \noise$. First note that for $t
\in \Rplus$
\begin{equation}
\label{obs}
\begin{aligned}
& \msgp^\Ga_t = (I_{\init'} \ot w^\Ga_t)^* k_t^\hil (\cdot) (I_{\init}
\ot w^\Ga_t), \text{ and} \\
& \msgp_t^{0,0} \otm \id_{B(\hil)} = (I_{\init'} \ot w^\Ga_0)^*
k_t^\hil (\cdot) (I_{\init} \ot w^\Ga_0),
\end{aligned}
\end{equation}
where $\hil = l^2(\Ix)$. Thus $\cbnorm{\msgp^\Ga_t -\msgp_t^{0,0} \otm
\id_{B(\hil)}} \Tends 0$ as $t \Tends 0$ by Lemma~\ref{lemma: F
Gamma}. The result follows by standard semigroup techniques, since
\[
\cbnorm{\msgp^{0,0}_t -\id_\Ov} = \cbnorm{(\id_\Ov \otm \om_{\w{0},
\w{0}}) \circ (k_t -k_0)}.
\]

(iii)$\Rightarrow$(i):
This follows as in part~\eqref{sc}.

(iv)$\Rightarrow$(iii):
Note that for any set $\Ix$, if $\hil = l^2(\Ix)$ then
\[
\msgp^{x,y}_t \otm \id_{B(\hil)} = (I_{\init'} \otm w^{\Ga_x}_t)^*
k^\hil_t (\cdot) (I_\init \otm w^{\Ga_y}_t)
\]
for the constant maps $\Ga_x, \Ga_y : \Ix \to \noise$, $\al \mapsto
x$, respectively $\al \mapsto y$. Thus $\lim_{t \Tends 0}
\cbnorm{\msgp^{x,y}_t -\msgp^{0,0}_t} = 0$ by similar arguments to
those used above.
\end{proof}

\begin{rems}
Proposition~\ref{auto CB reg} below shows the inappropriateness of
strong continuity on all of $\Ov \otm B(\hil)$ as an assumption for
$\msgp^\Ga$ in part~\eqref{sc}.

The hypothesis for (v) of part~\eqref{wc} rectifies an omission
in Proposition~5.4 of~\cite{father}.

Intermediate between~\eqref{sc} and~\eqref{cbc} there is the further
equivalence:
\begin{rlist}
\item  The map $\Rplus \To B(\Ov)$, $s \mapsto (\id_\Ov \otm \, \om)
\circ k_s$, is continuous, for all $\om \in B(\Fock)_*$.
\item $\msgp^{x,y}$ is norm continuous for all $x, y \in \noise$.
\item $\msgp^{x,y}$ is norm continuous for some $x, y \in \noise$.
\end{rlist}
\end{rems}

The analogous form of Theorem~\ref{M} for QS operator cocycles is
given next; it may be deduced via the natural correspondence between
bounded QS operator cocycles on $\init$ and completely bounded QS
mapping cocycles on $\ket{\init}$, the column-operator space of
$\init$ (\cite{spawn}, Proposition~5.5), or proved directly as
in~\cite{JTL}.

\begin{cor}
\label{cty for op cocycles}
Let $X \in \QSCb (\init, \noise)$ with locally bounded norm.
Then the following sets of equivalences hold\tu{:}
\begin{alist}
\item
\begin{rlist}
\item $X$ is strongly continuous.
\item $\osgp^{x,y}$ is a $C_0$-semigroup for all $x, y \in \noise$.
\item $\osgp^{x,y}$ is a $C_0$-semigroup for some $x, y \in \noise$.
\end{rlist}

\smallskip
\item
\begin{rlist}
\item The map $\Rplus \To B(\init)$, $t \mapsto (\id_{B(\init)} \vot
\om) (X_t)$ is continuous, for all $\om \in B(\Fock)_*$.
\item $X$ is elementary.
\item $\osgp^{x,y}$ is norm continuous for some $x, y \in \noise$.
\end{rlist}
\end{alist}
\end{cor}

We denote the former class of operator cocycles by $C_0$-$\QSC (\init,
\noise)$, and the latter class by $\El\QSC (\init, \noise)$.

\begin{rems}
For $X \in \QSCb (\init, \noise)$ with locally bounded norm, it is
easily seen that there are constants $M \ges 1$ and $\beta \in \Real$
such that $\norm{ X_t} \les M e^{\beta t}$ for all $t \in \Real$. Any
bounded operator QS cocycle that is weak operator continuous at $0$ is
necessarily locally norm bounded (\cite{JTL}, Proposition~2.1).
\end{rems}

We conclude this subsection with a result concerning continuity in the
arguments of $k_t$, $\msgp_t^{x,y}$ and $\msgp^\Ga_t$ for fixed $t \in
\Rplus$, when $\Ov$ is ultraweakly closed.

\begin{propn}
Let $k \in \QSCcb(\Ov, \noise)$, let $\Til$ any total subset of
$\noise$ that contains $0$, and suppose that $\Ov$ is ultraweakly
closed. Then, for $t \in \Rplus$, the following are equivalent\tu{:}
\begin{rlist}
\item
\label{kuw}
$k_t$ is ultraweakly continuous from $\Ov$ to $\Ov \vot B(\Fock)$.
\item
\label{Puw}
$\msgp^\Ga_t$ is ultraweakly continuous on $\Ov \vot B \big( l^2 (\Ix)
\big)$ for each set $\Ix$ and map $\Ga: \Ix \To \noise$.
\item
\label{Pxyuw}
$\msgp^{x,y}_t $ is ultraweakly continuous on $\Ov$ for each $x,y \in
\Til$.
\end{rlist}
\end{propn}

\begin{proof}
\eqref{kuw}\,$\Implies$\,\eqref{Puw}:
Suppose that~\eqref{kuw} holds. Since ultraweakly continuous
completely bounded maps ampliate to ultraweakly continuous maps,
$\msgp^\Ga_t$ is a composition of ultraweakly continuous maps,
so~\eqref{Puw} follows.

\eqref{Puw}\,$\Implies$\,\eqref{Pxyuw}:
This is obvious.

\eqref{Pxyuw}\,$\Implies$\,\eqref{kuw}:
This is covered by Proposition~4.2 of \cite{father}.
\end{proof}

\subsection{Strengthening continuity}

Theorem~\ref{M} identifies the relevant continuity conditions on a QS
cocycle $k$ that correspond to standard continuity conditions on its
associated semigroups. Thus if $\Ov$ is a $C^*$-algebra, for example,
then pointwise $\Fock$-ultraweak continuity is the appropriate
condition on $k$, whereas if $\Ov$ is a von Neumann algebra then
pointwise ultraweak continuity is the appropriate condition on $k$.
Our next two results show some of the hazards of moving outside these
conditions.

\begin{propn}
\label{O}
Let $k \in \QSCcp(\vN, \noise)$ for a von Neumann algebra $\vN$,
and suppose that $k$ is pointwise $\Fock$-ultraweakly continuous.
Then $k$ is elementary.
\end{propn}

\begin{proof}
Let $\{\msgp^{x,y}: x,y \in \noise\}$ be the associated semigroups of
$k$. By Theorem~\ref{M} the semigroup $\msgp^\bfx$ is strongly
continuous on $\Mat_n(\vN)$, for each $n \ges 1$ and $\bfx \in
\noise^n$. Since $\msgp^\bfx$ is completely positive (by
Proposition~\ref{kGaPGa}) this implies that $\msgp^\bfx$ is norm
continuous (\cite{strong-norm}).  Since norm-continuous completely
positive semigroups have completely bounded generators
(\cite{ChrisEvans}), it follows that $\msgp^\bfx$ is cb-norm
continuous, and thus each of its component semigroups is too.
\end{proof}

\begin{rem}
From this proposition it follows that a strongly continuous completely
positive QS contraction cocycle on $B(\init)$ necessarily has
cb-norm-continuous associated semigroups. Furthermore, the analysis
of~\cite{father} (or Sections~\ref{Section: gen}--\ref{sec: Mr CPC}
below) shows that such a cocycle is governed by the QS differential
equation~\eqref{QSDE} with a completely bounded stochastic
generator. This appears to have been overlooked in~\cite{AK}.
\end{rem}

Whilst Proposition~\ref{O} shows that insisting on strong continuity
for $k$ can result in $k$ being necessarily an elementary QS cocycle,
the following result and example shows that the converse is false:
there are elementary QS cocycles that are not strongly continuous.

\begin{propn}
\label{propn: 2.7}
Let $k$ be a strongly continuous ${}^*$-homomorphic QS cocycle on a
unital $C^*$-algebra $\Al$. Then $k$ is unital.
\end{propn}

\begin{proof}
If $k$ is nonunital then the projection-valued process $(1_\Al \ot
I_{\Fock} - k_t(1_\Al))_{t \ges 0}$ is nonzero for all $t > 0$, by
Proposition~5.8 of~\cite{spawn}. Thus $\lim_{t \Tends 0^+} \norm{k_t
(1_\Al) - 1_\Al \ot I_{\Fock}} = 1$ and so $k$ is not strongly
continuous.
\end{proof}

Since bounded derivations $\de$ on a unital $C^*$-algebra $\Al$
satisfy $\de (1_\Al) = 0$, norm continuity for a homomorphic semigroup
on $\Al$ implies unitality. This is not so for cocycles on $\Al$.

\begin{example}
Let $\pi$ be a nonunital endomorphism of a unital $C^*$-algebra $\Al$
and set $k = k^\phi$ for the completely bounded map
\[
\phi: \Al \To \Al \ot B(\Comp^2), \quad a \mapsto \begin{bmatrix}
0 & 0 \\ 0 & \pi(a) -a \end{bmatrix}.
\]
Here $k^\phi$ denotes the solution of the QS differential
equation~\eqref{QSDE} determined by $\phi$; it is a ${}^*$-homomorphic
elementary QS cocycle on $\Al$ with noise dimension space
$\Comp$. However $\phi (1_\Al) \neq 0$, so $k$ is not unital
(\cite{gran}, Theorem~5.1) thus $k$ cannot be strongly continuous (by
Proposition~\ref{propn: 2.7}). In this case the generator of the
global semigroup $\msgp^\eta$, with respect to the basis $\eta =
\{1\}$ for $\Comp$, is given by
\[
\begin{bmatrix} a & b \\ c & d \end{bmatrix} \mapsto \begin{bmatrix} 0
& -\tfrac{1}{2} b \\ -\tfrac{1}{2} c & \pi(d) -d \end{bmatrix},
\]
as follows from Theorem~\ref{Z} below.

For an example of this, with $\Al$ a separable unital algebra, take
$\Al$ to be the Toeplitz algebra, that is $C^*(T)$ where $T$ is the
right-shift operator on $l^2$, and $\pi$ to be the map
$a \mapsto TaT^*$.
\end{example}

The following lemma leads to a further automatic continuity result
that illustrates the problems that can arise if one seeks to enlarge
the domain on which the global semigroup acts in Theorem~\ref{M}\,(b).

\begin{lemma}
\label{cty on matrices}
Let $(\varphi_t)_{t > 0}$ a family of completely bounded maps on $\Ov$
and let $\Kil$ be an infinite dimensional Hilbert space. Then the
following are equivalent\tu{:}
\begin{rlist}
\item
$\varphi_t \otm \id_{B(\Kil)} \Tends 0$ strongly as $t \to 0$.
\item
$\cbnorm{\varphi_t} \Tends 0$ as $t \to 0$.
\end{rlist}
\end{lemma}

\begin{proof}
Clearly (ii) implies (i). Suppose therefore that $\cbnorm{\varphi_t}
\not\Tends 0$ as $t \to 0$. Then there is $\ep > 0$ and a decreasing
sequence $(t_n)_{n \ges 0}$ in $]0,\infty[$ with $\lim_{n \Tends
\infty} t_n = 0$ such that $\cbnorm{\varphi_{t_n }} \ges 2\ep$ for
each $n$. Choose a sequence $(\Kil_n)_{n\ges 0}$ of mutually
orthogonal finite-dimensional subspaces of $\Kil$ and a sequence
$(A_n)_{n \ges 0}$ of operators in $\Ov \ot B(\Kil_n)$ such that
$\norm{A_n} = 1$ and $\norm{(\varphi_{t_n} \ot \id_{B(\Kil_n)})(A_n)}
\ges \ep$. The block-diagonal operator $A = 0' \oplus \bigoplus A_n
\in \Ov \otm B(\Kil)$, where $0'$ denotes the zero element of $\Ov
\otm B (\Kil \ominus \bigoplus \Kil_n)$, satisfies
$\norm{\varphi_{t_n} \otm \id_{B(\Kil)} (A)} \ges \ep$ for all
$n$. Thus~(i) implies~(ii).
\end{proof}

\begin{propn}
\label{auto CB reg}
Let $k \in \QSCcb(\Ov, \noise)$ with locally bounded cb-norm, and
let $\Ga$ be a function from some set $\Ix$ to $\noise$.  Suppose
that $\Ix$ is infinite, $\Ga$ is bounded, and $\msgp^\Ga$ is strongly
continuous on $\Ov \otm B(l^2(\Ix))$. Then $k$ is elementary.
\end{propn}

\begin{proof}
In view of the identities~\eqref{obs}, $\msgp^{0,0}_t \otm
\id_{B(\hil)} \Tends \id_{\Ov \ot B(\hil)}$ strongly as $t \Tends 0$,
where $\hil = l^2(\Ix)$. Thus Lemma~\ref{cty on matrices} implies that
$\msgp^{0,0}$ is cb-norm continuous. The result therefore follows from
Theorem~\ref{M}.
\end{proof}

Finally we conclude with a result that is useful for establishing
quantum stochastic integrability of certain processes obtained from a
QS cocycle.

\begin{propn}
\label{ptwise SOT}
Let $k \in \QSCcb (\Ov, \noise)$ with locally bounded cb-norm, and
suppose that $k$ is pointwise continuous at $t=0$ with respect to the
weak operator topology. Then the following hold\tu{:}
\begin{alist}
\item
$k$ is pointwise right-continuous on all of $\Rplus$ in this
topology.
\item
If further $\Ov$ is a von Neumann algebra and $k$ is completely
positive and contractive then $k$ is pointwise right-continuous on
$\Rplus$ with respect to the strong operator topology.
\end{alist}
\end{propn}

\begin{proof}
From the cocycle identity~\eqref{coc defn} it is sufficient to show
that for any $a \in \Ov$ and $t \in \Rplus$
\[
\lim_{r \downarrow 0} \ip{u \ot \w{f}}{\big[ (\wh{k}_r \circ \si_r)
(A) -A \big] v \ot \w{g}} =0
\]
where $A = k_t(a)$. However, the definitions and identifications
in~\eqref{coc defn} yield the identity
\[
\ip{u \ot \w{f}}{(\wh{k}_r \circ \si_r) (A) v \ot \w{g}} = \ip{u \ot
\w{f_{[0,r[}}}{k_r (E^{\w{s^*_r f}} A E_{\w{s^*_r g}}) v \ot \w{g_{[0,r[}}}
\]
where $(s^*_r)_{r \ges 0}$ denotes the $C_0$-semigroup of left-shifts
on $L^2 (\Rplus; \noise)$. This fact, coupled with the assumptions on
$k$, gives the desired limit as $r \downarrow 0$.

Suppose now that $\Ov$ is a von Neumann algebra and $k$ completely
positive and contractive. Then the operator Schwarz inequality gives,
for $A = k_t(a)$ and $\xi \in \init \ot \Fock$,
\[
\norm{\big[ (\wh{k}_r \circ \si_r) (A) -A \big] \xi}^2 \les
\ip{\xi}{(\wh{k}_r \circ \si_r) (A^* A) \xi} -2 \re \ip{A
\xi}{(\wh{k}_r \circ \si_r) (A) \xi} +\norm{A \xi}^2,
\]
and the RHS tends to zero as $r \downarrow 0$ by the first part.
The second part follows.
\end{proof}

\section{QS cocycles and QS differential equations}
\label{Section: QSC QSDE}

As noted earlier, solutions of the QS differential
equation~\eqref{QSDE} provide a basic source of QS cocycles. We adopt
the notations~\eqref{flipped lift}, \eqref{eqn: iota}, \eqref{khat
  defn}, \eqref{vec fnl} and the definition
\begin{align}
&\chi: \noise \times \noise \to \Comp, \quad (x,y) \mapsto
\tfrac{1}{2} \norm{x}^2 + \tfrac{1}{2} \norm{y}^2 - \ip{x}{y},
\label{eqn: chi} \\
\intertext{so that}
&\ip{\w{f}}{\w{g}} =
\exp \Big( {-} \int \rd s \, \chi (f(s), g(s)) \Big),
\quad f,g \in L^2(\Rplus; \noise). \notag
\end{align}
An operator process $\Lambda_F = \big( \Lambda_F(t) \big)_{t \ges 0}$
is defined, for each operator $F \in B(\init \ot \khat; \init' \ot
\khat)$, through the identity
\[
\ip[\Big]{u \ot \w{f}}{\Lambda_F(t) v \ot \w{g}} =
\int^t_0 \rd s \, \ip[\big]{u \ot \wh{f}(s)}{F ( v \ot \wh{g}(s)}
\ip{\w{f}}{\w{g}}
\]
for all $v \in \init$, $u \in \init'$, $f, g \in \Step$ and
$t \in \Rplus$.  We also use the following notations for an operator
$\phi$ from $\Ov$ to $\Ov \otm B(\khat)$ with domain $\Ov_0$:
\begin{equation}
\label{xygen}
\phi_{x,y} := \big( \id_\Ov \otm\, \om_{\wh{x}, \wh{y}} \big) \circ
\phi -\chi(x,y) \id_\Ov, \quad x,y \in \noise,
\end{equation}
and $\Lambda_\phi$ for the process-valued map $\Lambda_{\phi(\cdot)}$.
For subsets $\Til$ and $\Til'$ of $\noise$ which are total and contain
$0$, a process $k$ on $\Ov$ is said to \emph{satisfy~\eqref{QSDE}
$\Exps_{\Til'}$-weakly for the domain $\init \aot \Exps_\Til$}, if
\[
\ip{u \ot \w{f}}{\big( k_t(a) -a \ot I_\Fock \big) v \ot \w{g}} =
\int^t_0 \rd s \, \ip{u \ot \w{f}}{k_s \big( E^{\wh{f}(s)} \phi(a)
E_{\wh{g}(s)} \big) v \ot \w{g}}
\]
for all $a \in \Ov_0, t \in \Rplus, u \in \init', v \in \init, f \in
\Step_{\Til'}$ and $g \in \Step_{\Til}$. Such a solution is
\emph{weakly regular} if, for all $f \in \Step_{\Til'}$ and $g \in
\Step_{\Til}$, $E^{\w{f}} k_t( \cdot ) E_{\w{g}}$ is bounded with
bounds that are locally uniform in $t$.

Our first result contains a basic criterion for uniqueness of
solutions of QS differential equations, extending those of~\cite{gran}
and~\cite{mother}. Recall that a \emph{pregenerator} of a
$C_0$-semigroup is a closable operator whose closure is a generator.

\begin{thm}
\label{thm: unique}
Let $\phi$ be an operator from $\Ov$ to $\Ov \otm B(\khat)$ with
domain $\Ov_0$, and let $\Til$ and $\Til'$ be total subsets of
$\noise$ containing $0$.
\begin{alist}
\item
Let $k \in \QSC(\Ov, \noise)$ and suppose that, for all $x \in \Til'$
and $y \in \Til$, its $(x,y)$-associated semigroup $\msgp^{x,y}$ is
strongly continuous and has pregenerator $\phi_{x,y}$. Then $k$ is an
$\Exps_{\Til'}$-weak solution of the QS differential equation
\begin{equation}
\label{QSDE}
\rd k_t = k_t \cdot \srd \La_\phi (t), \quad k_0 = \iota^\Ov_\Fock
\end{equation}
on $\Ov_0$ for the domain $\init \aot \Exps_\Til$.
\item
Conversely, suppose that for all $x \in \Til'$ and $y \in \Til$,
$\phi_{x,y}$ is a pregenerator of a $C_0$-semigroup $\msgp^{x,y}$ on
$\Ov$. Then~\eqref{QSDE} has at most one weakly regular
$\Exps_{\Til'}$-weak solution on $\Ov_0$ for the domain
$\init \aot \Exps_\Til$. Moreover, any such solution is a
\tu{(}weak\tu{)} QS cocycle on $\Ov$ whose $(x,y)$-associated
semigroup is $\msgp^{x,y}$, for each $x \in \Til'$ and $y \in \Til$.
\item
Let $k$ be an $\Exps_{\Til'}$-weak solution of~\eqref{QSDE} on $\Ov_0$
for the domain $\init \aot \Exps_\Til$.  Then, for any element $a$ of
$\Ov_0$, the following are equivalent\tu{:}
\begin{rlist}
\item
$k_t(a) = \iota^\Ov_\Fock(a)$ for all $t \in \Rplus$\tu{;}
\item
$\phi(a) = 0$.
\end{rlist}
In particular, if $\Ov$ is an operator system and $1 \in \Ov_0$ then
$k$ is unital if and only if $\phi(1) = 0$.
\end{alist}
\end{thm}

\begin{proof}
(a) This is a straightforward consequence of the semigroup
decomposition of QS cocycles.

(b) Let $f \in \Step_{\Til'}$ and $g \in \Step_{\Til}$, let
$[t_0, t_1[, [t_1, t_2[, \cdots , [t_n, t_{n+1}[$ be common intervals
of constancy of $f$ and $g$, with $t_0 = 0$ and $t_{n+1} = \infty$,
and set $\ups := ( \phi_{f(t), g(t)} )_{t \ges 0}$. Since $\ups$ is
pregenerator-valued, it follows from semigroup theory that the
integral equation
\begin{equation}
\label{eqn: Qta}
Q_t a = a + \int_0^t \rd s \, Q_s \ups_s(a),
\qquad
a \in \Ov_0, t \in \Rplus
\end{equation}
has unique strongly continuous $B(\Ov)$-valued solution $Q$, namely
the `piecewise semigroup evolution' generated by $\ol{\ups} := (
\ol{\phi}_{f(t), g(t)} )_{t \ges 0}$, where $\ol{\phi}_{x,y}$ denotes
the closure of $\phi_{x,y}$; it is given by
\begin{equation}
\label{eqn: Qt=}
Q_t = \msgp^{f(t_0),g(t_0)}_{t_1 - t_0} \circ \cdots \circ
\msgp^{f(t_{i-1}),g(t_{i-1})}_{t_i - t_{i-1}} \circ
\msgp^{f(t_i),g(t_i)}_{t - t_i}
\quad
\text{for } t \in [t_i, t_{i+1}[ \text{ and }  i \in \{ 0, \cdots n \}.
\end{equation}
The result follows since, for any weakly regular $\Exps_{\Til'}$-weak
solution $k$ of~\eqref{QSDE} on $\Ov_0$ for the domain $\init \aot
\Exps_\Til$, the $B(\Ov)$-valued family $\big( E^{\w{f_{[0,t[}}} k_t(
\cdot ) E_{\w{g_{[0,t[}}} \big)_{t \ges 0}$ is strongly continuous and
satisfies~\eqref{eqn: Qta}, and identity~\eqref{eqn: Qt=} then affirms
the (weak) cocycle property of $k$ (through semigroup decomposition).

(c) If~(i) holds then~(ii) follows since, for all $u \in \init', v \in
\init, x \in \Til', y \in \Til$ and $T > 0$,
\[
e^{ -\chi(x,y) T} \ip{u}{E^{\wh{x}} \phi(a) E_{\wh{y}} v} = \lim_{t
\To 0} t^{-1} \ip{u \otimes \w{x_{[0,T[}}}{(k_t(a) -\iota^\Ov_\Fock
(a)) v \otimes \w{y_{[0,T[}}} =0.
\]
The converse is even more straightforward to verify.
\end{proof}

\begin{rems}
(i) Another version of the uniqueness part of this theorem is proved
in~\cite{BWFK}.

(ii) The theorem has a corresponding version appropriate when $\Ov$ is
a von Neumann algebra.
\end{rems}

Specialising to completely bounded operators $\phi$, we next summarise
some quantum stochastic lore (see~\cite{Lgreifswald} and references
therein). Note that, for such $\phi$, each $\phi_{x,y}$ is completely
bounded and therefore generates a cb-norm-continuous semigroup so
Theorem~\ref{thm: unique} applies.

\begin{thm}
\label{CB generation}
Let $\phi \in CB(\Ov; \Ov \otm B(\khat))$. The QS differential
equation~\eqref{QSDE} has a unique weakly regular $\Exps$-weak
solution on the exponential domain $\init \aot \Exps$, denoted
$k^\phi$. Moreover $k^\phi$ is a strong solution and an elementary
\tu{(}weak\tu{)} QS cocycle on $\Ov$ whose $(x,y)$-associated
semigroup has generator $\phi_{x,y}$, for all $x,y\in\noise$.

If the QS cocycle $k^\phi$ is completely bounded then, for any Hilbert
space $\hil$, the cb-process $\big( (k^\phi_t)^{\hil} \big)_{t \ges 0}$
on $\Ov \otm B(\hil)$ equals the elementary QS cocycle $k^\Phi$, where
$\Phi = \phi^\hil$. Moreover it is $\ka$-decomposable for any
orthonormal basis $\ka$ of $\hil$.
\end{thm}

\begin{rem}
The resulting map
\begin{equation}
\label{eqn: QS gen}
\Phi_{\Ov, \noise}:
CB (\Ov; \Ov \otm B(\khat)) \To \El\QSC (\Ov, \noise), \quad \phi
\mapsto k^\phi.
\end{equation}
is injective, and is referred to as the \emph{QS generation map}
(\emph{on $\Ov$, with respect to the noise dimension space $\noise$});
$\phi$ is called the \emph{stochastic generator} of the QS cocycle
$k^\phi$. For $\phi \in CB (\Ov; \Ov \otm B(\khat))$, $\phi^\dagger :=
\bigl( a^* \mapsto \phi(a)^* \bigr) \in CB (\Ov^\dagger; \Ov^\dagger
\otm B(\khat))$, where $\Ov^\dagger$ is the adjoint operator space $\{
a^*: a \in \Ov \}$, and $k^\phi$ has hermitian conjugate QS cocycle
$k^{\phi^\dagger}$. When $\Ov^\dagger = \Ov$, we set
\begin{equation}
\label{eqn: hVk}
\genh(\Ov, \noise):=
\{ \phi \in CB (\Ov; \Ov \otm B(\khat)): \phi^\dagger = \phi \},
\end{equation}
thus (see~\eqref{eqn: QSCh}) $\Phi_{\Ov, \noise} ( \genh(\Ov, \noise) )
\subset \QSC_\h(\Ov,\noise)$. Also, when $\Ov$ is an operator system,
we set
\begin{equation}
\label{eqn: uVk}
\genu(\Ov, \noise):=
\{ \phi \in CB (\Ov; \Ov \otm B(\khat)): \phi(1) = 0 \},
\end{equation}
thus (see~\eqref{eqn: QSCu}), by part~(c) of Theorem~\ref{thm:
unique}, $\Phi_{\Ov, \noise} ( \genu(\Ov, \noise) ) \subset
\QSC_\unital(\Ov, \noise)$.

For a $C^*$-algebra, the QS generation map is shown to restrict to
various bijections of interest in Section~\ref{sec: Mr CPC}.
\end{rem}

We next give a useful invariance principle, whose straightforward
proof follows from the definition of $\msgp^\Ga$
(Proposition~\ref{kGaPGa}), the semigroup decomposition of QS
cocycles, and the fact that, for all $x,y\in\noise$, $\phi_{x,y}$ is
the generator of the $(x,y)$-associated semigroup of the QS cocycle
$k^\phi$.

\begin{propn}
\label{invariance}
Let $k \in \QSCcb(\Ov, \noise)$ and, for a set $\Ix$, let $\Ga: \Ix
\To \noise$ be a map whose range is total and contains $0$. Set
$\hil = l^2(\Ix)$. Then, for a closed subspace $\Ow$ of $\Ov$, the
following are equivalent\tu{:}
\begin{rlist}
\item
$k_t (\Ow) \subset \Ow \otm B(\Fock)$ for all $t \ges 0$.
\item
$\msgp^\Ga_t \big( \Ow \otm B( \hil ) \big) \subset \Ow \otm B(\hil)$
for all $t \ges 0$.
\end{rlist}
If these hold then $k$ defines a QS cocycle on $\Ow$ by restriction.

Further equivalences arise if either $k = k^\phi$ for some map $\phi
\in CB (\Ov; \Ov \otm B(\khat))$, or $\msgp^\Ga$ is norm continuous
with generator $\varphi$. Namely,
\begin{rlist}\setcounter{enumi}{2}
\item
$\phi (\Ow) \subset \Ow \otm B(\khat)$,
\end{rlist}
respectively,
\begin{rlist}\setcounter{enumi}{3}
\item
$\varphi \big( \Ow \otm B(\hil) \big) \subset \Ow \otm B(\hil)$.
\end{rlist}
\end{propn}

\begin{rem}
Suppose that the map $\Ga$ is bounded, $k = k^\phi$ where
$\phi \in CB (\Ov; \Ov \otm B(\khat))$, and $k$ is completely bounded
and locally bounded in cb-norm. Then, by Theorem~\ref{thm: KGagen}
below, $\msgp^\Ga$ is cb-norm continuous.
\end{rem}

\section{Transformations of CB generators}
\label{Section: gen}

In this section we obtain the affine transformation from (completely
bounded) stochastic generators of quantum stochastic cocycles to the
generators of the corresponding global semigroups, and establish its
injectivity in cases of interest. We also examine the question of
inverting the transform in case the global semigroup is determined by
an orthonormal basis of the noise dimension space. A basic tool for
the section is diagonal Weyl processes $W^\Ga$ associated with a
$\noise$-valued map $\Ga$ on an index set (defined in~\eqref{diag
Weyl}).

\begin{lemma}
\label{Wgen}
Let $\Ix$ be a set and $\Ga: \Ix \To \noise$ a map. Set $\hil =
l^2(\Ix)$ and let $(\de^\al)_{\al \in \Ix}$ be its standard
orthonormal basis. Then the diagonal Weyl process $W^\Ga$ is both a
left and a right unitary cocycle. It is elementary if and only if
$\Ga$ is bounded, in which case $W^\Ga$ is a strong solution of the QS
differential equation
\begin{equation}
\label{W HP eqn}
\rd X_t = X_t \srd \La_F (t), \quad X_0 = I_{ \hil \ot \Fock},
\end{equation}
for the coefficient operator $F = F_\Ga \in B( \hil \ot \khat)$ given
by
\[
F_\Ga = \begin{bmatrix} -\tfrac{1}{2} L^*_\Ga L_\Ga & - L^*_\Ga \\
L_\Ga & 0 \end{bmatrix},
\]
in which $L_\Ga \in B(\hil; \hil \ot \noise)$ is the bounded operator
determined by the prescription $\de^\al \mapsto \de^\al \ot
\Ga(\al)$. In this case, $j := \big( T \mapsto (W^\Ga_t)^* (T \ot
I_\Fock) W^\Ga_t \big)_{t \ges 0}$ defines an elementary QS flow $j$
on $B(\hil)$ with stochastic generator
\[
\Theta_\Ga (T) := F^*_\Ga (T \ot I_\khat) + (T \ot I_\khat) F_\Ga +
F^*_\Ga (T \ot \De) F_\Ga.
\]
\end{lemma}

\begin{note}
By a \emph{QS flow} we mean a unital ${}^*$-homomorphic QS cocycle.
\end{note}

\begin{proof}
The time-shift identity and Weyl relations,
\[
\si_t \big( W(f) \big) = W (s_t f) \quad \text{and} \quad W(x_{[0,r[})
W(x_{[r,r+t[}) = W(x_{[0,r+t[}) = W(x_{[r,r+t[}) W(x_{[0,r}),
\]
for $f \in L^2 (\Rplus; \noise)$, $x \in \noise$ and $r,t \in \Rplus$,
imply that the unitary process $W^\Ga$ is a left and right QS operator
cocycle. It follows from Lemma~\ref{lemma: F Gamma} and
Corollary~\ref{cty for op cocycles} that $W^\Ga$ is elementary if and
only if $\Ga$ is bounded. The unique solution $X^F$ of~\eqref{W HP
eqn} is a left QS operator cocycle. To show that $X^F = W^\Ga$ first
note that for any $\al \in \Ix$ and $x \in \noise$
\begin{equation}
\label{Ga ids}
F_\Ga (\de^\al \ot \wh{x}) = \de^\al \ot \begin{pmatrix} -\frac{1}{2}
\norm{\Ga(\al)}^2 -\ip{\Ga(\al)}{x} \\ \Ga(\al) \end{pmatrix}.
\end{equation}
From this one readily obtains
\[
\ip{\de^\al}{(E^{\wh{x}} F_\Ga E_{\wh{y}} -\chi(x,y)) \de^\be} =
\frac{\rd}{\rd t} \Big|_{t=0} \ip{\de^\al \ot \w{x_{[0,t[}}}{ W^\Ga_t
(\de^\be \ot \w{y_{[0,t[}}}
\]
which is enough to show that the associated semigroups of these two
cocycles coincide.

That $j = k^\phi$ for $\phi = \Theta_\Ga$ is a standard result
(see~\cite{Lgreifswald} or \cite{mother}, Theorem~7.4).
\end{proof}

Recall the notation $k^\Ga$ for the associated $\Ga$-cocycle
introduced in Proposition~\ref{kGaPGa}.

\begin{thm}
\label{thm: KGagen}
Let $\phi \in CB (\Ov; \Ov \otm B(\khat))$ and, for a set $\Ix$, let
$\Ga$ be a bounded map from $\Ix$ to $\noise$. Set $\hil = l^2(\Ix)$,
with standard orthonormal basis $(\de^\al)_{\al \in \Ix}$. In terms of
the operators given in \tu{Lemma~\ref{Wgen}}, define maps $\Phi_\Ga
\in CB (\Ov \otm B(\hil); \Ov \otm B(\hil \ot \khat))$ and
$\varphi_\Ga \in CB (\Ov \otm B(\hil))$, by
\begin{align*}
&
\Phi_\Ga := (\id_\Ov \otm \ups_\Ga) \circ \phi^\hil + \id_\Ov \otm
\Theta_\Ga, \text{ and }
\\
&
\varphi_\Ga := (\id_{\Ov \otm B(\hil)} \otm \, \om_{\wh{0}, \wh{0}})
\circ \Phi_\Ga.
\end{align*}
where $\ups_\Ga (A) := \wt{F}_\Ga^* A \wt{F}_\Ga$ for the operator
$\wt{F}_\Ga := I_{\hil \ot \khat} +(I_{\hil} \ot \De) F_\Ga$. Under
the identifications~\eqref{op matrix} and~\eqref{map matrix}, the
following hold\tu{:}
\begin{alist}
\item
$(e^{t \varphi_\Ga})_{t \ges 0}$ is the Schur-action semigroup on $\Ov
\otm B(\hil)$ comprised of associated semigroups of the
\tu{(}weak\tu{)} QS cocycle $k^\phi$\tu{:}
\[
e^{t \varphi_\Ga} = [\msgp_t^{\Ga(\al), \Ga(\be)}]_{\al, \be \in \Ix}
\, \schur \quad \text{ for } t \in \Rplus.
\]
\item
Suppose that the QS cocycle $k^\phi$ is completely bounded with
locally bounded cb-norm and let $\msgp^\Ga$ denote its global
$\Ga$-semigroup. Then, setting $k^{\phi, \Ga} = (k^\phi)^\Ga$,
\[
k^{\phi,\Ga} = k^{\Phi_\Ga} \ \text{ and } \ \msgp^\Ga = (e^{t
\varphi_\Ga})_{t\ges 0}.
\]
\item
Suppose that the range of $\Ga$ is total and contains $0$. Then the
affine-linear map $\phi \mapsto \varphi_\Ga$ is injective.
\end{alist}
\end{thm}

\begin{proof}
Let $B(\init; \init')$ be the ambient full operator space of $\Ov$.
For any $\al \in \Ix$ and $x \in \noise$, using~\eqref{Ga ids},
\begin{equation}
\label{FtildeGa}
\wt{F}_\Ga (\de^\al \ot \wh{x}) =  \de^\al \ot \begin{pmatrix} 1 \\ x
+\Ga(\al) \end{pmatrix},
\end{equation}
and therefore
\begin{equation}
\label{phiGatheta}
\ip{u' \ot \de^\al}{\varphi_\Ga (A) u \ot \de^\be} =
\ip{u'}{\phi_{\Ga(\al), \Ga(\be)} (E^{e_\al} A E_{e_\be}) u},
\end{equation}
for all $A \in \Ov \otm B(\hil)$, $u \in \init$, $u' \in \init'$ and
$\al, \be \in \Ix$, where $\phi_{x,y}$ denotes the semigroup generator
defined in~\eqref{xygen}. Consequently (a) follows; in
particular~\eqref{phiGatheta} shows that $\varphi_\Ga$ has
Schur-action.

To prove~(b) note that, by Theorem~\ref{CB generation}, the weak QS
cocycle $k^{\Phi_\Ga}$ is elementary. By assumption on $k^\phi$ and
Proposition~\ref{kGaPGa}, $k^{\phi,\Ga}$ is a well-defined completely
bounded cocycle whose vacuum-expectation semigroup is $\msgp^\Ga$. In
particular $\msgp^\Ga$ is cb-norm continuous by part~(a), and so
$k^{\phi,\Ga}$ is also elementary by Theorem~\ref{M}. So to verify that
$k^{\phi,\Ga} = k^{\Phi_{\Ga}}$ it suffices to check that their
cb-norm-continuous associated semigroups are the same, which can be
established by showing that
\[
\frac{\rd}{\rd t} \Big|_{t=0} \ip{u' \ot \de^\al \ot \w{x_{[0,t[}}}
{(k^{\Phi_\Ga}_t -k^{\phi,\Ga}_t) (A) u \ot \de^\be \ot \w{y_{[0,t[}}}
= 0
\]
for all $u \in \init$, $u' \in \init'$, $\al, \be \in \Ix$,
$x, y \in \noise$ and $A \in \Ov \otm B(\hil)$. This follows from a
straightforward if somewhat lengthy computation using~\eqref{Ga ids},
\eqref{FtildeGa} and~\eqref{xygen}.

For~(c) let $\phi, \phi' \in CB(\Ov; \Ov \otm B( \hil )$, and suppose
that $\varphi_\Ga = \varphi'_\Ga$.  Then, by~\eqref{xygen}
and~\eqref{phiGatheta}, $E^{\wh{\Ga(\al)}} (\phi -\phi')(a)
E_{\wh{\Ga(\be)}} = 0$ for all $\al, \be \in \Ix$ and $a \in \Ov$.
The totality of $\{\wh{\Ga(\al)}: \al \in \Ix\}$ in $\khat$ gives
$\phi = \phi'$.
\end{proof}

\begin{rem}
The two results above have been proved using the `matrix of cocycles'
idea from Section~5.1 of~\cite{spawn}. For example, the matricial QS
cocycle $k^{\phi,\Ga}$ takes the form
\[
\Big( \big[ (I_{\init'} \ot W(\Ga(\al)_{[0,t[}))^* k^\phi_t (\cdot)
(I_\init \ot W(\Ga(\be)_{[0,t[})) \big]_{\al,\be \in \Ix} \Big)_{t
\ges 0}.
\]
\end{rem}

\emph{For the rest of the section} we specialise to the case where
$\Ga$ is determined by the choice of some orthonormal basis $\eta =
(d_i)_{i \in \Ix_0}$ of $\noise$. Thus $l^2 (\Ix_0)$ is identified
with $\noise$, $l^2 (\Ix)$ is identified with $\khat$, and, as usual,
the convention~\eqref{eta bar defn} is in operation. Operators such as
$F_\Ga$ from above will be rebranded as $F_\eta$, etc. For example,
the operator $L_\eta: \khat \To \khat \ot \noise$ from
Lemma~\ref{Wgen} is now a decapitated version of the Schur isometry
$S_\etaaug$ from Section~\ref{MS+decomp}.

The spaces $CB (\Ov; \Ov \otm B(\khat))$ and $CB_{\etaaugdec} (\Ov
\otm B(\khat))$ can both be viewed as subspaces of $\Mat_\Ix \big(
CB(\Ov) \big)$ through the identifications~\eqref{genkamat}
and~\eqref{map matrix} respectively. For the given choice of basis
$\eta$ the map $\phi \mapsto \varphi_\eta$ extends to a map on
$\Mat_\Ix \big( CB(\Ov) \big)$ which is then manifestly
bijective. Indeed, in terms of their matrix components, the
transformation and its inverse are as follows:
\begin{equation}
\label{cpttransfs}
\begin{aligned}
\varphi^0_0 &= \phi^0_0, \\
\varphi^i_0 &= \phi^0_0 +\phi^i_0 -\tfrac{1}{2} \id_\Ov, \\
\varphi^0_j &= \phi^0_0 +\phi^0_j -\tfrac{1}{2} \id_\Ov, \\
\varphi^i_j &= \phi^0_0 +\phi^i_0 +\phi^0_j +\phi^i_j +(\de^i_j -1)
\id_\Ov,
\end{aligned}
\quad \text{and} \quad
\begin{aligned}
\phi^0_0 &= \varphi^0_0, \\
\phi^i_0 &= \varphi^i_0 -\varphi^0_0 +\tfrac{1}{2} \id_\Ov, \\
\phi^0_j &= \varphi^0_j -\varphi^0_0 +\tfrac{1}{2} \id_\Ov, \\
\phi^i_j &= \varphi^i_j -\varphi^i_0 -\varphi^0_j +\varphi^0_0
-\de^i_j \id_\Ov
\end{aligned}
\end{equation}
for $i, j \in \Ix_0$. The difficult task is to determine if, on
restriction, we have a surjection $CB (\Ov; \Ov \otm B(\khat)) \To
CB_{\etaaugdec} (\Ov \otm B(\khat))$. In Section~\ref{sec: Mr CPC} we
determine the image of various subspaces of interest, but for now we
must make do with a partially defined inverse. This is given in terms
of truncations of $\phi$ and uses operators that are defined in terms
of finite subsets $\Jx_0$ of $\Ix_0$. For such subsets, in the
notation~\eqref{eta bar defn}, we set $\Jx = \{0\} \cup \Jx_0$. The
requisite operators are:
\[
C_{\Jx_0} = \sum_{i \in \Jx_0} \ket{d_i} \in \ket{\noise}, \qquad
C_\Jx = \sum_{\al \in \Jx} \ket{e_\al}
\in \ket{\khat}, \qquad \, Q_{\Jx_0} = \sum_{i
\in \Jx_0} \dyad{d_i}{d_i} \in B(\noise),
\]
and the following operators in $B(\khat)$:
\begin{equation}
\label{L,M and N}
\begin{gathered}
\BoxJ := C_\Jx C^*_\Jx, \qquad \De_\Jx := \begin{bmatrix} 0 & 0 \\ 0 &
Q_{\Jx_0} \end{bmatrix}, \qquad Q_\Jx := \begin{bmatrix} 1 & 0 \\ 0 &
Q_{\Jx_0} \end{bmatrix} = \sum_{\al \in \Jx} \dyad{e_\al}{e_\al}, \\
A_\Jx :=
\begin{bmatrix} 0 & -\frac{1}{2} C^*_{\Jx_0} \\
-\frac{1}{2} C_{\Jx_0} & Q_{\Jx_0} - C_{\Jx_0} C^*_{\Jx_0} \end{bmatrix}
\quad \text{and} \quad
B_\Jx := \begin{bmatrix} 1 & - C^*_{\Jx_0} \\ 0 &
I_\noise \end{bmatrix} = I_\khat -\sum_{i \in \Jx_0} \dyad{e_0}{e_i}.
\end{gathered}
\end{equation}

For each $\Jx_0 \fsubset \Ix_0$ and $T \in B(\init \ot \khat)$ we
define
\[
T_{[\Jx]} := Q_\Jx T Q_\Jx,
\]
so that $T_{[\Jx]}$ is a truncation of $T$ when viewed as a matrix
according to~\eqref{op matrix}. Note that $Q_{\Jx_0}$, $Q_\Jx$ and
$\De_\Jx$ are truncations of the bounded operators $I_\noise$,
$I_\khat$ and $\De$, whereas the other operators defined above are
truncations of sesquilinear forms, \emph{i.e.} the norms are unbounded
as $\Jx_0$ grows.  Recall the left inverse $\Ups_\etaaug$ of the Schur
homomorphism $\Si_\etaaug$ from Section~\ref{MS+decomp}.

\begin{thm}
\label{Z}
Let $\phi \in CB \bigl(\Ov; \Ov \otm B(\khat) \bigr)$ and let $\varphi
\in CB_{\etaaugdec} \big(\Ov \otm B(\khat) \big)$. Then the following
are equivalent\tu{:}
\begin{rlist}
\item
\label{Za}
$\varphi = (\id_{\Ov \otm B(\khat)} \otm\, \om_{\wh{0}, \wh{0}}) \circ
\Phi_\eta$, as defined in \tu{Theorem~\ref{thm: KGagen}}.
\item
\label{Zb}
$\varphi = (\id_\Ov \otm \wt{\Ups}_\eta) \circ \phi^{\khat} + \id_\Ov
\otm\, \Ups'_\eta$, for the maps
\begin{align*}
\wt{\Ups}_\eta &: B(\khat \ot \khat) \To B(\khat), \quad
T \mapsto \wt{S}_\eta^* T \wt{S}_\eta \ \text{ where } \ \wt{S}_\eta
:= S_\etaaug + \De \ot \ket{ \wh{0} }, \ \text{and} \\
\Ups'_\eta &: B(\khat) \To B(\khat), \quad T \mapsto
\Ups_\etaaug (T \ot \De) - \tfrac{1}{2} (T \De +\De T).
\end{align*}
\item
\label{Zc}
For each $a \in \Ov$ and $\Jx_0 \fsubset \Ix_0$, the truncation
$\phi_{[\Jx]}: a \mapsto \phi (a)_{[\Jx]}$ is given by
\begin{equation}
\label{thetatrunc}
\phi_{[\Jx]}: a \mapsto (I_{\init'} \ot B_\Jx)^* \big(\varphi (a \ot
\BoxJ) -a \ot A_\Jx \big) (I_\init \ot B_\Jx).
\end{equation}
\end{rlist}
\end{thm}

\begin{proof}
Set $\varphi_\eta := (\id_{\Ov \otm B(\khat)} \otm\, \om_{\wh{0}, \wh{0}})
\circ \Phi_\eta$. The identities
\[
F_\eta E_{\wh{0}} = (I_\khat \ot \De) S_\etaaug -\tfrac{1}{2} \De
\ot \ket{\wh{0}} \quad \text{and} \quad \wt{F}_\eta E_{\wh{0}} =
\wt{S}_\eta
\]
are verified using the relations~\eqref{Ga ids} and~\eqref{FtildeGa};
the equivalence of~\eqref{Za} and~\eqref{Zb} then follows.

The relations
\[
\Ups_\etaaug (\BoxJ \ot \De) = \De_\Jx \quad \text{and} \quad A_\Jx =
\De_\Jx -\tfrac{1}{2} (\De \BoxJ +\BoxJ \De)
\]
show that
\begin{equation}
\label{UpsdashBox}
\Ups'_\eta (\BoxJ) = A_\Jx.
\end{equation}
Moreover, since $(C^*_\Jx \ot \bra{e_\be}) \wt{S}_\eta B_\Jx =
\bra{e_\be} Q_\Jx$ it follows that
\begin{equation}
\label{BStBox}
B^*_\Jx \wt{S}^*_\eta (\BoxJ \ot T) \wt{S}_\eta B_\Jx = T_{[\Jx]},
\end{equation}
first for $T = \dyad{e_\al}{e_\be}$ and thus for any $T \in B(\khat)$
by linearity and ultraweak continuity.

Suppose that~\eqref{Zb} holds. Then~\eqref{Zc} follows
from~\eqref{UpsdashBox},~\eqref{BStBox} and the totality, with respect
to the weak operator topology, of the set of simple tensors in
$\Ov \otm B(\khat)$.

Finally, suppose that $\varphi$ is such that~\eqref{Zc} holds. We
already know that~\eqref{Zc} holds for $\varphi_\eta$ and so
\begin{align*}
(I_{\init'} \ot B_\Jx)^* \varphi(a \ot \BoxJ) (I_\init \ot B_\Jx) &=
\phi_{[\Jx]} (a) +a \ot B^*_\Jx A_\Jx B_\Jx \\
&= (I_{\init'} \ot B_\Jx)^* \varphi_\eta (a \ot \BoxJ) (I_\init \ot
B_\Jx).
\end{align*}
Invertibility of $B_\Jx$ shows that $\varphi (a \ot \BoxJ) =
\varphi_\eta (a \ot \BoxJ)$ for all $a \in \Ov$ and $\Jx_0 \fsubset
\Ix_0$. Since $\phi$ and $\phi_\eta$ are both $\etaaug$-decomposable,
this is sufficient to deduce that $\varphi = \varphi_\eta$, so (i)
holds.
\end{proof}

\section{Generation of completely positive QS cocycles}
\label{CPC cocs}

In this section we consider the stochastic generation of completely
positive QS cocycles.  Let $\Oe$ be an operator system acting
nondegenerately on $\init$ (so that $1_\Oe = I_\init$), and fix an
orthonormal basis $\eta = (d_i)_{i \in \Ix_0}$ for $\noise$ and, as
usual, set $\Ix := \{ 0 \} \cup \Ix_0$ and $\etaaug = (e_\al)_{\al
\in \Ix}$ with the convention~\eqref{eta bar defn} in operation, as
always.  As well as the truncated operators $\De_\Jx$, $\BoxJ$,
$A_\Jx$ etc.\ from the previous section, we define, for $\Jx_0
\fsubset \Ix_0$ and $t \ges 0$, the operator:
\[
\BoxJt := \begin{bmatrix} 1 \\ e^{-t/2}   C_{\Jx_0} \end{bmatrix}
\begin{bmatrix} 1 \\ e^{-t/2} C_{\Jx_0} \end{bmatrix}^* +(1-e^{-t})
\De_\Jx \in B(\khat).
\]
Thus $\square_{\Jx, 0} = \BoxJ$. We focus on the global semigroup
$\msgp^\eta$ on $\Oe \otm B(\khat)$ determined by a QS cocycle on
$\Oe$ and the given choice of basis. The following summarises
Proposition~6.2 and Theorem~6.4 of~\cite{spawn}, written now in terms
of this single semigroup by means of a re-indexing as used in the
proof of Theorem~6.5 of that paper. (There we worked instead with the
infinite family of semigroups $\{\msgp^\bfx: \bfx \in \Til^n, n \in
\Nat \}$ for a total subset $\Til$ of $\noise$ containing $0$.)
Recall~\eqref{eta bar defn}, the notation $\etaaug$ for the
orthonormal basis $(\wh{d}_\al)_{\al \in \Ix}$ of $\khat$.

\begin{thm}[\cite{spawn}]
\label{recon}
Let $\msgp$ be a semigroup on $\Oe\ot \Cpt(\khat)$, or on $\Oe \otm
B(\khat)$, comprised of completely bounded $\etaaug$-decomposable
maps.
\begin{alist}
\item
Suppose that $\msgp = \msgp^\eta$, the global $\eta$-semigroup
associated to a cocycle $k \in \QSCcb(\Oe, \noise)$. Then the
following hold\tu{:}
\begin{rlist}
\item
$k$ is completely positive if and only if $\msgp$ is.
\item
Suppose that $k$ is completely positive. Then $k$ is contractive if
and only if $\msgp$ satisfies
\begin{equation}
\label{cont+ineq}
\msgp_t (1_\Oe \ot \BoxJ) \les 1_\Oe \ot \BoxJt \text{ for all } \Jx_0
\fsubset \Ix_0 \text{ and } t \ges 0,
\end{equation}
and $k$ is unital if and only if~\eqref{cont+ineq} holds with
equality.
\end{rlist}
\item
Conversely, suppose that $\msgp$ is completely positive and
satisfies~\eqref{cont+ineq}. Then there is a unique cocycle $k \in
\QSCcpc (\Oe, \noise)$ whose global $\eta$-semigroup is $\msgp$.
\end{alist}
\end{thm}

\begin{rems}
In part (a)\,(i) $\Oe$ may alternatively be assumed to be a (not
necessarily unital) $C^*$-algebra. Parts~(a)\,(ii) and~(b) can also be
reformulated to cover the case of nonunital $C^*$-algebras by using
operator intervals. See Theorem~6.7 of \cite{spawn} for details.

Note that no continuity in $t$ is assumed in Theorem~\ref{recon}. The
flexibility as to whether we look at semigroups on $\Oe \ot \Cpt
(\khat)$ or on $\Oe \otm B(\khat)$ follows from Proposition~\ref{new
propn}, however when considering continuity questions one should
heed cautionary results such as Proposition~\ref{auto CB reg}.

Let $\Jx_0 \subset \Lx_0 \fsubset \Ix_0$. If the
inequality~\eqref{cont+ineq} holds for $\Lx = \{ 0 \} \cup \Lx_0$
then, by $\etaaug$-decomposability, it holds for $\Jx$ too. Thus if
$\noise$ is finite dimensional then~\eqref{cont+ineq} may be replaced
by the inequalities
\[
\msgp_t (1_\Oe \ot \square_\Ix) \les 1_\Oe \ot
\square_{\Ix, t}, \quad t \ges 0.
\]
\end{rems}

We next derive an infinitesimal adjunct to part~(a) of
Theorem~\ref{recon}.

\begin{propn}
\label{CPC+char}
Let $\msgp$ be a $C_0$-semigroup on $\Oe \ot \Cpt(\khat)$ consisting
of completely positive, $\etaaug$-decomposable maps, and let $\varphi$
be its generator. Let $\Jx_0 \fsubset \Ix_0$ and assume that
$1_\Oe \ot \BoxJ \in \Dom \varphi$. Then the following are
equivalent\tu{:}
\begin{rlist}
\item
\label{semi ineq}
$\msgp_t (1_\Oe \ot \BoxJ) \les 1_\Oe \ot \BoxJt$ for all $t \in
\Rplus$.
\item
\label{gen ineq}
$\varphi (1_\Oe \ot \BoxJ) \les 1_\Oe \ot A_\Jx$ where $A_\Jx$ is the
operator defined in~\eqref{L,M and N}.
\end{rlist}
Furthermore,~\eqref{semi ineq} holds with equality if and only
if~\eqref{gen ineq} does.
\end{propn}

\begin{proof}
Let us abbreviate $1_\Oe \ot \BoxJ, 1_\Oe \ot \De_\Jx$ and $1_\Oe \ot
\De$, to $ \BoxJ, \De_\Jx$ and $\De$ respectively. In view of the
$\etaaug$-decomposability of $\msgp$,
\[
\wt{\msgp}_t := R_t \msgp_t (\cdot) R_t \text{ where } R_t = \De^\perp
+ e^{t/2} \De,
\]
defines a completely positive, $\eta$-decomposable $C_0$-semigroup
$\wt{\msgp}$ with generator $\wt{\varphi}$ given by
\[
\Dom \wt{\varphi} = \Dom \varphi, \quad \wt{\varphi}(T) = \varphi (T)
+\tfrac{1}{2} (\De T + T \De),
\]
and the proposition is established once it is proved that
\[
\wt{\msgp}_t (\BoxJ) \les \BoxJ +(e^t-1) \De_\Jx \text{ for all } t
\ges 0 \ \Iff \ \wt{\varphi} (\BoxJ) \les \De_\Jx.
\]
Equality holds on the left when $t=0$, so the left to right
implication is obtained by differentiation. Suppose therefore that
$\wt{\varphi}$ satisfies the right-hand inequality. Since the
semigroup $\wt{\msgp}$ is positive and $\eta$-decomposable
\begin{align*}
\wt{\msgp}_t (\BoxJ) &= \BoxJ +\int^t_0 \!\rd t_1 \,
(\wt{\msgp}_{t_1} \circ \wt{\varphi}) (\BoxJ) \\
&\les \BoxJ +\int^t_0 \!\rd t_1 \, \wt{\msgp}_{t_1} (\De_\Jx) = \BoxJ
+\int^t_0 \!\rd t_1 \, \De_\Jx \schur \wt{\msgp}_{t_1} (\BoxJ)
\end{align*}
where $\schur$ denotes the $\eta$-Schur product. Since $\De_\Jx \in
1_\Oe \ot B(\khat)_+$, iteration is sanctioned by
Lemma~\ref{positivity for Schur}, giving
\begin{align*}
\wt{\msgp}_t (\BoxJ) &\les \BoxJ +\int^t_0 \!\rd t_1 \, \De_\Jx \schur
\bigg( \BoxJ +\int^{t_1}_0 \!\rd t_2 \ \De_\Jx \schur \wt{\msgp}_{t_2}
(\BoxJ) \bigg) \\
&= \BoxJ +t \De_\Jx +\int^t_0 \!\rd t_1 \int^{t_1}_0 \!\rd t_2 \
\De_\Jx \schur \wt{\msgp}_{t_2} (\BoxJ).
\end{align*}
Further iteration gives
\[
\wt{\msgp}_t (\BoxJ) \les \BoxJ +\sum^n_{k=1} \frac{t^k}{k!} \De_\Jx
+\ep_n (t),
\]
where
\[
\ep_n (t) = \int^t_0 \!\rd t_1 \int^{t_1}_0 \!\rd t_2 \, \cdots
\int^{t_n}_0 \!\rd t_{n+1} \, \wt{\msgp}_{t_{n+1}} (\De_\Jx).
\]
Since $\ep_n (t) \To 0$ as $n \To \infty$, the result follows, noting
that the case of equality in~\eqref{semi ineq} and~\eqref{gen ineq}
follows by the same argument.
\end{proof}

\begin{rem}
One cannot work with $\msgp$ and $\varphi$ directly since $A_\Jx$ is
not positive.
\end{rem}

Let $\phi$ be an operator from $\Oe$ to $\Oe \otm B(\khat)$ with
domain $E_0$. If $k$ is a completely positive process on $\Oe$
satisfying~\eqref{QSDE} $\Exps_{\Til(\eta)}$-\emph{weakly} for the
domain $\init \aot \Exps_{\Til(\eta)}$ and $1_\Oe \in \Oe_0$ then $k$
is contractive (respectively unital) if and only if $\phi (1_\Oe) \les
0$ (respectively $\phi (1_\Oe) = 0$); the proof is given in
Proposition~5.1 of~\cite{mother}. Moreover, in favourable
circumstances (\emph{e.g.}~\cite{mother}, Theorem~3.1) such a weak
solution $k$ of~\eqref{QSDE} enjoys a semigroup decomposition that, in
conjunction with Proposition~5.1 of~\cite{spawn}, shows that $k$ is a
cocycle. In the case that $\phi$ is unbounded, the missing ingredient
at this stage is some hypothesis that ensures that certain affine
combinations of components of $\phi$ are pregenerators of
$C_0$-semigroups.

We now turn this line of argument around and instead prove existence
of a solution of~\eqref{QSDE} by going via global semigroups and their
associated cocycles. For this we make use of the necessary conditions
on $\phi$ for a CP process weakly satisfying~\eqref{QSDE} to be
contractive.

To do this, first note the following identity for $\Jx_0 \fsubset
\Ix_0$:
\begin{equation}
\label{eqn: 4.2}
\sum_{\al, \be \in \Jx} \chi( d_\al, d_\be ) \dyad{e_\al}{e_\be} = -A_\Jx,
\end{equation}
where the map $\chi$ is defined in~\eqref{eqn: chi},
$\Jx = \{ 0 \} \cup \Jx_0$ as usual, and $A_\Jx$ is the operator
defined in~\eqref{L,M and N}. Also, we further specialise the notation
$\msgp^\Ga$ for associated semigroups and cocycles as follows: for
$\Jx_0 \fsubset \Ix_0$, $\msgp^{[\Jx]}$ denotes the semigroup
associated to the map $\Ga: \Jx \To \noise$, $\al \mapsto d_\al$. Note
that $\msgp^{[\Jx]}$ has Schur-action when we make the identification
$\Oe \otm B(l^2(\Jx)) \cong \Mat_\Jx (\Oe)$. Recall the notations
$\Til(\eta)$ and $\etaaug$ introduced in~\eqref{eqn: 2.3}

\begin{thm}
\label{thm: ABC}
Let $\phi$ be an operator from $\Oe$ to $\Oe \otm B(\khat)$ with
domain $\Oe_0$ satisfying the following conditions\tu{:}
\begin{itemize}
\item[(1)]
$1_\Oe \in \Oe_0$ and $\phi(1_\Oe) \les 0$.
\item[(2)]
For all $\al, \be \in \Ix$, the operator $\phi_{d_\al, d_\be}$
on $\Oe$ with domain $\Oe_0$ \tu{(}defined in~\eqref{xygen}\tu{)} is a
pregenerator of a $C_0$-semigroup $\msgp^{(\al, \be)}$ on $\Oe$.
\item[(3)]
For all $\Jx_0 \fsubset \Ix_0$, the semigroup $\msgp^{[\Jx]}$
\tu{(}defined above\tu{)} is completely positive.
\end{itemize}
Then there is a unique process $k \in \QSCb (\Oe, \noise)$ which is
locally norm bounded and such that for all $\al, \be \in \Ix$,
$\msgp^{(\al, \be)}$ is its $(d_\al, d_\be)$-associated
semigroup. Moreover,
\begin{alist}
\item
$k$ is completely positive and contractive, and pointwise
$\Fock$-ultraweakly continuous. Furthermore it is unital if and only
if $\phi (1_\Oe) = 0$.
\item
$k$ is the unique weakly regular $\Exps_{\Til(\eta)}$-weak solution of
the QS differential equation~\eqref{QSDE} on $\Oe_0$ for the domain
$\init \aot \Exps_{\Til(\eta)}$.
\item
$k$ strongly satisfies~\eqref{QSDE} on $\Oe_0$ for the domain $\init
\aot \Exps_{\Til(\eta)}$, under the following further condition\tu{:}
\begin{itemize}
\item[(4)]
$\init$ and $\noise$ are separable.
\end{itemize}
\end{alist}
\end{thm}

\begin{proof}
The uniqueness part is immediate since every bounded QS cocycle
is uniquely determined by its $(x,y)$-associated semigroups for $x$
and $y$ running through any total subset of $\noise$ containing
$0$. We therefore address existence.

First note that, for all $\al \in \Ix$, $\phi_{d_\al, d_\al} =
E^{\wh{d_\al}} \phi( \cdot ) E_{\wh{d_\al}}$ since $\chi (d_\al,
d_\al) = 0$, and so $\phi_{d_\al, d_\al} (1_\Oe) \les 0$
by~(1). Moreover, since $\msgp^{[\Jx]}$ has Schur-action it follows
that each $\msgp^{(\al, \al)}$ is positive, and hence is contractive
since
\[
0 \les \msgp^{(\al, \al)}_t (1_\Oe) = 1_\Oe +\int^t_0 \rd s \,
\msgp^{(\al, \al)}_s \big( \phi_{d_\al, d_\al} (1_\Oe) \big) \les
1_\Oe.
\]
Thus $\msgp^{[\Jx]}_t (1_{\Mat_\Jx (\Oe)}) \les 1_{\Mat_\Jx (\Oe)}$
for each $\Jx_0 \fsubset \Ix_0$, showing that the completely positive
semigroup $\msgp^{[\Jx]}$ is contractive.

Linear extension of the prescription
\[
\msgp^0_t: a \ot \dyad{e_\al}{e_\be} \mapsto \msgp^{(\al, \be)}_t (a)
\ot \dyad{e_\al}{e_\be}
\]
defines an operator $\msgp^0_t$ on the operator system $\Oe \ot \Cpt
(\khat)$ with dense domain $\Oe_0 \aot B_{00} (\etaaug)$, where
$B_{00} (\etaaug) := \Lin \big\{ \dyad{e_\al}{e_\be}: \al, \be \in
\Ix \big\} $, which inherits complete positivity and contractivity
from the family $( \msgp^{[\Jx]}_t )_{\Jx_0 \fsubset \Ix_0}$. Let
$\msgp_t$ denote its continuous extension to all of $\Oe \ot
\Cpt(\khat)$. Then $\msgp := (\msgp_t)_{t \ges 0}$ is a completely
positive and contractive $\etaaug$-decomposable $C_0$-semigroup on
$\Oe \ot \Cpt(\khat)$. Its generator $\psi$ satisfies
\begin{align*}
&\Dom \psi \supset \Oe_0 \aot B_{00} (\etaaug), \ \text{ and} \\
&\psi (a \ot \dyad{e_\al}{e_\be}) = \phi_{d_\al, d_\be} (a) \ot
\dyad{e_\al}{e_\be}, \qquad a \in \Oe_0, \, \al, \be \in \Ix.
\end{align*}
In particular, using identity~\eqref{eqn: 4.2}, for all $\Jx_0
\fsubset \Ix_0$
\begin{align*}
\psi (1_\Oe \ot \BoxJ) &= \sum_{\al, \be \in \Jx} \phi_{d_\al, d_\be}
(1_\Oe) \ot \dyad{e_\al}{e_\be} \\
&= \sum_{\al, \be \in \Jx} \big( E^{\wh{d_\al}} \phi (1_\Oe)
E_{\wh{d_\be}} -\chi (d_\al, d_\be) 1_\Oe \big) \ot
\dyad{e_\al}{e_\be} \\
&= D^*_\Jx \, \phi(1_\Oe) D_\Jx +1_\Oe \ot A_\Jx \les 1_\Oe \ot A_\Jx,
\end{align*}
where $D_\Jx = I_\init \ot \sum_{\al \in \Jx}
\dyad{\wh{d_\al}}{e_\al}$. Therefore, by Proposition~\ref{CPC+char},
\[
\msgp_t (1_\Oe \ot \BoxJ) \les 1_\Oe \ot \BoxJt \qquad \Jx_0 \fsubset
\Ix_0, \, t \in \Rplus,
\]
and so, by Theorem~\ref{recon}, there is a unique cocycle $k \in
\QSCcpc (\Oe, \noise)$ whose global $\eta$-semigroup on $\Oe \ot
\Cpt(\khat)$ is $\msgp$. In particular the $(d_\al, d_\be)$-associated
semigroup of $k$ is $\msgp^{(\al, \be)}$.

(a) We have already proved that $k$ is completely positive and
contractive. The rest follows from Theorems~\ref{M}
and~\ref{thm: unique}.

(b) This follows from Theorem~\ref{thm: unique}.

(c) Let $K$ denote the QS cocycle $\bigl( (k_t)^{\khat} \bigr)_{t \ges
0}$ on $\Oe \otm B(\khat)$. To see that our weak solution is actually
a strong solution involves showing that, for each $a \in \Oe_0$, the
process $\bigl( K_t (\phi(a)) \bigr)_{t \ges 0}$ is QS integrable on
$(\init \ot \khat) \aot \Exps$ which, because $k$ is completely
contractive, amounts to proving strong measurability of the
vector-valued process $\bigl( K_t (\phi(a)) \xi \bigr)_{t \ges 0}$ for
each $\xi \in (\init \ot \khat) \aot \Exps$; under assumption~(4) this
follows from Pettis' Theorem.
\end{proof}

\begin{rems}
(i) In Section~\ref{sec: structure} we give sufficient conditions on
$\phi$ for assumption~(3) to hold.

(ii) The uniqueness here extends that of the QS cocycle generated by a
mapping in $CB (\Ov; \Ov \otm B(\khat))$ for an operator space $\Ov$
(see Theorem~\ref{CB generation}).

(iii) An alternative condition to~(4), which implies that a weak
solution is in fact a strong solution, is that $\Oe$ is a von Neumann
algebra so that Proposition~\ref{ptwise SOT} applies. However, as
noted in the remark following Proposition~\ref{O}, strong continuity
is an inappropriate assumption in those circumstances.
\end{rems}

\section{Completely positive elementary QS cocycles}
\label{sec: Mr CPC}

For this section let $\Al$ be a $C^*$-algebra acting nondegenerately
on $\init$. The main goal of the section is to reveal the structure of
the stochastic generator of a completely positive elementary QS
cocycle. To this end, for $R \in B(\init; \init \ot \khat)$ and
$\phi \in CB( \Al; \Al \otm B(\khat))$, define maps
\begin{align}
&&& \psi_R \in CB( B(\init); B(\init \ot \khat) ), && T \mapsto R T
E^{\wh{0}} +E_{\wh{0}} T R^* -T \ot \De, \ \text{ and } &&
\label{psiR} \\
&&& \chi_{\phi, R} \in CB (\Al; B(\init \ot \khat)), && \chi_{\phi, R}
:= \phi -\psi_R. && \notag
\end{align}
Set
\[
\gencp(\Al, \noise) := \big\{ \phi \in CB (\Al; \Al \otm B(\khat)):
\chi_{\phi,R} \text{ is CP for some } R \in B(\init; \init \ot \khat)
\big\}.
\]
Note that $\gencp(\Al, \noise) \subset \genh(\Al, \noise)$.

\begin{propn}
\label{CPcoc}
Let $\phi \in \gencp (\Al, \noise)$, suppose that $\Al$ is unital and
that the QS cocycle $k^\phi$ is completely bounded. Then $k^\phi$ is
completely positive.
\end{propn}

\begin{proof}
Fix an orthonormal basis $\eta$ for $\noise$. Let $\varphi_\eta$ be
the generator of the global $\eta$-semigroup $\msgp^\eta$ of the
elementary cocycle $k^\phi$, as given by Theorem~\ref{Z} whose
notations we continue to use. Choose $R$ such that $\chi_{\phi, R}$ is
completely positive, then
\[
\phi^\khat (A) = \chi^\khat_{\phi, R} (A) +\Pi (R \ot I_\khat) A
(I_{\init \ot \khat} \ot \bra{\wh{0}}) +(I_{\init \ot \khat} \ot
\ket{\wh{0}}) A (R^* \ot I_\khat) \Pi -A \ot \De,
\]
where $\Pi$ is the unitary tensor flip on $\init \ot \khat \ot \khat$
exchanging the two copies of $\khat$. The identities
\begin{align*}
(\id_\Al \otm \wt{\Ups}_\eta) (A \ot \De) &= (\id_\Al \otm
\Ups_\etaaug) (A \ot \De) \ \text{ and} \\
\big( I_{\init \ot \khat} \ot \bra{\wh{0}} \big) (I_\init \ot
\wt{S}_\eta) &= I_{\init \ot \khat}
\end{align*}
therefore imply that
\[
\varphi_\eta (A) = (I_\init \ot \wt{S}_\eta)^* \chi_{\phi,R}^{\khat}
(A) (I_\init \ot \wt{S}_\eta) + \wt{R}A + A\wt{R}^*
\]
for the operator $\wt{R} := (R \ot I_\khat) \Pi (I_\init \ot
\wt{S}_\eta) -\tfrac{1}{2} I_\init \ot \De$. Applying Stinespring's
Theorem to the map $\chi_{\phi,R}$, and then applying
Theorem~\ref{decomp CE} to $\varphi_\eta$, it follows that
$\msgp^\eta$ is completely positive. Thus $k^\phi$ is completely
positive by Theorem~\ref{recon}.
\end{proof}

\begin{rem}
In case $\phi \in \gencp (\Al, \noise)$ but $k^\phi$ is not completely
bounded, the global semigroup $\msgp^\eta$ still exists and is
completely positive. This implies that $\msgp^{\bfx}$ is completely
positive for all $\bfx \in \bigcup_{n \in \Nat} \noise^n$, which still
implies that the (weak) QS cocycle $k^\phi$ is completely positive,
but now in the following sense: for all $n \in \Nat$, $A = [a^i_j] \in
\Mat_n (\Al)_+$ and $\zeta \in (\init \aot \Exps)^n$,
\[
\sum_{i,j = 1}^n \ip{\zeta^i}{k^\phi_t (a^i_j) \zeta^j} \ges 0.
\]
\end{rem}

For $\Al$ unital, we define the following subsets of $\gencp (\Al,
\noise)$:
\begin{align*}
&
\gencpbeta (\Al, \noise)  := \big\{ \phi \in \gencp (\Al,\noise):
\phi(1_\Al) \les \be \De^\perp \big\},
\\
&
\gencpc (\Al, \noise) := \gencpzero (\Al, \noise),
\\
&
\gencpqc (\Al, \noise) := \bigcup_{\be \in \Real} \gencpbeta (\Al,
\noise), \qquad
\text{and}
\\
&
\gencpu (\Al, \noise) := \gencp (\Al, \noise) \cap \genu (\Al,
\noise).
\end{align*}

\begin{rems}
For $\be \in \Real$, the prescriptions
\[
k \mapsto \big( e^{-\be t} k_t \big)_{t \ges 0} \quad \text{and} \quad
\phi \mapsto \phi -\be \de^\perp,
\]
in which $\de^\perp \in CB( \Al; \Al \otm B(\khat) )$ denotes the map
$a \mapsto a \ot \De^\perp$, define bijections of $\QSC (\Al, \noise)$
and of $CB (\Al; \Al \otm B(\khat))$ respectively. For a QS cocycle
$k$ and real number $\beta$, the QS cocycle $(e^{-\be t} k_t)_{t \ges
0}$ is elementary if and only if $k$ is, and the bijections are
compatible, in the sense that, for $\phi \in CB (\Al; \Al \otm
B(\khat))$,
\[
\big( e^{-\be t} k^\phi_t \big)_{t \ges 0} = k^{\phi -\be \de^\perp}.
\]
This is easily verified by appealing to uniqueness of weak solutions
to the QS differential equation~\eqref{QSDE}, in which $\phi -\be
\de^\perp$ replaces $\phi$, or by making use of the semigroup
decomposition of such solutions (\cite{mother}, Theorem~3.1). The
remarks so far apply equally to stochastically generated cocycles on
\emph{any} operator space $\Ov$.

The above bijections restrict to bijections
\[
\QSCcpbeta (\Al, \noise) \To \QSCcpc (\Al, \noise) \quad
\text{and} \quad \gencpbeta (\Al, \noise) \To \gencpc (\Al,
\noise),
\]
noting that, in the notation~\eqref{psiR}, $\psi_R -\be \de^\perp =
\psi_{R'}$ for $R' = R -\frac{1}{2} \be E_{\wh{0}}$.
\end{rems}

For the proof of the next result, we require the following easily
verified identities associated with an orthonormal basis $\eta =
(d_i)_{i \in \Ix_0}$ for $\noise$ and the operators defined
in~\eqref{L,M and N}: for $\Jx_0 \fsubset \Ix_0$ and Hilbert spaces
$\Kil$ and $\Kil'$, setting $\Ix := \{0\} \cup \Ix_0$, $\Jx := \{0\}
\cup \Jx_0$ as usual and $\etaaug = (e_\al)_{\al \in \Ix}$,
\begin{align}
& B^*_\Jx C_\Jx = \ket{\wh{0}} \label{N*LN=Mcopy} \\
& B_\Jx^* A_\Jx B_\Jx = \De_\Jx - \tfrac{1}{2} \big( \ket{\wh{0}}
C_\Jx^* \De B_\Jx +B_\Jx^* \De C_\Jx \bra{\wh{0}} \big), \text{
and} \label{phitrunccopy} \\
& X \ot \BoxJ = (I_{\Kil'} \ot C_\Jx) X (I_{\Kil} \ot C_\Jx)^*, \text{
for } X \in B (\Kil; \Kil'). \label{Rtrunccopy}
\end{align}
For a supplementary Hilbert space $\Hil$, define the transformation
\begin{align*}
&\al_{\Jx, \Hil}: B(\Hil; \init) \vot \Diag_\etaaug (\khat) \To
B(\Hil; \init) \vot \ket{\khat} = B(\Hil; \init \oplus (\init \ot
\noise)), \\
&X = \begin{bmatrix} x^0 & 0 \\ 0 & X^1 \end{bmatrix} \mapsto (I_\init
\ot B_\Jx)^* X (I_\init \ot C_\Jx) = \begin{bmatrix} x^0 \\ X^1
(I_\Hil \ot C_{\Jx_0}) -(I_\init \ot C_{\Jx_0}) x^0 \end{bmatrix},
\end{align*}
which enjoys the following properties:
\begin{itemize}
\item[(1)]
For any operator space $\Ow$ in $B(\Hil; \init)$,
\[
\al_{\Jx, \Hil} \bigl(\Ow \otm \Diag_\etaaug (\khat) \bigl)
\subset \Ow \otm \ket{\khat}.
\]
\item[(2)]
If $\Jx_0 \subset \Lx_0 \fsubset \Ix_0$ then
\[
Q_\Jx \, \al_{\Lx, \Hil} (\cdot) = \al_{\Jx, \Hil}
\]
since $Q_{\Jx_0} C_{\Lx_0} = C_{\Lx_0}$ and $Q_{\Jx_0} X^1 = X^1
Q_{\Jx_0}$ for all $X^1 \in B(\Hil; \init) \vot \Diag_\eta (\noise)$.
\end{itemize}

Recall the QS generation map defined in~\eqref{eqn: QS gen}.

\begin{propn}
\label{unital qcont cocs}
Suppose that $\Al$ is unital, and let $\be \in \Real$. Then
$\Phi_{\Al, \noise}$ restricts to bijections
\begin{align*}
&
\gencpbeta (\Al, \noise) \To \El\QSCcpbeta (\Al, \noise), \text{ and}
\\
&
\gencpu (\Al, \noise) \To \El\QSCcpu (\Al, \noise).
\end{align*}
Moreover, for $\phi \in \gencpqc (\Al, \noise)$, there is some $R \in
\Al'' \, \vot \ket{\khat}$ for which $\chi_{\phi,R}$ is completely
positive.
\end{propn}

\begin{proof}
By the above remarks we may assume without loss that $\be =0$ and so
it is enough to show that $\Phi_{\Al, \noise}$ restricts to a
bijection $\gencpc (\Al, \noise) \To \El\QSCcpc (\Al, \noise)$. The
second part then follows from the first since for \emph{any}
$\phi \in CB( \Al; \Al \otm B(\khat))$, $k^\phi$ is unital if and only
if $\phi (1_\Al) =0$ by Theorems~\ref{thm: unique} and~\ref{CB
generation}.

Fix an orthonormal basis $\eta$ for $\noise$, and let $\phi \in
\gencpc (\Al, \noise)$. Choose $R \in B(\init; \init \ot \khat)$ such
that $\chi_{\phi, R}$ is completely positive. Then $\varphi_\eta$
(defined by Theorem~\ref{thm: KGagen}) is the generator of a
completely positive semigroup, as shown in the proof of
Proposition~\ref{CPcoc}. Now $\phi^\khat (1_\Al \otm \BoxJ) \les 0$
since $\phi (1_\Al) \les 0$, so Theorem~\ref{Z} and the
identity~\eqref{UpsdashBox} imply that
\[
\varphi_\eta (1_\Al\ot \BoxJ) \les 1_\Al \ot \Ups_\eta' (\BoxJ) =
1_\Al \ot A_\Jx.
\]
Therefore, by Proposition~\ref{CPC+char} and Theorem~\ref{recon},
there is a completely positive QS contraction cocycle $k$ whose global
$\eta$-semigroup $\msgp^\eta$ is $\eta$-decomposable and has generator
$\varphi_\eta$, in particular $\msgp^\eta$ is cb-norm continuous.
Moreover, by Theorem~\ref{thm: KGagen}, the component semigroups of
$\msgp^\eta$ are associated semigroups of the QS cocycle $k^\phi$, and
so, by the semigroup decomposition of QS cocycles, $k = k^\phi$. It
follows that $\msgp^{0,0}$, the vacuum-expectation semigroup of
$k^\phi$, is cb-norm continuous and so \mbox{Theorem~\ref{M}} implies
that the cocycle $k^\phi$ is elementary, as required.

Conversely, suppose that $k \in \El\QSCcpc (\Al, \noise)$. Then
$\msgp^{0,0}$, the vacuum-expectation semigroup of $k$, is cb-norm
continuous. Theorem~\ref{M} therefore implies that the global
$\eta$-semigroup of $k$ is cb-norm continuous; let $\varphi \in
CB_{\etaaugdec} (\Al \otm B(\khat))$ be its generator. By
Theorem~\ref{decomp CE} there is a unital representation $(\pi, \Kil)$
of $\Al$ and $\etaaug$-diagonal operators $T \in B(\Kil; \init) \vot
\Diag_\etaaug (\khat)$ and $N \in \Al'' \, \vot \Diag_\etaaug (\khat)$
such that $\varphi (A) =T (\pi \otm \id_{B(\khat)}) (A) T^* + NA +
AN^*$.

Our goal now is to show that $\varphi$ is of the form $\varphi_\eta$,
as defined in Theorem~\ref{thm: KGagen}, for a map $\phi \in CB(\Al;
\Al \otm B(\khat))$. Taking our cue from Theorem~\ref{Z} define, for
each $\Jx_0 \fsubset \Ix_0$, the map
\[
\phi_\Jx: \Al \To \Al \otm B(\khat), \quad a \mapsto (I_\init \ot
B_\Jx)^* \big( \varphi (a \ot \BoxJ) -a \ot A_\Jx \big) (I_\init \ot
B_\Jx).
\]
Applying identities~\eqref{N*LN=Mcopy}--\eqref{Rtrunccopy} to our
formula for $\varphi$ then gives
\[
\phi_\Jx (a) = V_\Jx \pi (a) V^*_\Jx +R_\Jx a E^{\wh{0}} +E_{\wh{0}} a
R^*_\Jx -a \ot \De_\Jx
\]
for the operators
\[
V_\Jx = \al_{\Jx, \Kil} (T) \in B(\Kil; \init \ot \khat) \quad
\text{and} \quad R_\Jx = \al_{\Jx, \init} (N +\tfrac{1}{2} I_\init \ot
\De) \in B(\init) \, \vot \ket{\khat}.
\]
Property~(1) of the transformations $\al_{\Jx, \Hil}$ implies that
$R_\Jx \in \Al'' \vot \ket{\khat}$ for $\Jx_0 \fsubset
\Ix_0$. Property~(2) entails the following compatibility: for $\Jx_0
\subset \Lx_0 \fsubset \Ix_0$,
\begin{gather*}
V_\Jx = (I_\init \ot Q_\Jx) V_\Lx, \quad R_\Jx = (I_\init \ot Q_\Jx)
R_\Lx \quad \text{and so} \\
\phi_\Jx = (I_\init \ot Q_\Jx) \phi_\Lx (\, \cdot \,) (I_\init \ot
Q_\Jx).
\end{gather*}

Up to this point we have used the complete positivity of $k$ and local
boundedness of its cb-norms, but not its contractivity in any explicit
way. Now, by Theorem~\ref{recon} and Proposition~\ref{CPC+char},
$\phi_\Jx (1_\Al) \les 0$ for each $\Jx_0 \fsubset \Ix_0$, in other
words
\[
V_\Jx V^*_\Jx \les -R_\Jx E^{\wh{0}} -E_{\wh{0}} R^*_\Jx +I_\init \ot \De_\Jx.
\]
Thus, in terms of the block decompositions
\[
V_\Jx = \begin{bmatrix} v^0 \\ V^1_{\Jx_0} \end{bmatrix} \quad
\text{and} \quad R_\Jx = \begin{bmatrix} r^0 \\
R^1_{\Jx_0} \end{bmatrix},
\]
we have the inequalities
\[
0 \les \begin{bmatrix} v^0 {v^0}^* & v^0 {V^1_{\Jx_0}}^* \\
V^1_{\Jx_0} {v^0}^* & V^1_{\Jx_0} {V^1_{\Jx_0}}^* \end{bmatrix}
\les \begin{bmatrix} -r^0 -{r^0}^* & - {R^1_{\Jx_0}}^* \\ -R^1_{\Jx_0}
& Q_{\Jx_0} \end{bmatrix}.
\]
Since $\norm{Q_{\Jx_0}} \les 1$, it follows that for each $\Jx_0
\fsubset \Ix_0$ we have $\norm{V^1_{\Jx_0}} \les 1$ and
$\norm{R^1_{\Jx_0}} \les \norm{r^0+{r^0}^*}^{1/2}$. Compatibility now
implies that the nets $(V^1_{\Jx_0})_{\Jx_0 \fsubset \Ix_0}$ and
$(R^1_{\Jx_0})_{\Jx_0 \fsubset \Ix_0}$ converge in the strong operator
topology to column operators $V^1 \in B(\init; \Kil) \vot
\ket{\noise}$ and $R^1 \in \Al'' \, \vot \ket{\noise}$ respectively,
with $V^1_{\Jx_0}$ being the $\Jx_0$-truncation of $V^1$, and
similarly for $R^1_{\Jx_0}$. Therefore, in terms of the operators
\[
V := \begin{bmatrix} v^0 \\ V^1 \end{bmatrix}
\ \text{ and } \
R := \begin{bmatrix} r^0 \\ R^1 \end{bmatrix},
\]
$\phi_\Jx$ is the $\Jx$-truncation of the map $\phi := \chi +\psi_R$,
where $\chi$ is the completely positive map $a \mapsto V \pi (a) V^*$.

Now $\phi$ is completely bounded, and is the pointwise limit of the
net $(\phi_\Jx)_{\Jx_0 \fsubset \Ix_0}$ in the weak operator
topology. Since $E^{e_\al} \phi (a) E_{e_\be} = E^{e_\al} \phi_\Jx (a)
E_{e_\be}$ as soon as $\Jx$ contains $\al$ and $\be$, $\phi$ is $\Al
\otm B(\khat)$-valued. Moreover
\[
R = \text{w.o.-}\lim R_\Jx \in \Al'' \, \vot \ket{\khat}.
\]
Theorems~\ref{thm: KGagen} and~\ref{Z} therefore imply that $k =
k^\phi$. Finally, since
\[
\phi (1_\Al) = \text{w.o.-}\lim \phi_\Jx (1_\Al) \les 0,
\]
$\phi \in \gencpc (\Al, \noise)$, as required.
\end{proof}

We now turn to the case of \emph{nonunital} $C^*$-algebras. Let $\Al$
be such an algebra acting non\-de\-generately on $\init$, and denote its
unitisation by
\[
\uAl = \{a +\la I_\init: a \in \Al, \la \in \Comp\}.
\]
For each $\phi \in CB (\Al; \Al \otm B(\khat))$ and $\be \in \Real$
define an extension by
\[
\phi_\be: \uAl \To \uAl \otm B(\khat), \qquad a + \la I_\init \mapsto
\phi(a) + \be \la \De^\perp,
\]
and set
\begin{align*}
&
\gencpbeta (\Al, \noise) := \big\{ \phi \in CB(\Al; \Al \otm
B(\khat)): \phi_\be \in \gencpbeta (\uAl, \noise) \big\},
\\
&
\gencpc (\Al, \noise) := \gencpzero (\Al, \noise), \quad
\text{and}
\\
&
\gencpqc (\Al, \noise) :=
\bigcup_{\be \in \Real} \gencpbeta (\Al, \noise).
\end{align*}
Thus $\phi \in \gencpbeta (\Al, \noise)$ if $\phi \in CB(\Al; \Al \otm
B(\khat))$ and $\chi_{\phi_\be, R} \in CB (\uAl; B(\init \ot \khat))$
is completely positive for some choice of $R \in B(\init; \init \ot
\khat)$.

\begin{propn}
\label{unitising cocs}
Suppose that $\Al$ is nonunital, and let $k \in \QSCcpbeta
(\Al,\noise)$. Then the prescription
\[
\bek_t (a +\la I_\init) = k_t(a) +e^{\be t} \la I_{\init \ot \Fock},
\quad a \in \Al, \la \in \Comp, t \in \Real,
\]
defines a cocycle $\bek \in \QSCcpbeta (\uAl, \noise)$. Moreover, if
$k$ is elementary then so is $\bek$, and the stochastic generator of
$\bek$ \tu{(}according to \tu{Proposition~\ref{unital qcont cocs})} is
of the form $\phi_\be$ for a map $\phi \in \gencpbeta (\Al, \noise)$.
\end{propn}

\begin{proof}
By the remarks preceding Proposition~\ref{unital qcont cocs}, we may
suppose without loss that $\be = 0$. Proposition~2.3, and (the proof
of) Theorem~6.7 of~\cite{spawn} imply that the given prescription
defines a cocycle $\zk \in \QSCcpu (\uAl, \noise)$. For $x, y \in
\noise$, the $(x,y)$-associated semigroup of the QS cocycles $\zk$ and
$k$ are related by
\[
\zmsgp^{x,y}_t (a +\la I_\init) = \msgp^{x,y}_t (a) +\la \exp (-t
\chi(x,y)) I_\init, \quad a \in \Al, \la \in \Comp.
\]
Suppose now that $k$ is elementary. Since the map $\uAl \To \Comp$,
$a + \la I_\init \mapsto \la$ is ${}^*$-homomorphic and thus
contractive, we have $\norm{a} \les 2 \norm{a +\la I_\init}$. Thus the
semigroup $\zmsgp^{x,y}$ inherits norm continuity from
$\msgp^{x,y}$. Continuity with respect to the cb-norm follows
similarly, thus $\zk$ is elementary. Therefore, by
Proposition~\ref{unital qcont cocs}, $\zk = k^{\phi'}$ for some $\phi'
\in \gencpu (\uAl, \noise)$. In view of the invariance principle
(Proposition~\ref{invariance}), $\phi'$ induces a map $\phi \in
CB(\Al; \Al \otm B(\khat))$ by restriction. In particular, $\phi' (a +
\la I_\init) = \phi (a)$ for all $a \in \Al$ and $\la \in \Comp$, so
$\phi' = \phi_\be$ for $\be = 0$. Thus $\phi \in \gencpc (\Al,
\noise)$ and the result follows.
\end{proof}

\begin{thm}
\label{thm: nonu}
For any $C^*$-algebra $\Al$, the QS generation map $\Phi_{\Al, \noise}$
restricts to a bijection
\begin{equation}
\label{eqn: bij}
\gencpqc (\Al, \noise) \To \El\QSCcpqc (\Al, \noise)
\end{equation}
and, for $\phi \in \gencpqc (\Al, \noise)$, there is some
$R \in \Al'' \vot \ket{\khat}$ for which $\chi_{\phi, R}$ is
completely positive. In more detail, the following bijections are
incorporated in~\eqref{eqn: bij}\tu{:}
\begin{align*}
&
\gencpbeta (\Al, \noise) \To \El\QSCcpbeta (\Al, \noise)
\quad  \text{for each } \be \in \Real, \text{ in particular,}
\\
&
\gencpc (\Al, \noise) \To \El\QSCcpc (\Al, \noise) \quad
\text{and, if $\Al$ is unital, also}
\\
&
\gencpu (\Al, \noise) \To \El\QSCcpu (\Al, \noise).
\end{align*}
\end{thm}

\begin{proof}
Courtesy of Propositions~\ref{unital qcont cocs} and~\ref{unitising
cocs}, it only remains to assume that $\Al$ is nonunital and show that
the image of $\gencpqc (\Al, \noise)$ is contained in $\El\QSCcpqc
(\Al, \noise)$.  Again by the remark preceding Proposition~\ref{unital
qcont cocs}, it suffices to verify that the image of $\gencpc(\Al,
\noise)$ is contained in $\El\QSCcpc(\Al,\noise)$. Let $\phi \in
\gencpc (\Al, \noise)$, and set $k' := k^{\phi'}$ where $\phi' :=
\phi_\be$ for $\be =0$. Thus $k' \in \El\QSCcpu (\uAl, \noise)$, by
Proposition~\ref{unital qcont cocs}. Since $\phi'$ extends $\phi$,
Proposition~\ref{invariance} implies that $k'$ restricts to a cocycle
$k \in \QSCcpc (\Al, \noise)$. By uniqueness of weak solutions for the
QS differential equation~\eqref{QSDE}, it follows that $k = k^\phi$,
in particular the elementary QS cocycle $k^\phi$ is completely
positive and contractive. Since $(\uAl)'' = \Al''$, the possibility of
choosing $R$ to be of the desired form follows from
Proposition~\ref{unital qcont cocs}.
\end{proof}

\begin{rems}
Theorem~\ref{thm: nonu} completes the quantum stochastic extension of
the Christensen--Evans characterisation of the generators of
norm-continuous, completely positive semigroups on $C^*$-algebras
(\cite{gran},~\cite{father}). This is so because norm-continuous
semigroups are automatically quasicontractive.

It is currently not known whether or not completely positive
elementary QS cocycles are automatically quasicontractive.
\end{rems}

\subsection{${}^*$-homomorphic elementary QS cocycles}

We now incorporate the known characterisation of stochastic generators
of \emph{weakly multiplicative} elementary QS cocycles into the
conclusions of Theorem~\ref{thm: nonu}. Thus let $\QSC_\wmult (\Al,
\noise)$ denote the set of $k$ in $\QSC(\Al, \noise)$ that satisfy
\[
\Dom k_t(a)^* \supset \init \aot \Exps
\quad \text{and} \quad
\bigl( k_t(a)^*|_{\init \aot \Exps} \bigr)^* k_t(b) = k_t(ab),
\qquad
t \in \Rplus, \ a,b \in \Al
\]
and, recalling the notations~\eqref{eqn: QSCh} and~\eqref{eqn: QSCu},
set
\begin{align*}
&
\QSC_\starhom(\Al, \noise) :=
\QSC_\wmult(\Al, \noise) \cap \QSC_\h(\Al, \noise)
\quad \text{and, when $\Al$ is unital},
\\
&
\QSF(\Al, \noise) :=
\QSC_\starhom(\Al, \noise) \cap \QSC_\unital(\Al, \noise),
\end{align*}
referred to as the \emph{${}^*$-homomorphic quantum stochastic
cocycles}, respectively \emph{quantum stochastic flows} (on $\Al$
with respect to $\noise$).

The following known automatic continuity result is relevant.

\begin{lemma}
\label{lemma: 6.5}
Let $\al \in L\big( \Al; \Op( \domain; \Hil) \big)$ where $\domain$ is
a dense subspace of a Hilbert space $\Hil$ and $\Op( \domain; \Hil)$
denotes the space of operators on $\Hil$ with domain
$\domain$. Suppose that $\al$ is hermitian and weakly multiplicative,
equivalently, $\al$ satisfies
\[
\al(a)^* \supset \al(a^*)
\quad \text{and} \quad
\al(a^*)^* \al(b) = \al(ab),
\qquad
a,b \in \Al.
\]
Then $\al$ is bounded operator valued. Moreover the map $\Al \to
B(\Hil)$, $a \mapsto \ol{\al(a)} = \al(a)^{**}$, is a ${}^*$-algebra
morphism.
\end{lemma}

\begin{proof}
In case $\Al$ is not unital, the prescription $a +\la 1 \mapsto \al(a)
+\la I_\Hil$ defines a linear extension of $\al$ to the unitisation of
$\Al$ and it is readily verified that the extension is also hermitian
and weakly multiplicative. We may therefore suppose without loss of
generality that $\Al$ is unital.

For $a \in \Al$ and $\zeta \in \domain$ with $\norm{a} \les 1$,
$\norm{\al(1) \zeta}^2 = \ip{\zeta}{\al(1) \zeta} =
\ip{\zeta}{\al(a^*a) \zeta} +\ip{\zeta}{\al(1 -a^*a) \zeta} =
\norm{\al(a) \zeta}^2 +\norm{\al( (1 -a^*a)^{1/2}) \zeta}^2$, so
$\norm{\al(a) \zeta}^2 \les \ip{\zeta}{\al(1) \zeta}$.  Taking $a = 1$
this shows that the operator $\al(1)$ is bounded with norm at most
one; feeding this back then gives $\norm{\al(a) \zeta}^2 \les
\norm{\zeta}^2$ for such $a$ and $\zeta$, so $\al(a)$ is likewise
bounded.  It follows that $\al$ is bounded operator valued. The second
part is straightforward to verify.
\end{proof}

It follows from Lemma~\ref{lemma: 6.5} that
\[
\QSC_\starhom(\Al, \noise)
=
\bigl\{ k \in \QSC_\bd(\Al, \noise): \ol{k_t( \cdot )} \text{ is a
  ${}^*$-algebra morphism for all } t \in \Rplus \bigr\},
\]
in particular,
\begin{equation}
\label{eqn: incl}
\QSC_\starhom(\Al, \noise) \subset \QSC_\cpc(\Al, \noise).
\end{equation}

Now let $\genwmult(\Al,\noise)$ denote the collection of maps
$\phi \in CB( \Al; \Al \otm B(\khat))$ which satisfy the following
higher-order structure relations:
\begin{equation}
\label{eqn: h.o.struct}
\phi_n(ab) =
\sum \phi_{\#\la}(a)(\la;n) \Delta[\la\cap\mu;n]
\phi_{\#\la}(b)(\mu;n),
\qquad
n \in \Nat, a,b \in \Al,
\end{equation}
in which the sum is over all pairs of subsets $(\la,\mu)$ of
$\{1, \cdots , n \}$ whose union is $\{1, \cdots , n \}$. The notation
here is as follows: for $k \in \Nat$, $\phi_k := \phi^{(k)} \circ
\cdots \circ \phi^{(1)}$ where $\phi^{(i)} := \phi^\Hil$ for $\Hil =
\khat^{\ot (i-1)}$ (see~\eqref{flipped lift}) and, for a subset $\la$ of
$\{1, \cdots , n\}$ with cardinality $j$, $\phi_j(a)(\la;n)$ denotes
the ampliation of $\phi_j(a)$ to $B(\init \ot \khat^{\ot n})$ acting
as the identity on each copy of $\khat$ labelled by an index from the
set $\{1, \cdots, n\} \setminus \la$ (see~\cite{LWhom}, for a fuller
explanation). Set $\genstarhom(\Al, \noise) := \genwmult(\Al, \noise)
\cap \genh(\Al, \noise)$ and, when $\Al$ is unital, $\genflow(\Al,
\noise) := \genstarhom(\Al, \noise) \cap \genu(\Al, \noise)$
(recall~\eqref{eqn: hVk} and~\eqref{eqn: uVk}).

\begin{thm}
\label{thm: 6.6}
The QS generation map $\Phi_{\Al, \noise}$ restricts to bijections
\begin{align*}
&
\genstarhom(\Al, \noise) \To
\El\QSC_\starhom(\Al,\noise)
\quad \text{and, when $\Al$ is unital},
\\
&
\genflow(\Al, \noise) \to
\El\QSF(\Al, \noise).
\end{align*}
\end{thm}

\begin{proof}
Write $\Phi$ for $\Phi_{\Al, \noise}$.

(a) Let $\phi \in CB(\Al; \Al \otm B(\khat))$. By Theorem~3.4
of~\cite{LWhom}, $k^\phi$ is weakly multiplicative if and only if
$\phi \in \genwmult(\Al, \noise)$. (From the fully coordinate-free
perspective adopted here, the proof given there is valid without the
stated separability assumption on the noise dimension space.)

(b) By the remark following Theorem~\ref{CB generation},
$\Phi(\genh(\Al, \noise)) \subset \QSC_\h(\Al, \noise)$, and
$\Phi(\genu(\Al, \noise)) \subset \QSC_\unital(\Al, \noise)$ when
$\Al$ is unital.

(c) If $k \in \El\QSC_\starhom(\Al, \noise)$ then, by
Theorem~\ref{thm: nonu} and the inclusion~\eqref{eqn: incl},
$k \in \Ran \Phi$ and so, by~(a) and~(b), $k = k^\phi$ where
$\phi \in \genwmult(\Al, \noise) \cap \genh(\Al, \noise) =
\genstarhom(\Al, \noise)$.

(d) If $\Al$ is unital and $k \in \El\QSF(\Al, \noise)$ then,
by~(c) and~(b), $k = k^\phi$ where $\phi \in \genstarhom(\Al, \noise)
\cap \genu(\Al, \noise) = \genflow(\Al, \noise)$.

From~(a) and~(b) it follows that $\Phi (\genstarhom(\Al, \noise))
\subset \El\QSC_\starhom(\Al, \noise)$ and, when $\Al$ is unital,
$\Phi (\genflow(\Al, \noise)) \subset \El\QSF(\Al, \noise)$; from~(c)
and~(d) it follows that these inclusions are equalities. The result
therefore follows from the injectivity of $\Phi$.
\end{proof}

We end this section by relating the above to the next section. Set
\begin{multline*}
\genstarhomone(\Al, \noise) := \\
\bigl\{
\phi \in L(\Al; \Al \otm B(\khat)):
\phi(a^*a) =
\phi(a)^* \iota(a) + \iota(a)^* \phi(a) + \phi(a)^* \Delta \phi(a)
\text{ for all } a \in \Al
\bigr\},
\end{multline*}
in which $\iota := \iota^\Al_{\khat}$, the ampliation introduced
in~\eqref{eqn: iota}, and set
\[
\genstruct( \Al, \noise) :=
\bigl\{
\phi \in \genstarhomone(\Al, \noise): \
\phi(1) = 0 \text{ if $\Al$ is unital}
\bigr\}.
\]

\begin{rems}
The following are equivalent conditions on a map $\phi \in L(\Al; \Al
\otm B(\khat))$ for it to be in $\genstarhomone(\Al, \noise)$:
\begin{rlist}
\item
$\phi^\dagger = \phi$ and
\begin{equation}
\label{eqn: carre}
\phi(ab) =
\phi(a) \iota(b) +\iota(a) \phi(b) +\phi(a) \Delta \phi(b),
\qquad a,b \in \Al.
\end{equation}
\item
$\phi$ has block matrix form $\left[\begin{smallmatrix} \mathcal{L} &
\delta^\dagger \\ \delta & \pi - \iota \end{smallmatrix}\right]$ where
$\iota = \iota^\Al_\noise$, $\pi$ is a ${}^*$-algebra morphism from
$\Al$ to $\Al \otm B(\noise)$, $\delta$ is a $\pi$-derivation from
$\Al$ to $\Al \otm \ket{\noise}$ (in other words $\delta(ab) =
\delta(a) b + \pi(a) \delta(b)$ for all $a,b \in \Al$) and
$\mathcal{L} \in L(\Al)$ is hermitian and satisfies
\begin{equation}
\label{eqn: Lin delta}
\mathcal{L}(ab) -\mathcal{L}(a)b -a\mathcal{L}(b) = \delta^\dagger(a)
\delta(b), \qquad a,b \in \Al.
\end{equation}
\end{rlist}
The equivalence with (i) follows from polarisation and the fact that
$\Al_\h = \Real$-$\Lin \Al_+$.

In particular, we see that if $\Al$ is unital then
\[
\phi(1) = 0 \text{ if and only if } \pi \text{ is unital}.
\]

The condition~\eqref{eqn: carre} is the case $n=1$ of~\eqref{eqn:
h.o.struct}, so $\genstarhom (\Al, \noise) \subset \genstarhomone
(\Al, \noise)$. Conversely, if $\phi \in \genstarhomone (\Al, \noise)$
then $\phi$ is necessarily completely bounded (as shown in the next
theorem), and further $\phi \in \genstarhom (\Al, \noise)$ when either
of the following two conditions obtain: ($\alpha$)\ $\phi(\Al) E_\zeta
\subset \Al \ot \ket{\khat}$ for all $\zeta \in \khat$, or ($\beta$)\
$\Al$ is a von Neumann algebra and $\phi$ is ultraweakly continuous
(\cite{LWhom},~\cite{LWhat}). In particular, $\genstarhom (\Al,\noise)
= \genstarhomone (\Al,\noise)$ if either $\noise$ or $\Al$ is finite
dimensional.
\end{rems}

Recall the definition of $\psi_R$ from~\eqref{psiR}.

\begin{thm}
\label{thm: starhomone}
Let $\phi \in L(\Al; \Al \otm B(\khat))$. Then the following are
equivalent\tu{:}
\begin{rlist}
\item
$\phi \in \genstarhomone (\Al, \noise)$.
\item
$\phi$ has block matrix form $\left[\begin{smallmatrix} \mathcal{L} &
\delta^\dagger \\ \delta & \pi - \iota \end{smallmatrix}\right]$
where, for some $l \in \Al''\, \vot \ket{\noise}$ and $h \in
\Al''_\h$, $\delta = \delta_{\pi, l}: a \mapsto la -\pi(a)l$ and
$\mathcal{L} = \mathcal{L}_{\pi, l, h}: a \mapsto l^*\pi(a)l
-\frac{1}{2} (l^*la +al^*l) + \I (ha -ah)$ in which $\pi$ is a
${}^*$-algebra morphism from $\Al$ to $\Al \otm B(\noise)$ and $\iota
= \iota^\Al_\noise$.
\item
$\phi = L^* \pi( \cdot ) L +\psi_R$ where $\pi$ is a ${}^*$-algebra
morphism from $\Al$ to $\Al \otm B(\noise)$, $L := \begin{bmatrix} l &
-I_{\init \ot \noise} \end{bmatrix} \in \Al'' \, \vot B(\khat;
\noise)$ and $R := \left[ \begin{smallmatrix} \I h - \tfrac{1}{2} l^*l
\\ l \end{smallmatrix}\right] \in \Al'' \, \vot \ket{\khat}$, for some
$l \in \Al'' \, \vot \ket{\noise}$ and $h \in \Al''_\h$.
\end{rlist}
Moreover, when these hold, if $\Al$ is unital then $\phi(1) = - T^*
\pi(1)^\perp T$ where $T:= \begin{bmatrix} l & I_{\init \ot
\noise} \end{bmatrix}$.
\end{thm}

\begin{proof}
(i)$\implies$(ii):
Let $\phi \in \genstarhomone (\Al, \noise)$, and let
$\left[\begin{smallmatrix} \mathcal{L} & \delta^\dagger \\ \delta &
\pi -\iota \end{smallmatrix}\right]$ be its block matrix form. Since
$\delta$ is a $\pi$-derivation and $\Ran \delta \subset \Al'' \, \vot
\ket{\noise}$, Theorem~2.1 of~\cite{ChrisEvans} implies that $\delta =
\delta_{\pi,l}$ for some $l \in \Al'' \, \vot \ket{\noise}$. Set
$\mathcal{L}_{\pi,l}: a \mapsto l^*\pi(a)l -\tfrac{1}{2}(l^*l a +a
l^*l)$. Then $\mathcal{L}_{\pi, l}$ is a hermitian map satisfying
$\Ran \mathcal{L}_{\pi, l} \subset \Al''$ and $\mathcal{L}_{\pi, l}
(ab) -\mathcal{L}_{\pi, l}(a) b -a \mathcal{L}_{\pi, l}(b) =
\delta^\dagger(a) \delta(b)$ for all $a,b \in \Al$. It follows that
$\mathcal{L} -\mathcal{L}_{\pi, l}$ is a hermitian derivation with
range in $\Al''$ and so, by another application of Theorem~2.1
of~\cite{ChrisEvans}, there is $k \in \Al''$ such that $\delta
=\delta_k: a \mapsto ka -ak$. Since $\delta_k$ is hermitian, $(ka^* -
a^*k)^* = ka-ak$ for all $a \in \Al$, or $k+k^* \in \Al' \cap
\Al''$. It follows that $\delta_k = \delta_{\I h}$ for a hermitian
element $h$ of $\Al''$.

(ii)$\implies$(i):
Let $l \in \Al'' \, \vot \ket{\noise}$ and $h \in (\Al'')_\h$ be such
that $\phi$ has block matrix form $\left[\begin{smallmatrix}
\mathcal{L}_{\pi, l, h} & \delta^\dagger_{\pi, l} \\ \delta_{\pi,l} &
\pi -\iota \end{smallmatrix}\right]$. Then $\delta_{\pi, l}$ is a
$\pi$-derivation and it is easily verified that $\mathcal{L}_{\pi, l,
h}$ is a hermitian map satisfying~\eqref{eqn: Lin delta}, so $\phi \in
\genstarhomone (\Al, \noise)$ by the first remark above.

The equivalence of~(ii) and~(iii) is readily verified, as is the given
consequence of $\Al$ being unital, and so the proof is complete.
\end{proof}

\begin{rems}
The above result was proved in~\cite{mother}, under the restriction
that $\Al$ is unital and $\noise$ is separable.
\end{rems}

From~(iii) we deduce the second of the following inclusions, which
should be seen in the light of the inclusion~\eqref{eqn: incl} and
Theorems~\ref{thm: nonu} and~\ref{thm: 6.6}.

\begin{cor}
\label{cor: struct in cpc}
$\genstruct (\Al, \noise) \subset \genstarhomone (\Al, \noise) \subset
\gencpc (\Al, \noise)$.
\end{cor}

\section{From structure map to completely positive QS cocycle}
\label{sec: structure}

For this section fix a unital $C^*$-algebra $\Al$ acting
nondegenerately on $\init$. Our goal is to identify sufficient
conditions on an operator $\phi$ from $\Al$ to $\Al \otm B(\khat)$
with dense domain $\Al_0$ for Theorem~\ref{thm: ABC} to apply and
yield a completely positive and unital QS cocycle $k$ on $\Al$
governed by the QS differential equation
\begin{equation}
\label{QSDE again}
\rd k_t = k_t \cdot \rd \La_\phi(t), \quad k_0 = \iota^\Al_\Fock.
\end{equation}
Our sufficient conditions on $\phi$ include the necessary algebraic
conditions for $k$ to be a \emph{QS flow on} $\Al$, that is a unital
and ${}^*$-homomorphic QS cocycle on $\Al$ (cf.\ the results of the
previous section in which $\phi$ was defined on \emph{all} of $\Al$).

\medskip
A \emph{structure map for} $(\Al, \noise)$ \emph{with domain} $\Al_0$
is an operator $\phi$ from $\Al$ to $\Al \otm B(\khat)$ with domain
$\Al_0$ such that
\begin{itemize}
\item[(S1)]
$\Al_0$ is a dense ${}^*$-subalgebra of $\Al$,
\item[(S2)]
$\phi$ satisfies $\phi = \phi^\dagger$ and $\phi(ab) = \phi(a) \iota
(b) + \iota(a) \phi(b) + \phi(a) \De \phi(b)$, for all $a,b \in
\Al_0$, where $\iota = \iota^\Al_\khat$ in the notation~\eqref{eqn:
iota},
\item[(S3)]
$\Al_0$ contains $1_\Al$ and satisfies $\phi (1_\Al) = 0$.
\end{itemize}

\begin{rems}
(i) The notation here is that $\phi^\dagger$ is the map on $\Al$ given
by $a \mapsto \phi(a^*)^*$, thus~(S2) includes the condition that
$\phi$ is hermitian.

(ii) If a QS flow on $\Al$ strongly satisfies~\eqref{QSDE again} on a
dense ${}^*$-subalgebra $\Al_0$ of $\Al$ containing $1_\Al$ then, it
follows from the quantum It\^{o} product formula that $\phi$ is
necessarily a structure map.

(iii) If a QS flow $j$ on $\Al$ is elementary then $j = k^\phi$ for a
completely bounded structure map $\phi$ (with domain $\Al$) by
Theorem~5.10 of \cite{father}. Conversely, necessary and sufficient
conditions for a completely bounded map $\phi$ to stochastically
generate a QS flow are given in~\cite{LWhom}. In favourable cases
(notably, if $\phi(\Al) E_\xi \subset \Al \ot \ket{\khat}$ for each
$\xi \in \khat$, or if $\Al$ is a von Neumann algebra and $\phi$ is
ultraweakly continuous) these conditions reduce to $\phi$ simply being
a structure map.
\end{rems}

\begin{lemma}
\label{lemma: X}
Let $\mu = (\id_\Al \otm \, \om) \circ \phi$ where $\phi$ is a
structure map for $(\Al, \noise)$ with domain $\Al_0$ and $\om \in
B(\khat)_{*, +}$. Then $\mu$ is hermitian, vanishes at $1_\Al$, and
satisfies
\[
\mu (a^* a) - \mu (a)^* a - a^* \mu (a) \ges 0, \qquad a \in \Al_0.
\]
\end{lemma}

\begin{proof}
Being a composition of hermitian maps, $\mu$ is hermitian. Since
$\phi$ vanishes at $1_\Al$, $\mu$ does too. Let $a \in \Al_0$, then
\begin{align*}
\mu(a^*a) &= (\id_\Al \otm\, \om) \big( \phi(a)^* \iota(a)
+\iota(a)^* \phi(a) +\phi(a)^* \De \phi(a) \big) \\
&= \mu(a)^* a +a^* \mu(a) + (\id_\Al \otm \, \om) \big( \phi(a)^*
\De \phi(a) \big).
\end{align*}
The result follows since $\id_\Al \otm \, \om$ is a positive map.
\end{proof}

\begin{rem}
It follows that the map $\mu$ is conditionally positive: if $a \in
\Al_0$ and $b \in \Al$ are such that $ab=0$ then $b^* \mu(a^*a) b \ges
0$. This is a necessary condition for $\mu$ to be the generator of a
positive semigroup on $\Al$. It is also a sufficient condition if
$\mu$ is bounded with domain $\Al$.
\end{rem}

If $k$ is a QS flow then for any map $\Ga: \Ix \To \noise$ it follows
from Proposition~\ref{kGaPGa} that $k^\Ga$ is a cocycle on (the
operator system) $\Al \otm B(\hil)$, where $\hil = l^2(\Ix)$, and that
it is multiplicative when restricted to the ${}^*$-subalgebra $\Al \ot
B(\hil)$.  Thus if $k = k^\phi$ for some $\phi \in CB (\Al; \Al \otm
B(\khat))$ and $\Ga$ is bounded, then the map $\Phi_\Ga$ from
Theorem~\ref{Z} ought to be a structure map.

We want to establish this directly when starting with an unbounded
structure map, and to that end we must extend the use of the notation
introduced in~\eqref{flipped lift}. Let $\nu$ be an operator from
$\Al$ to $\Al \otm B(\kil_1; \kil_2)$ with domain $\Al_0$, and let
$\hil$ be a Hilbert space. Then denote by $\nu^\hil$ the operator from
$\Al \ot B(\hil)$ to $(\Al \ot B(\hil)) \otm B(\noise_1; \noise_2)$
with domain $\Al_0 \aot B(\hil)$ determined by
\[
\nu^\hil ( a \ot T ) := \Pi (\ups(a) \ot T)
\]
where $\Pi$ is the tensor flip from $B(\init \ot \kil_1 \ot \hil;
\init \ot \kil_2 \ot \hil)$ to $B(\init \ot \hil \ot \hil_1; \init \ot
\hil \ot \kil_2)$. When $\phi$ is completely bounded, the resulting
operator is simply a restriction of that defined in~\eqref{flipped
lift}.

\begin{propn}
\label{stmapprops}
Let $\phi$ be a structure map for $(\Al, \noise)$ with domain
$\Al_0$. Then the following hold.
\begin{alist}
\item
Given a Hilbert space $\Hil$, $\phi^\Hil$ is a structure map for
$(\Al \ot B(\Hil), \noise)$ with domain $\Al_0 \aot B(\Hil)$.

\item
Let $\Ga: \Ix \To \noise$ be a bounded map and set $\hil :=
l^2(\Ix)$. Define $\Phi_\Ga$ to be the operator from $\Al \otm
B(\hil)$ to $(\Al \otm B(\hil)) \otm B(\khat)$ with domain $\Al_0 \aot
B(\hil)$ given by
\[
\Phi_\Ga = (\id_\Al \otm \ups_\Ga) \circ \phi^\hil +\id_\Al \otm
\Theta_\Ga
\]
for $\ups_\Ga$ and $\Theta_\Ga$ as in \tu{Theorem~\ref{thm:
KGagen}}. Then $\Phi_\Ga$ is a structure map.
\end{alist}
\end{propn}

\begin{proof}
(a) That $\phi^\Hil$ is hermitian and satisfies $\phi^\Hil (1_{\Al \ot
B(\Hil)}) =0$ is clear, moreover the identity
\[
\phi^\Hil (AB) = \phi^\Hil (A) \iota(B) + \iota(A) \phi^\Hil (B)
+\phi^\Hil (A) \De \phi^\Hil (B),
\]
is readily verified for simple tensors $A$ and $B$ in $\Al_0 \aot
B(\Hil)$. Part~(a) therefore follows by bilinearity.

(b) For ease of reading we suppress all subscripts and ampliations for
the rest of the proof. By definition, $\Phi = \Om +\Theta$ where, for
$A \in \Al_0 \aot B(\hil)$,
\begin{gather*}
\Om (A) = \wt{F}^* \phi^\hil(A) \wt{F} \quad \text{and} \quad \Theta
(A) = F^* \iota(A) +\iota(A) F +F^* (A \ot \De) F, \quad \text{for} \\
F := \begin{bmatrix} -\frac{1}{2} L^*L & -L^* \\ L & 0 \end{bmatrix}
\quad \text{and} \quad
\wt{F} := I +\De F = \begin{bmatrix} I & 0 \\ L &
I \end{bmatrix},
\end{gather*}
in which $L$ is the $\init$-ampliation of the operator in
$B(\hil; \hil \ot \noise )$ defined as in Lemma~\ref{Wgen}, so that
\[
\Theta (A) = \begin{bmatrix} L^* (A \ot I) L -\frac{1}{2} L^*L A
-\frac{1}{2} A L^*L & L^* (A \ot I) -AL^* \\ (A \ot I) L -LA &
0 \end{bmatrix}.
\]
Note that $\wt{F} \De = \De = \De \wt{F}^*$ and so
\[
\iota(A) \wt{F} - \wt{F} \iota(A) = \begin{bmatrix} 0 & 0 \\ (A \ot
I)L -LA & 0 \end{bmatrix} = \De \Theta(A) = \wt{F} \De \Theta(A).
\]
Since $\phi^\hil$ is a structure map this implies that, for $A \in
\Al_0 \aot B(\hil)$,
\begin{multline*}
\Om (A^*A) -\Om (A^*) \iota (A) -\iota (A)^* \Om (A) \\
\begin{aligned}
&= \wt{F}^* \phi^\hil(A)^* \big( \iota(A) \wt{F} -\wt{F} \iota(A)
\big) +\big( \wt{F}^* \iota (A)^* -\iota(A)^* \wt{F}^* \big) \phi^\hil
(A) \wt{F} +\wt{F}^* \phi^\hil (A)^* \De \phi^\hil (A) \wt{F} \\
&= \wt{F}^* \phi^\hil (A)^* \wt{F} \De \Theta (A) +\Theta (A)^* \De
\wt{F}^* \phi^\hil (A) \wt{F} +\wt{F}^* \phi^\hil (A)^* \wt{F} \De
\wt{F}^* \phi^\hil (A) \wt{F} \\
&= \Om (A)^* \De \Theta (A) +\Theta (A)^* \De \Om (A) +\Om (A)^* \De
\Om (A).
\end{aligned}
\end{multline*}
Furthermore $\Theta$ is a structure map since $F +F^* +F^* \De F = 0 =
F +F^* +F \De F^*$ (\emph{cf}.~Lemma~\ref{Wgen}). Thus, since $\Phi
=\Om +\Theta$, this implies that
\[
\Phi (A^*A) -\Phi (A^*) \iota(A) -\iota(A)^* \Phi (A) = \Phi (A)^* \De
\Phi(A)
\]
for $A \in \Al_0 \aot B(\hil)$, as required. That $\Phi$ is hermitian
and $\Phi (1_{\Al \ot B(\hil)}) =0$ are easily verified.
\end{proof}
	
Now fix an orthonormal basis $\eta = (d_i)_{i \in \Ix_0}$ for $\noise$
and recall our convention~\eqref{eta bar defn}. Continuing with the
notation introduced before Theorem~\ref{thm: ABC}, we write
$\Phi_{[\Jx]}$ and $\varphi_{[\Jx]}$ as special cases of $\Phi_\Ga$
and $\varphi_\Ga$. Also recall the $\phi_{x,y}$ notation introduced
in~\eqref{xygen}.

\begin{propn}
\label{propn: XYZ}
Let $\phi$ be a structure map for $(\Al, \noise)$ with domain $\Al_0$.
Suppose that the following hold\tu{:}
\begin{itemize}
\item[(1)]
For each $n \in \Nat$, $\Mat_n(\Al_0)$ is \emph{square-root closed},
that is
\[
\Mat_n(\Al_0)_+ = \big\{ A^2: A \in \Mat_n(\Al_0)_+ \big\},
\ \text{ where }
\Mat_n(\Al_0)_+:= \Mat_n(\Al)_+ \cap \Mat_n(\Al_0).
\]
\item[(2)]
For all $\al, \be \in \Ix$, $\phi_{d_\al, d_\be}$ is a pregenerator of
a $C_0$-semigroup $\msgp^{(\al, \be)}$ on $\Al$.
\end{itemize}
Then for all $\Jx_0 \fsubset \Ix_0$ the Schur-action semigroup
$\msgp^{[\Jx]}$ on $\Mat_\Jx (\Al)$ is completely positive and
unital.
\end{propn}

\begin{proof}
Let $\Jx_0 \fsubset \Ix_0$. Define a Schur-action operator on
$\Mat_\Jx (\Al)$ with domain $\Mat_\Jx (\Al_0)$ by
\[
\varphi_{[\Jx]} := \big[ \phi_{d_\al, d_\be} \big]_{\al,\be \in \Jx}
\schur
\]
Then $\varphi_{[\Jx]}$ is a pregenerator of the $C_0$-semigroup
$\msgp^{[\Jx]}$, with $\varphi_{[\Jx]} = (\id_{\Mat_\Jx (\Al)} \otm\,
\om_{\wh{0}, \wh{0}}) \circ \Phi_{[\Jx]}$. However, $\Phi_{[\Jx]}$ is
a structure map on $\Mat_\Jx (\Al)$ by
Proposition~\ref{stmapprops}\,(b), and so
\[
\varphi_{[\Jx]} (A^*A) -\varphi_{[\Jx]} (A^*) A -A^*
\varphi_{[\Jx]} (A) \ges 0, \qquad A \in \Mat_\Jx (\Al_0)
\]
by Lemma~\ref{lemma: X}. It follows from Proposition~3.2.22
of~\cite{BrR} that $\varphi_{[\Jx]}$ is dissipative. Thus its closure
$\ol{\varphi_{[\Jx]}}$ is dissipative by Proposition~3.1.15
of~\cite{BrR}. However, $\ol{\varphi_{[\Jx]}}$ is the generator of the
$C_0$-semigroup $\msgp^{[\Jx]}$, so $\id_{\Mat_n (\Al)} -\al \,
\ol{\varphi_{[\Jx]}}$ is surjective for some $\al > 0$ by the
Hille--Yosida Theorem (\cite{BrR}, Theorem~3.1.10).  Since its
generator is dissipative, $\msgp^{[\Jx]}$ is contractive by the
Lumer--Phillips Theorem (\cite{BrR}, Theorem~3.1.16). Moreover
$\msgp^{[\Jx]}$ is unital on the $C^*$-algebra $\Mat_{[\Jx]} (\Al)$
since $\varphi_{[\Jx]}$ vanishes at the identity. Putting these two
properties together shows that $\msgp^{[\Jx]}$ is positive
(\cite{Paulsen}, Proposition~2.11).

To get complete positivity of $\msgp^{[\Jx]}$, note that for any $n
\in \Nat$ the operator $\phi^{\Comp^n}$ is a structure map for $(\Mat_n
(\Al), \noise)$ by Proposition~\ref{stmapprops}\,(a), and that
conditions~(1) and~(2) hold for this map. If $\msgq^{[\Jx]}$ denotes
the contraction semigroup on $\Mat_\Jx \big( \Mat_n (\Al) \big)$
obtained by running the argument above for this lifted structure map
then $\msgq^{[\Jx]} = \Pi \circ (\msgp^{[\Jx]} \otm \id_{\Mat_n
(\Comp)})$ where $\Pi: \Mat_n \big( \Mat_\Jx (\Al) \big) \To \Mat_\Jx
\big( \Mat_n (\Al) \big)$ is the natural flip isomorphism. The result
follows.
\end{proof}

Combined with Theorem~\ref{thm: ABC}, Proposition~\ref{propn: XYZ}
yields the following stochastic generation theorem for completely
positive and unital QS cocycles.

\begin{thm}
\label{thm: stoch gen}
Let $\phi$ be a structure map for $(\Al, \noise)$ with domain
$\Al_0$. Suppose that the following hold\tu{:}
\begin{itemize}
\item[(1)]
For all $n \in \Nat$, $\Mat_n(\Al_0)$ is square-root closed.
\item[(2)]
There is an orthonormal basis $(d_i)_{i \in \Ix_0}$ for $\noise$ such
that, for all $\al, \be \in \Ix:= \{0\} \cup \Ix_0$,
$\phi_{d_\al, d_\be}$
is a pregenerator of a $C_0$-semigroup $\msgp^{(\al, \be)}$ on
$\Al$.
\end{itemize}
Then there is a unique process $k \in \QSCb (\Al, \noise)$ which is
locally norm bounded and such that, for all $\al, \be \in \Ix$,
$\msgp^{(\al, \be)}$ is its $(d_\al, d_\be)$-associated
semigroup. Moreover, $k$ is completely positive and unital, and it is
the unique weakly regular $\Exps_{\Til(\eta)}$-weak solution of the QS
differential equation~\eqref{QSDE again} on $\Al_0$ for the domain
$\init \aot \Exps_{\Til(\eta)}$\tu{;} if both $\init$ and $\noise$ are
separable then it is a strong solution.
\end{thm}

\begin{rems}
Assumption~(1) clearly applies when $\Al_0$ is the dense
${}^*$-subalgebra associated with an AF algebra. In the next section
we give such an example in which assumption~(2) also holds.

Interesting recent work on the construction of QS flows is nicely
complementary to ours, in that it takes as \emph{standing hypothesis}
the existence of a completely positive and contractive QS cocycle
weakly satisfying~\eqref{QSDE again} on a dense ${}^*$-subalgebra
containing the identity of the algebra, for a structure map $\phi$
(see~\cite{DGS}, Theorem~3.1 and Definitions~2.5 and~2.7).
\emph{Warning}: their use of the terminology QS flow is different to
ours.
\end{rems}

\section{The quantum exclusion process}
\label{sec: asym. qep}

Symmetric quantum exclusion processes have recently been considered by
a number of authors. The original paper to tackle the challenge of
extending the classical theory of exclusion processes to the quantum
domain was~\cite{Reb}. This, like subsequent work, has focused on the
underlying Markov semigroup. Subsequently conditions were found for
the construction of symmetric quantum exclusion \emph{processes}
(\cite{BeW}), using the theory of multiple quantum Wiener integrals
developed in~\cite{LWhom}. In this final section we demonstrate how
Theorem~\ref{thm: stoch gen} may be employed in the construction of
quantum exclusion processes. Our approach is complementary to that
of~\cite{BeW} and~\cite{BWFK}.

Fixing a nonempty set $\Rx$, let $\Al = CAR(\Rx)$ denote the CAR
algebra over $\Rx$ in its Fock representation (\cite{BrR2}). A useful
description arises by putting a total order on $\Rx$. The
antisymmetric Fock space over $l^2(\Rx)$ may be naturally identified
with the Hilbert space $\init := l^2(\Ga_\Rx)$ where $\Ga_\Rx := \{\si
\subset \Rx: \# \sigma < \infty\}$, with the Fock space Fermi
annihilation and creation operators given by
\[
(b_r F) (\si)= \ind_{r \notin \si} \ep (\si, r) F (\si \cup r)
\ \text{ and } \
(b_r^* F) (\si)= \ind_{r \in \si} \ep (\si, r) F(\si \setminus r),
\]
in which the singleton set $\{r\}$ is abbreviated to $r$. The
notation here is as follows:
\[
\ep (\si, r) := (-1)^{ n(\sigma, r) }
\quad \text{where} \quad
n (\sigma, r) := \# \{s \in \si: \, s > r\}
\quad \text{for} \quad
\si \in \Ga_\Rx \ \text{and} \ r \in \Rx.
\]
The anticommutation relations
\[
b_r b_s + b_s b_r = 0,
\quad
b_r b_s^* + b_s^* b_r = \begin{cases} 1_\Al & \text{if } r = s, \\
0 & \text{if } r \neq s, \end{cases}
\]
imply that $b_r b_r^* b_r = (1_\Al - b_r^* b_r) b_r = b_r$ ($r \in
\Rx$), so each $b_r$ is a (nonzero) partial isometry. Let $\Al_0$
denote ${}^*$-$\Alg \{b_r: r \in \Rx\}$. Thus $\Al_0$ is a dense
${}^*$-subalgebra of $\Al$ containing $1_\Al$; it is weak operator
dense in $B(\init)$. In terms of the standard basis $(e_\si)_{\si \in
\Ga_\Rx}$ for $\init$,
\[
b_r e_\si = \ind_{r \in \si} \, \ep (\si, r) \, e_{\si \setminus r}
\quad \text{and} \quad
b_r^* e_\si = \ind_{r \notin \si} \, \ep (\si, r) \, e_{\si \cup r},
\quad
\si \in \Ga_\Rx, r \in \Rx.
\]

For $\Sx \subset \Rx$, set $\Al_\Sx := {}^*\text{-}\Alg \{b_s: s \in
\Sx\}$. Then, by the anticommutation relations, each $\Al_\Sx$ is a
${}^*$-subalgebra of $\Al$ containing $1_\Al$ with linear basis
$\{b_\si^* b_\tau: \si, \tau \fsubset \Sx\}$, where
\[
b_\tau := \prod_{t \in \tau}^{\longrightarrow} b_t
\quad \text{and} \quad
b_\si^* := (b_\si)^* = \prod_{s \in \si}^{\longleftarrow} b_s^*,
\]
with the convention $b_\emptyset = 1_\Al$. Thus
\[
b_\si^* \, e_\emptyset = e_\si, \text{ for } \si \in \Ga_\Rx,
\quad\
b_\tau e_\emptyset = 0 \text{ for } \tau \in \Ga_\Rx \setminus
\{\emptyset\},
\]
in particular $e_\emptyset$ is a cyclic vector for the $C^*$-algebra
$\Al$ (it is the Fermi Fock vacuum vector), and
\[
\Al_0 = \bigcup_{\Sx \fsubset \Rx} \Al_\Sx.
\]
Moreover, for each $\Sx \fsubset \Rx$, $\Al_\Sx$ is finite dimensional
and thus a closed subspace of $\Al$; it is a $C^*$-algebra isomorphic
to $B \bigl( l^2(\Ga_\Sx) \bigr)$. Thus $\Al_0$ is square-root closed,
as is $\Mat_n (\Al_0)$ for each $n \in \Nat$. In addition, $\Al$ is
separable and an AF algebra if and only if $\Rx$ is countable.

The elements of $\Rx$ are used to label sites at which Fermionic
particles may exist, with the operator $b_r$ representing the
annihilation of a particle at site $r$, and $b_r^*$ its creation.

Let $\{\al_{r,s}: r,s \in \Rx\}$ be a fixed collection of (complex)
\emph{amplitudes}, and set
\[
\inter (r) := \{s \in \Rx: \al_{r,s} \neq 0 \text{ or } \al_{s,r} \neq
0\}, \text{ and } \inter^+(r) := \inter(r) \cup \{r\}.
\]
Thus $\inter(r)$ is the collection of sites that interact with site
$r$; $\#\inter(r)$ is termed the \emph{valency of the site} $r$. We
make the finite valency assumption
\begin{equation}
\label{finite interactions}
\#\inter(r) < \infty \text{ for all } r \in \Rx.
\end{equation}
The transport of a particle from site $r$ to site $s$ with amplitude
$\al_{r,s}$ is described by the operator $t_{r,s} := \al_{r,s} b^*_s
b_r$. Also let $\{h_r: r \in \Rx\}$ be a fixed set of (real)
\emph{site energies}. Define bounded operators $\rho_{r,s}$,
$\tau_{r,s}$ and $\delta_r$ on $\Al$ by
\begin{align*}
&
\rho_{r,s} (a) := [t_{r,s},a] = \al_{r,s} (b^*_s b_r a -a b^*_s b_r),
\\
&
\tau_{r,s} (a) := -\frac{1}{2} \big( t^*_{r,s} \rho_{r,s} (a)
+\rho^\dagger_{r,s} (a) t_{r,s} \big) = -\frac{1}{2} \big( t^*_{r,s}
[t_{r,s},a] +[a,t^*_{r,s}] t_{r,s} \big) \ \text{ and}
\\
&
\delta_r(a) := \I h_r [b_r^* b_r, a].
\end{align*}
Thus, for each $r, s \in \Rx$, $\rho_{r,s}$ and $\delta_r$ are
derivations. In particular $\delta_r(1_\Al) = \rho_{r,s}(1_\Al) = 0$
so also $\tau_{r,s} (1_\Al) = 0$, and the following identity holds:
$\tau_{r,s}(ab) -\tau_{r,s}(a) b -a \tau_{r,s}(b) =
\rho^\dagger_{r,s}(a) \rho_{r,s}(b)$ for all $a,b \in \Al$.  Moreover,
each $\tau_{r,s}$ is hermitian, as is each $\delta_r$.  Noting that,
for $\Sx \fsubset \Rx$,
\begin{equation}
\label{comm relation}
[b^*_s b_r, a] = 0 \ \text{ for all } a \in \Al_\Sx, r, s \notin \Sx,
\end{equation}
we see that, for all $\Sx \fsubset \Rx$ and $a \in \Al_\Sx$,
$\delta_r(a) = 0$ unless $r \in \Sx$, in which case $\delta_r(a) \in
\Al_\Sx$, whereas $\rho_{r,s}(a) = 0$ unless either $r \in \Sx$ and $s
\in \inter (r)$, in which case $\rho_{r,s}(a) \in \Al_{\Sx \cup \inter
(r)}$, or $s \in \Sx$ and $r \in \inter (s)$, in which case
$\rho_{r,s}(a) \in \Al_{\Sx \cup \inter (s)}$, and similarly for
$\tau_{r,s}(a)$.  It follows that, under the finite valency
assumption~\eqref{finite interactions}, for all $\Sx \fsubset \Rx$ and
$a \in \Al_\Sx$, the sets $\{ r \in \Rx: \delta_r(a) \neq 0 \}$, $\{q
\in \Rx \times \Rx: \rho_q(a) \neq 0\}$ and $\{q \in \Rx \times \Rx:
\tau_q(a) \neq 0\}$ are all finite. Therefore, setting $\noise:= l^2
(\Rx \times \Rx)$ and letting $\eta = (f_q)_{q \in \Rx \times \Rx}$ be
its standard orthonormal basis, there is a well-defined operator
$\phi := \left[\begin{smallmatrix} \Lindbladian & \rho^\dagger \\ \rho
& 0 \end{smallmatrix}\right]$ from $\Al$ to $\Al \ot B(\khat)$
with domain $\Al_0$ given by $\Lindbladian := \delta + \tau$,
\[
\delta(a) := \sum_{r \in \Rx} \delta_r(a),
\quad
\tau(a) := \!\!\sum_{q \in \Rx \times \Rx} \tau_q(a)
\quad \text{and} \quad
\rho(a) := \!\!\sum_{q \in \Rx \times \Rx} \rho_q(a) \ot
\ket{f_q}.
\]
For $\Sx \fsubset \Rx$, set $\Sx^+ := \bigcup_{r \in \Sx} \inter^+(r)$
and $\noise_\Sx := \Lin \{f_q: q \in \Sx \times \Sx\}$, then
\[
\de (\Al_\Sx) \subset \Al_\Sx,
\qquad
\tau (\Al_\Sx) \subset \Al_{\Sx^+}
\quad
\text{and}
\quad
\rho (\Al_\Sx) \subset \Al_{\Sx^+} \aot \ket{\noise_{\Sx^+}},
\]
so that $\phi$ enjoys the approximate invariance property $\phi
(\Al_\Sx) \subset \Al_{\Sx^+} \aot \Lin \big\{ \dyad{x}{y}: x,y \in
\wh{ \noise_{\Sx^+} } \big\}$. In particular $\Ran \phi \subset \Al_0
\aot B_{00} (\khat)$, where $B_{00}$ denotes bounded finite rank
operators. Since $\tau$ and $\delta$ are hermitian, $\delta$ is a
derivation, $\rho$ an $\iota^\Al_\noise$-derivation, and $\tau(ab)
-\tau(a)b -a\tau(b) = \rho^\dagger(a) \rho(b)$ for all $a,b \in
\Al_0$, it follows that $\phi$ is a structure map for $(\Al, \noise)$
with domain $\Al_0$.

\begin{thm}
\label{thm: 7.4}
Let $\Al = CAR(\Rx)$, in its Fock representation, for a nonempty set
$\Rx$ and let $\phi := \left[\begin{smallmatrix} \Lindbladian &
\rho^\dagger \\ \rho & 0 \end{smallmatrix}\right]$ be the structure
map defined as above, in terms of a set of complex amplitudes
$\{\al_{r,s}: r,s \in \Rx\}$ satisfying the finite valency
condition~\eqref{finite interactions} and a set of real site energies
$\{h_r: r \in \Rx\}$. Assume that $\Lindbladian$ is a pregenerator of
a $C_0$-semigroup on $\Al$. Then there is a unique process $k^\phi$ in
$\QSCcpu (\Al, \noise)$ such that $\phi_{x,y} := (\id_\Al \ot
\omega_{\wh{x}, \wh{y}}) \circ \phi -\chi(x,y) \id_\Al$ is a
pregenerator of its $(x,y)$-associated semigroup, for all $x, y \in
\{0\} \cup \{f_q: q \in \Rx \times \Rx\}$. Moreover, $k^\phi$ is also
the unique weakly regular $\Exps_{\Til(\eta)}$-weak solution of the QS
differential equation~\eqref{QSDE again} on $\Al_0$ for the domain
$\init \aot \Exps_{\Til(\eta)}$, satisfying the equation strongly if
the set $\Rx$ is countable.
\end{thm}

\begin{proof}
By identities~\eqref{cpttransfs},
\begin{equation*}
\phi_{0,f_q}
= \Lindbladian +\rho^\dagger_q -\tfrac{1}{2} \id_\Al,
\quad
\phi_{f_q,0}
= (\phi_{0,f_q})^\dagger
\quad  \text{and} \quad
\phi_{f_q,f_{q'}}
= \Lindbladian +\rho_q +\rho_{q'}^\dagger +(\de_{q,q'} -1) \id_\Al
\end{equation*}
for all $q,q' \in \Rx \times \Rx$. Therefore, for all $x, y \in \{0\}
\cup \{f_q: q \in \Rx \times \Rx\}$, $\phi_{x,y}$ is a bounded
perturbation of $\Lindbladian$, and so is a pregenerator of a
$C_0$-semigroup on $\Al$. As noted already, $\Mat_n (\Al_0)$ is
square-root closed, since any $A \in \Mat_n (\Al_0)_+$ belongs to the
$C^*$-algebra $\Mat_n (\Al_\Sx)$ for some $\Sx \fsubset
\Rx$. Therefore the result follows from Theorem~\ref{thm: stoch gen}.
\end{proof}

\begin{rems}
If an orthonormal basis $\eta'$ is chosen for $\noise$ other than the
standard one then the resulting QS cocycle $k'$ is unitarily conjugate
to $k^\phi$ via the second quantisation operator $\Gamma
(I_{L^2(\Rplus)} \ot V)$, $V$ being the unitary on $\noise$ which
exchanges the bases.

Strong unital ${}^*$-homomorphic solutions for the QS differential
equation~\eqref{QSDE again} (where, as is necessary, $\phi$ is assumed
to be a structure map) are sought in~\cite{BeW}. Existence is proved
assuming uniform boundedness of the valencies and energies, and
modulus-symmetry of the amplitudes: $|\al_{s,r}| = |\al_{r,s}|$ for
all $r,s \in \Rx$, along with certain coupled conditions. These latter
constraints, on magnitudes of amplitudes and sizes of valencies, are
forged from growth restrictions, on iterates of $\phi$ applied to the
generators $b_r$, required for the convergence of sums of relevant
multiple quantum Wiener integrals. Under strengthened conditions the
solution is shown to be a QS cocycle. In the preprint~\cite{BWFK}
these symmetry and valency restrictions are loosened, and existence
theorems for weak completely positive contractive solutions
of~\eqref{QSDE again} are established by means of Feynman--Kac type
perturbations of the ${}^*$-homomorphic processes constructed
in~\cite{BeW}. The cocycle property for these is then recovered under
conditions allowing the application of Theorem~\ref{thm:
ABC}/\ref{thm: stoch gen}, as in the current paper.

In Rebolledo's work, and subsequent study (\cite{questions}), the
coefficients $\al_{r,s}$ are assumed to be real; moreover a large part
of the analysis is carried out in the $W^*$-category, focusing on the
minimal quantum dynamical semigroup on $B(\init)$. In that context,
our results (like those of~\cite{BeW} and~\cite{BWFK}) deliver Feller
properties for the resulting semigroups --- that is, invariance of the
$C^*$-algebra $\Al = CCR(\Rx)$ and strong continuity there.
\end{rems}

We finish with a special case of the quantum exclusion process in
which the index set $\Rx$ is an integer lattice, and where the fact
that $\Lindbladian$ is a pregenerator is established using the
following result.

\begin{thm}[\cite{BrK}, Theorem~4.2.1]
\label{BrKgen}
Let $L$ be a dissipative operator on a Banach space $\X$. Suppose
that $(\X(n))_{n=1}^\infty$ is an increasing sequence of closed
subspaces of $\X$ contained in $\Dom L$ such that
\begin{alist}
\item
$\X_0 := \bigcup_{n=1}^\infty \X(n)$ is dense in $\X$, and
\item
there are constants $M \ges 0$ and $\al > 0$ and, for each $n \in
\Nat$ and $m \in \Int_+$, an operator $L^{n,m}$ from $\X(n)$ to $\X$
such that, for all such $n$ and $m$,
\[
\Ran L^{n,m} \subset \X(n+m)
\quad \text{and} \quad
\norm{L|_{\X(n)} - L^{n,m}} \les M n e^{-\al m}.
\]
\end{alist}
Then $L$ is a pregenerator of a \tu{(}contractive\tu{)}
$C_0$-semigroup on $\X$ and $\X_0$ is a core for $\ol{L}$.
\end{thm}

\begin{example}
Let $\Al = CAR(\Rx)$ acting on the Hilbert space $l^2(\Ga_\Rx)$, where
$\Rx = \Int^d$ which we view as a metric subspace of
$l^\infty(\{ 1, \cdots , d \})$, and let $\phi$ be the structure map
for $(\Al, \noise)$ with domain $\Al_0$ defined as above, in terms of
a given set of amplitudes and energies. We claim that $\Lindbladian$
is a pregenerator of a $C_0$-semigroup on $\Al$, so that
Theorem~\ref{thm: 7.4} applies, under the following assumptions.
\begin{enumerate}
\item[I]
There is $D \in \Nat$ such that $\al_{r,s} = 0$ whenever the pair
$(r,s)$ satisfies $\norm{r-s}_\infty > D$.
\item[II]
There is $K \in \Rplus$ such that, for all  $(r,s) \in \Rx \times \Rx
\setminus \{ (0,0) \}$,
\[
\abs{\al_{r,s}}^2 \les K \big( \max \{\norm{r}_\infty,
\norm{s}_\infty\} \big)^{2-d}.
\]
\end{enumerate}
Assumption~I strengthens that of~\eqref{finite interactions} to there
being a uniform limit on the range of interaction and hence a uniform
bound on the valencies. In dimensions $d=1$ and $2$, assumption~II
covers the physically reasonable situation of uniformly bounded
amplitudes. Unlike in~\cite{BeW}, no symmetry condition is imposed on
the amplitudes and no bounds on the energies are required.

To see the validity of our claim, first note that $\Lindbladian$ is
conditionally completely positive and so, since its domain is
square-root closed and contains $1_\Al$, Proposition~3.2.22 of
\cite{BrR} implies that $\Lindbladian$ is dissipative. Let $n \in
\Nat$ and set $\Al(n) := \Al_{\Rx_n}$ where $\Rx_n := \{r \in \Rx:
\norm{r}_\infty \les nD\}$; since $\# \Rx_n = (1+2nD)^d < \infty$,
$\Al(n)$ is finite dimensional and thus a closed subspace of
$\Al$. Also note that $\bigcup_{n=1}^\infty \Al(n) = \Al_0$ and, by
assumption~I, $\Rx^+_n \subset \Rx_{n+1}$ so $\Al_{\Rx_n^+} \subset
\Al(n+1)$ for all $n \in \Nat$. For $m \in \Int_+$ define the map
$\Lindbladian^{n,m}: \Al(n) \to \Al$ by
\[
\Lindbladian^{n,m}: \Al(n) \to \Al, \quad
a \mapsto \begin{cases} \delta(a) +\sum_{r,s \in \Rx_n}
\tau_{r,s}(a) & \text{if } m=0, \\ \Lindbladian(a) & \text{if } m \ges
1. \end{cases}
\]

Fix $n \in \Nat$ and $a \in \Al(n)$ such that $\norm{a} \les 1$. Since
$\Al(n+1) \subset \Al(n+m)$ for all $m \in \Nat$, the claim follows by
appeal to Theorem~\ref{BrKgen} once we have shown that
\begin{equation}
\label{eqn: Lnzero}
\Lindbladian^{n,0}(a) \in \Al(n)
\quad \text{and} \quad
\Lindbladian(a) \in
\Al(n+1),
\end{equation}
and have found a constant $M = M(d, D, K)$, independent of $n$ and
$a$, such that
\begin{equation}
\label{eqn: LMna}
\norm{\Lindbladian(a) -\Lindbladian^{n,0}(a)} \les M n.
\end{equation}
Moreover the relations~\eqref{eqn: Lnzero} follow from the facts that
$\delta(a) \in \Al(n)$, $\tau_{r,s}(a) \in \Al(n)$ for $r,s \in \Rx_n$
and, by assumption~I, $\tau_{r,s}(a) \in \Al(n+1)$ for all $r, s \in
\Rx$.

Assumption~I and the commutation relation~\eqref{comm relation} imply
that $\tau_{r,s}(a) = 0$ unless $\norm{r-s}_\infty \les D$ and $(r,s)
\in (\Rx_n \times \Rx_{n+1}) \cup (\Rx_{n+1} \times \Rx_n)$. Therefore
$\Lindbladian(a) -\Lindbladian^{n,0}(a) = \sum_{(r,s) \in \Sx_n}
\tau_{r,s}(a)$ where
\begin{align*}
  \Sx_{n}
  :=&\, \bigl\{ (r,s) \in [\Rx_n \times (\Rx_{n+1} \setminus \Rx_n)]
    \cup [(\Rx_{n+1} \setminus \Rx_n) \times \Rx_n]: \norm{r-s}_\infty
    \les D \bigr\} \\
  \subset&\, \bigl\{ (r,s): r \in \Rx_n \setminus \Rx_{n-1}, s \in (\Rx_1 +r)
    \big\} \cup \bigl\{ (r,s): s \in \Rx_n \setminus \Rx_{n-1}, r \in
    (\Rx_1 +s) \bigr\}.
\end{align*}
Now $\# (\Rx_k +t) = (1 +2kD)^d$, for $k \in \Nat$ and $t \in \Rx$, so
\begin{align*}
  \# \Sx_n
  &\les 2\Bigl( \bigl( 1 +2nD \bigr)^d -\bigl( 1 +2(n-1) D \bigr)^d
    \Bigr) (1+2D)^d \\
  &\les 4Dd \bigl( 1 +2nD \bigr)^{d-1} (1 +2D)^d \les 4 Dd (1
    +2D)^{2d-1} n^{d-1}.
\end{align*}
Moreover, for all $(r,s) \in \Sx_n$, $\norm{\tau_{r,s}(a)} \les 2
\abs{\al_{r,s}}^2$ and so, by assumption~II,
\[
\norm{\tau_{r,s}(a)} \les 2K ((n+1)D)^{2-d} = 2K D^{2-d} \bigl( 1
+\tfrac{1}{n} \bigr)^{2-d} n^{2-d}.
\]
Therefore, since $(1 +\tfrac{1}{n})^{2-d} \les 1 +\de_{1,d}$
(Kronecker delta), these estimates combine to yield~\eqref{eqn: LMna},
as required, with $M = 8Kd (1 +\de_{1,d}) D^{3-d} (1 +2D)^{2d-1}$.
\end{example}


\begin{thebibliography}{PMQ}

\bibitem[Acc]{Accardi}
L.~Accardi,
On the quantum Feynman--Kac formula,
\emph{Rend.\ Sem.\ Mat.\ Fis.\ Milano}
\textbf{48} (1978), 135--180 (1980).


\bibitem[AFL]{AFL}
L.~Accardi, A.~Frigerio and J.T.~Lewis,
Quantum stochastic processes,
\emph{Publ.\ Res.\ Inst.\ Math.\ Sci.}
\textbf{18} (1982), no. 1, 97-–133.

\bibitem[AcK]{AK}
L.~Accardi and S.V.~Kozyrev,
On the structure of Markov flows,
\emph{Chaos, Solitons \& Fractals}
\textbf{12} (2001), nos.~14--15, 2639--2655.



\bibitem[BW$_1$]{BeW}
A.C.R.~Belton and S.J.~Wills,
An algebraic construction of quantum stochastic flows with unbounded
generators,
\emph{Ann.\ Inst.\ H.\ Poincar\'{e} Probab.\ Statist.}\
\textbf{51} (2015), no.~1, 349--375.

\bibitem[BW$_2$]{BWFK}
--- --- ,
Feynman--Kac perturbation of $C^*$ quantum-stochastic flows,
\emph{Preprint}.

\bibitem[BlL]{BlL}
D.P.~Blecher and C.~Le Merdy,
``Operator Algebras and Their Modules: An Operator Space Approach,''
OUP, Oxford 2005.

\bibitem[BrK]{BrK}
O.~Bratteli and A.~Kishimoto,
Generation of semigroups, and two-dimensional quantum lattice systems,
\emph{J.\ Funct.\ Anal.}\
\textbf{35} (1980), no.\ 3, 344--368.

\bibitem[BR$_1$]{BrR}
O.~Bratteli and D.W.~Robinson,
``Operator Algebras and Quantum Statistical Mechanics I:
$C^*$- and $W^*$-Algebras, Symmetry Groups, Decompositions of States,''
Corrected 2nd Edition, Springer-Verlag, Heidelberg, 2002.

\bibitem[BR$_2$]{BrR2}
--- --- ,
``Operator Algebras and Quantum Statistical Mechanics II:
Equilibrium States. Models in Quantum Statistical Mechanics,''
Corrected 2nd Edition, Springer-Verlag, Heidelberg, 2002.


\bibitem[ChE]{ChrisEvans}
E.~Christensen and D.E.~Evans,
Cohomology of operator algebras and quantum dynamical semigroups,
\emph{J. London Math.\ Soc.}\
\textbf{20} (1979), no.~2, 358--368.

\bibitem[DGS]{DGS}
B.~Das, D.~Goswami and K.B.~Sinha,
A homomorphism theorem and a Trotter product formula
for quantum stochastic flows with unbounded coefficients,
\emph{Comm.\ Math.\ Phys.}\
\textbf{330} (2014), no.~2, 435--467.

\bibitem[DLT]{sesquiBanach}
B.K.~Das, J.M.~Lindsay and O.~Tripak,
Sesquilinear quantum stochastic analysis in Banach space,
\emph{J. Math.\ Anal.\ Appl.}\
\textbf{409} (2014), no.~2, 1032--1051.

\bibitem[Dav]{Davies}
E.B.~Davies,
``Linear Operators and their Spectra,''
Cambridge University Press, Cambridge, 2007.

\bibitem[EfR]{EfR}
E.G.~Effros and Z.-J.~Ruan,
``Operator Spaces,''
OUP, Oxford 2000.

\bibitem[Ell]{strong-norm}
G.A.~Elliott,
On the convergence of a sequence of completely positive maps to the
identity,
\emph{J. Austral.\ Math.\ Soc.\ Ser.\ A}
\textbf{68} (2000), no.~3, 340--348.


%



\bibitem[EvK]{QuTome}
D.E.~Evans and Y.~Kawahigashi,
``Quantum Symmetries on Operator Algebras,''
OUP, Oxford 1998.

\bibitem[Fag]{FagPryc}
F.~Fagnola,
Quantum Markov semigroups and quantum flows,
\emph{Proyecciones}
\textbf{18} (1999) no.~3, 144 pp.

\bibitem[GKS]{GKS}
V. Gorini, A. Kossakowski and E.C.G. Sudarshan,
Completely positive dynamical semigroups on $N$-level systems,
\emph{J. Math.\ Phys.}\
\textbf{17} (1076), no. 5, 821--825.

\bibitem[HuP]{HuP}
R.L.~Hudson and K.R.~Parthasarathy,
Quantum It\^{o}'s formula and stochastic evolutions,
\emph{Comm.\ Math.\ Phys.}\
\textbf{93} (1984), no. 3, 301-–323.

\bibitem[KuM]{KuM}
B.~K\"{u}mmerer and H.~Maassen,
The essentially commutative dilations of dynamical semigroups on
$M_n$,
\emph{Comm.\ Math.\ Phys.}\
\textbf{109} (1987), no. 1, 1--22.


\bibitem[Lin]{Lindblad}
G. Lindblad,
On the generators of quantum dynamical semigroups,
\emph{Comm.\ Math.\ Phys.}\
\textbf{48} (1976), no. 2, 119--130.

\bibitem[L]{Lgreifswald}
J.M.~Lindsay,
Quantum stochastic analysis --- an introduction,
\emph{in}
``Quantum Independent Increment Processes I,''
\emph{eds.\ M.~Sch\"{u}rmann \& U.~Franz},
Lecture Notes in Mathematics \textbf{1865},
Springer, Heidelberg 2005.

\bibitem[LiP]{gran}
J.M.~Lindsay and K.R.~Parthasarathy,
On the generators of quantum stochastic flows,
\emph{J. Funct.\ Anal.}\
\textbf{158} (1998), no.~2, 521--549.

\bibitem[LW$_1$]{mother}
J.M.~Lindsay and S.J.~Wills,
Existence, positivity, and contractivity for quantum stochastic flows
with infinite dimensional noise,
\emph{Probab.\ Theory Related Fields}
\textbf{116} (2000), no.~4, 505--543.

\bibitem[LW$_2$]{father}
--- --- ,
Markovian cocycles on operator algebras, adapted to a Fock filtration,
\emph{J. Funct.\ Anal.}\
\textbf{178} (2000), no.~2, 269--305.

%

\bibitem[LW$_3$]{LWhom}
--- --- ,
Homomorphic Feller cocycles on a $C^*$-algebra,
\emph{J.\ London Math.\ Soc.\ (2)}
\textbf{68} (2003), no.~1, 255--272.

\bibitem[LW$_4$]{LWhat}
--- --- ,
Multiplicativity by a hat trick,
\emph{in},
``Quantum Probability \& Infinite-dimensional Analysis''
QP--PQ \textbf{15}
(\emph{Ed}.\ W.~Freudenberg),
World Scientific, 2003, pp.~181--193.

\bibitem[LW$_5$]{spawn}
--- --- ,
Quantum stochastic cocycles and completely bounded semigroups on
operator spaces
\emph{Int.\ Math.\ Res.\ Not.\ IMRN}
(2014), no.~11, 3096--3139.



\bibitem[PMQ]{questions}
L.~Pantale\'{o}n-Martinez and R.~Quezada,
The asymmetric exclusion quantum Markov semigroup,
\emph{Infin.\ Dimens.\ Anal.\ Quantum Probab.\ Relat.\ Top.}\
\textbf{12} (2009), no. 3, 367--385.

\bibitem[Par]{Partha}
K.R.~Parthasarathy,
``Introduction to Quantum Stochastic Calculus'',
[2012 reprint of the 1992 original]
\emph{Modern Birkh\"{a}user Classics},
Birkh\"{a}user/Springer Basel AG, Basel, 1992.

\bibitem[Pau]{Paulsen}
V.~Paulsen,
``Completely bounded maps and operator algebras,''
CUP, Cambridge 2002.

\bibitem[Pis]{Pisier}
G.~Pisier,
``Introduction to Operator Space Theory,''
London Mathematical Society Lecture Note Series
\textbf{294},
CUP, Cambridge 2003.

\bibitem[Reb]{Reb}
R.~Rebolledo,
Decoherence of quantum Markov semigroups,
\emph{Ann.\ Inst.\ H.\ Poincar\'{e} Probab.\ Statist.}\
\textbf{41} (2005), no.~3, 349--373.

\bibitem[SSS]{SSS}
L.~Sahu, M.~Sch\"urmann and K.B.~Sinha,
Unitary processes with independent increments and representations of
Hilbert tensor algebras,
\emph{Publ.\ Res.\ Inst.\ Math.\ Sci.}\
\textbf{45} (2009), no.~3, 745-–785.

\bibitem[Ske]{Skeide}
 M.~Skeide,
Indicator functions of intervals are totalizing in the symmetric Fock
space $\Gamma(L^2(\Rplus))$,
\emph{in}
``Trends in Contemporary Infinite Dimensional Analysis and Quantum
Probability. Volume in Honour of Takeyuki Hida,''
(\emph{Eds}.\ L.~Accardi, H.-H.~Kuo, N.~Obata, K.~Saito, Si~Si and
L.~Streit), Istituto Italiano di Cultura, Kyoto, 2000.


\bibitem[W]{JTL}
S.J.~Wills,
On the generators of operator Markovian cocycles,
\emph{Markov Proc.\ Related Fields}
\textbf{13} (2007), no.~1, 191--211.
\end{thebibliography}
\end{document}